\title{Spectral Measures and Generating Series for Nimrep Graphs in Subfactor Theory}
\author{
        David E. Evans and Mathew Pugh \\ \\
        School of Mathematics, \\
        Cardiff University, \\
        Senghennydd Road, \\
        Cardiff, CF24 4AG, \\
        Wales, U.K.
}
\date{\today}

\documentclass[12pt]{article}
\usepackage{amssymb}
\usepackage{graphicx}

\newtheorem{definition}{Definition}[section]
\newtheorem{proposition}[definition]{Proposition}
\newtheorem{lemma}[definition]{Lemma}
\newtheorem{corollary}[definition]{Corollary}
\newtheorem{theorem}[definition]{Theorem}

\begin{document}
\maketitle

\begin{abstract}
We determine spectral measures for some nimrep graphs arising in subfactor theory, particularly those associated with $SU(3)$ modular invariants and subgroups of $SU(3)$. Our methods also give an alternative approach to deriving the results of Banica and Bisch for $ADE$ graphs and subgroups of $SU(2)$ and explain the connection between their results for affine $ADE$ graphs and the Kostant polynomials. We also look at the Hilbert generating series of associated pre-projective algebras.
\end{abstract}

\section{Introduction}

Banica and Bisch \cite{banica/bisch:2007} studied the spectral measures of bipartite graphs, particularly those of norm less than two, the $ADE$ graphs, and those of norm two, their affine versions associated with subgroups of $SU(2)$. Here and in a sequel \cite{evans/pugh:2009vi} we look at such spectral measures in a wider context, particularly from the viewpoint of associating spectral measures to nimreps (non-negative integer matrix representations).
$ADE$ graphs appear in the classification of non-negative integer matrices with norm less than two \cite{goodman/de_la_harpe/jones:1989}. Their affine version $A^{(1)}, D^{(1)}, E^{(1)}$ classify the finite subgroups of $SU(2)$. The $ADE$ graphs are also relevant for the classification of subfactors with Jones index less than 4, but only $A, D_{even}, E_6, E_8$ appear as principal graphs (\cite{ocneanu:1988,izumi:1991,kawahigashi:1995,bion-nadal:1991,izumi:1994} or see \cite{evans/kawahigashi:1998} and references therein). However all appear in the classification of $SU(2)$ modular invariants by Cappelli, Itzykson and Zuber \cite{cappelli/itzykson/zuber:1987ii}, and in their realisation by $SU(2)$ braided subfactors \cite{ocneanu:2000i,xu:1998,bockenhauer/evans/kawahigashi:2000}.

The Verlinde algebra of $SU(n)$ at level $k$ is represented by a non-degenerately braided system of endomorphisms ${}_N \mathcal{X}_N$ on a type $\mathrm{III}_1$ factor $N$, whose fusion rules $\{ N_{\lambda \nu}^{\mu} \}$ reproduce exactly those of the positive energy representations of the loop group of $SU(n)$ at level $k$,
$N_{\lambda} N_{\mu} = \sum_{\nu} N_{\lambda \nu}^{\mu} N_{\nu}$
and whose statistics generators $S$, $T$ obtained from the braided tensor category ${}_N \mathcal{X}_N$ match exactly those of the Ka\u{c}-Peterson modular $S$, $T$ matrices which perform the conformal character transformations \cite{wassermann:1998}.
This family $\{ N_{\lambda} \}$ of commuting normal matrices can be simultaneously diagonalised:
\begin{equation} \label{eqn:verlinde_formula}
N_{\lambda} = \sum_{\sigma} \frac{S_{\sigma, \lambda}}{S_{\sigma,1}} S_{\sigma} S_{\sigma}^{\ast},
\end{equation}
where $1$ is the trivial representation.
The intriguing aspect being that the eigenvalues $S_{\sigma, \lambda}/S_{\sigma,1}$ and eigenvectors $S_{\sigma} = \{ S_{\sigma, \mu} \}_{\mu}$ are described by the modular $S$ matrix. A braided subfactor is an inclusion $N \subset M$ where the dual canonical endomorphism decomposes as a finite combination of elements of the Verlinde algebra, endomorphisms in ${}_N \mathcal{X}_N$. Such subfactors yield modular invariants through the procedure of $\alpha$-induction
which allows two extensions of $\lambda$ on $N$, depending on the use of the braiding or its opposite, to endomorphisms
$\alpha^{\pm}_{\lambda}$ of $M$, so that the matrix $Z_{\lambda,\mu} = \langle \alpha_{\lambda}^+, \alpha_{\mu}^- \rangle$
is a modular invariant \cite{bockenhauer/evans/kawahigashi:1999,bockenhauer/evans:2000,evans:2003}.
The classification of Cappelli, Itzykson and Zuber of $SU(2)$ modular invariants is understood via the action of the $N$-$N$ sectors ${}_N \mathcal{X}_N$ on the $M$-$N$ sectors ${}_M \mathcal{X}_N$ and produces a nimrep
$G_{\lambda} G_{\mu} = \sum_{\nu} N_{\lambda \nu}^{\mu} G_{\nu}$
whose spectrum reproduces exactly the diagonal part of the modular invariant, i.e.
\begin{equation} \label{eqn:verlinde_formulaG}
G_{\lambda} = \sum_i \frac{S_{i,\lambda}}{S_{i,1}} \psi_i \psi_i^{\ast},
\end{equation}
with the spectrum of $G_{\lambda} = \{ S_{\mu, \lambda}/S_{\mu,1}$ with multiplicity $Z_{\mu,\mu} \}$ \cite{bockenhauer/evans/kawahigashi:2000}.
Every $SU(2)$ modular invariant can be realised by $\alpha$-induction for a suitable braided subfactor. Evaluating the nimrep $G$ at the fundamental representation $\rho$, we obtain for each such inclusion a matrix $G_{\rho}$ which recovers the $ADE$ classification of Cappelli, Itzykson and Zuber.
Since these $ADE$ graphs can be matched to the affine Dynkin diagrams, the McKay graphs of the finite subgroups of $SU(2)$, di Francesco and Zuber \cite{di_francesco/zuber:1990} were guided to find candidates for classifying graphs for $SU(3)$ modular invariants by first considering the McKay graphs of the finite subgroups of $SU(3)$ to produce a candidate list of graphs whose spectra described the diagonal part of the modular invariant. Ocneanu claimed \cite{ocneanu:2002} that all $SU(3)$ modular invariants were realised by subfactors and this was shown in \cite{evans/pugh:2009ii}.
The nimrep associated to the conjugate Moore-Seiberg modular invariant $Z_{\mathcal{E}_{MS}^{(12)}}$ was not computed however in \cite{evans/pugh:2009ii}.
To summarize, we can  realize all $SU(3)$ modular invariants, but there is mismatch between the list of nimreps associated to each modular invariant and the McKay graphs of the finite subgroups of $SU(3)$ which are also the nimreps of the representation theory of the group. Both of these kinds of nimreps will play a role in this paper and its sequel \cite{evans/pugh:2009vi}. These nimreps also have a diagonalisation as in (\ref{eqn:verlinde_formula}) with diagonalising matrix $S = \{ S_{ij} \}$ usually non-symmetric, where $i$ labels conjugacy classes, and $j$ the irreducible characters (see \cite[Section 8.7]{evans/kawahigashi:1998} and Section \ref{Sect:Spectral_Measures_for_subgroup_SU(2)}).

We compute here the spectral measures of nimreps of braided subfactors associated to $SU(2)$ and $SU(3)$ and nimreps for the representations of subgroups of $SU(2)$. The case of subgroups of $SU(3)$ will be treated separately \cite{evans/pugh:2009vi}.
Suppose $A$ is a unital $C^{\ast}$-algebra with state $\varphi$.
If $b \in A$ is a normal operator then there exists a compactly supported probability measure $\mu_b$ on the spectrum $\sigma(b) \subset \mathbb{C}$ of $b$, uniquely determined by its moments
\begin{equation} \label{eqn:moments_normal_operator}
\varphi(b^m b^{\ast n}) = \int_{\sigma(b)} z^m \overline{z}^n \mathrm{d}\mu_b (z),
\end{equation}
for non-negative integers $m$, $n$.
If $a$ is self-adjoint (\ref{eqn:moments_normal_operator}) reduces to
\begin{equation} \label{eqn:moments_selfadjoint}
\varphi(a^m) = \int_{\sigma(a)} x^m \mathrm{d}\mu_a (x),
\end{equation}
with $\sigma(a) \subset \mathbb{R}$, for any non-negative integer $m$. The generating series of the moments of $a$ is the Stieltjes transform $\sigma(z)$ of $\mu_a$, given by
\begin{equation} \label{Stieltjes_transform}
\sigma(z) = \sum_{m=0}^{\infty} \varphi(a^m)z^m = \sum_{m=0}^{\infty} \int_{\sigma(a)} x^m z^m \mathrm{d}\mu_a (x) = \int_{\sigma(a)} \frac{1}{1-x z} \mathrm{d}\mu_a (x).
\end{equation}
What we compute are such spectral measures and generating series when $b$ is the normal operator $\Delta = G_{\rho}$ acting on the Hilbert space of square summable functions on the graph.

In particular we can understand the spectral measures for the torus $\mathbb{T}$ and $SU(2)$ as follows.
If $w_Z$ and $w_N$ are the self adjoint operators arising from the McKay graph of the fusion rules of the representation theory of $\mathbb{T}$ and $SU(2)$, then the spectral measures in the vacuum state can be describe  in terms of semicircular law, on the interval $[-2,2]$ which is the spectrum of either
as the image of the map $z \in \mathbb{T} \rightarrow z + z^{-1}$:
$$\mathrm{dim}\left( \left(\otimes^k M_2 \right)^{\mathbb{T}} \right) \;\; = \;\; C^{2k}_k \;\; = \;\; \varphi(w_Z^{2k}) \;\; = \;\; \frac{1}{\pi} \int_{-2}^2 x^{2k} \frac{1}{\sqrt{4-x^2}} \; \mathrm{d}x \, ,$$
$$\mathrm{dim}\left( \left(\otimes^k M_2 \right)^{SU(2)} \right) \;\; = \;\; \frac{1}{k+1}C^{2k}_k \;\; = \;\; \varphi(w_N^{2k}) \;\; = \;\; \frac{1}{2\pi} \int_{-2}^2 x^{2k} \sqrt{4-x^2} \; \mathrm{d}x \, ,$$
where $C^{r}_s$ and $C^{2k}_k/(k+1)$ denote Binomial coefficients and Catalan numbers respectively.
The spectral weight for $SU(2)$ arises from the Jacobian of a change of variable between the interval $[-2,2]$ and the circle.
Then for $\mathbb{T}^2$ and $SU(3)$,  the deltoid $\mathfrak{D}$ in the complex plane which is the image of the two-torus under the map $(\omega_1, \omega_2) \rightarrow \omega_1 + \omega_2^{-1} + \omega_1^{-1} \omega_2$ is the spectrum of the corresponding normal operators on the Hilbert spaces of the fusion graphs. The corresponding spectral measures are then described by a corresponding Jacobian or discriminant as:
\begin{eqnarray*}
\mathrm{dim}\left( \left(\otimes^k M_3 \right)^{\mathbb{T}^2} \right) & = & \sum_{j=0}^k C^{2j}_j (C^k_j)^2 \;\; = \;\; \varphi(|v_Z|^{2k}) \\
& = & \frac{3}{\pi^2} \int_{\mathfrak{D}} |z|^{2k} \frac{1}{\sqrt{27 - 18z\overline{z} + 4z^3 + 4\overline{z}^3 -z^2 \overline{z}^2}} \; \mathrm{d}z \, ,
\end{eqnarray*}
$$\mathrm{dim}\left( \left(\otimes^k M_3 \right)^{SU(3)} \right) \;\; = \;\; \varphi(|v_N|^{2k}) \;\; = \;\; \frac{1}{2\pi^2} \int_{\mathfrak{D}} |z|^{2k} \sqrt{27 - 18z\overline{z} + 4z^3 + 4\overline{z}^3 -z^2 \overline{z}^2} \; \mathrm{d}z \, ,$$
where $\mathrm{d}z:=\mathrm{d}\,\mathrm{Re}z \; \mathrm{d}\,\mathrm{Im}z$ denotes the Lebesgue measure on $\mathbb{C}$.
Then for the other graphs, the quantum graphs, the spectral measures distill onto very special subsets of the semicircle/circle, torus/deltoid and the theory of nimreps allows us to compute these measures precisely. For the case of finite subgroups, this nimrep approach clearly shows why Banica and Bisch were recovering the Kostant polynomials for finite subgroups of $SU(2)$.

We are also going to compute various Hilbert series of dimensions associated to $ADE$ models. In the $SU(2)$ case this corresponds to the study of the McKay correspondence \cite{reid:2002}, Kostant polynomials of \cite{kostant:1984}, the $T$-series of \cite{banica/bisch:2007}, and the study of pre-projective algebras \cite{brenner/butler/king:2002,malkin/ostrik/vybornov:2006}. The classical McKay correspondence relates finite subgroups $\Gamma$ of $SU(2)$ with the algebraic geometry of the quotient Kleinian singularities $\mathbb{C}^2/\Gamma$ but also with the classification of $SU(2)$ modular invariants, classification of subfactors of index less than 4, and quantum subgroups of $SU(2)$. The corresponding $SU(3)$ theory is related to the AdS-CFT correspondence and the Calabi-Yau algebras arising in the geometry of Calabi-Yau manifolds.

We take the superpotentials built on the $ADE$ Perron-Frobenius weights and the $\mathcal{ADE}$ cells \cite{ocneanu:2000ii, evans/pugh:2009i} and corresponding associated algebraic structures and study the Hilbert series of dimensions of corresponding algebras.
If $H_n$ is the matrix of dimensions of paths of length $n$ in a graph $\mathcal{G}$ in the pre-projective algebra $\Pi$, with indices labeled by the vertices, then the matrix Hilbert series $H$ of the pre-projective algebra is defined as $H(t) = \sum H_n t^n$. Let $\Delta$ be the adjacency matrix of $\mathcal{G}$.
Then if $\mathcal{G}$ is a finite (unoriented) graph which is not an $ADET$ graph (where $T$ denotes the tadpole graph $\mathrm{Tad}_n$), then $H(t) = (1 - \Delta t + t^2)^{-1}$, whilst if $\mathcal{G}$ is an $ADET$ graph, then $H(t) = (1 + Pt^h)(1 - \Delta t + t^2)^{-1}$, where $h$ is the Coxeter number of $\mathcal{G}$ and $P$ is the permutation matrix corresponding to an involution of the vertices of $\mathcal{G}$ \cite{malkin/ostrik/vybornov:2006}.

The dual $\Pi^{\ast} = \mathrm{Hom}(\Pi, \mathbb{C})$ is a $\Pi$-$\Pi$ bimodule, not usually identified with $_{\Pi}\Pi_{\Pi}$ or $_1\Pi_1$ with trivial right and left actions but with $_1\Pi_\nu$ with trivial left action and the right action twisted by an automorphism, the Nakayama automorphism $\nu$. The Nakayama automorphism measures how far away $\Pi$ is from being symmetric. In the case of a pre-projective algebra of Dynkin quiver, this Nakayama automorphism is identified with an involution on the underlying Dynkin diagram. More precisely it is trivial in all cases, except for the Dynkin diagrams $A_n$, $D_{2n+1}$, $E_6$ where it is the unique non-trivial involution. Bocklandt \cite{bocklandt:2008} has studied the types of quivers and relations (superpotentials) that appear in graded Calabi-Yau algebras of dimension 3. Indeed he also points out that the zero-dimensional case consists of semi-simple algebras, i.e. quivers without arrows, the one dimensional case consists of direct sums of one-vertex-one-loop quivers. Moreover, a Calabi-Yau algebra of dimension 2 is the pre-projective algebra of a non-Dynkin quiver.
The examples coming from finite subgroups of $SU(3)$ give Calabi-Yau algebras of dimension three \cite[Theorem 4.4.6]{ginzburg:2006}.

We are not only interested in the fusion graphs of finite subgroups of $SU(3)$, whose adjacency matrices have norm 3, but in the fusion $\mathcal{ADE}$ nimrep graphs arising in our subfactor setting as describing the $SU(3)$ modular invariants through $M$-$N$ systems which have norm less than 3.
The figures for the complete list of the $\mathcal{ADE}$ graphs are given in \cite{behrend/pearce/petkova/zuber:2000,evans/pugh:2009i}.
Unlike for $SU(2)$, there is no precise relation between finite subgroups of $SU(3)$ and $SU(3)$-modular invariants. For $SU(2)$ an affine Dynkin diagram describing the McKay graph of a finite subgroup gives rise to a Dynkin diagram describing a nimrep of a modular invariant by removing one vertex and the edges which have this vertex as an endpoint. For $SU(3)$, di Francesco and Zuber \cite{di_francesco/zuber:1990} used this procedure as a guide to find nimreps for some $SU(3)$-modular invariants by removing vertices from some McKay graphs of finite subgroups of $SU(3)$. However, not all finite subgroups were utilised, and not all nimreps or modular invariants can be found from a finite subgroup.

The spectral measures for the $ADE$ graphs were computed in terms of probability measures on the circle $\mathbb{T}$ in \cite{banica/bisch:2007}. In Section \ref{sect:spec_measure-ADE} we recover their results via a different method, which depends on the fact that the $ADE$ graphs are nimrep graphs. This method can then be generalized to $SU(3)$, which we do in Section \ref{sect:spec_measureSU(3)ADE}, and in particular obtain spectral measures for the infinite graphs $\mathcal{A}^{(6\infty)}$ and $\mathcal{A}^{(\infty)}$ corresponding to the representation graphs of the fixed point algebra of $\bigotimes_{\mathbb{N}} M_3$ under the action of $\mathbb{T}^2$ and $SU(3)$ respectively. We also obtain the spectral measure for the finite graphs $\mathcal{A}^{(n)}$, $\mathcal{A}^{(n)\ast}$, $n \geq 4$, and $\mathcal{D}^{(3k)}$, $k \geq 2$.
Finally, in Section \ref{sect:Hilbert-SU(3)} we consider the Hilbert series of the dimensions of the associated pre-projective algebras.

The final section depends on the existence of the cells \cite{ocneanu:2000ii,ocneanu:2002} (essentially the square roots of the Boltzmann weights) and to some degree on their explicit values computed in \cite{evans/pugh:2009i}. The theory of modular invariants constructed from braided subfactors \cite{bockenhauer/evans:1998,bockenhauer/evans:2000,bockenhauer/evans/kawahigashi:1999,bockenhauer/evans/kawahigashi:2000} also provides us with nimreps associated to $SU(3)$ modular invariants. It was announced by Ocneanu \cite{ocneanu:1988} and shown in \cite{evans/pugh:2009ii} that every $SU(3)$ modular invariant is realised by a braided subfactor.

\section{$SU(2)$ Case}

In this section we will compute the spectral measures for the $ADE$ Dynkin diagrams and their affine counterparts. We will present a method for computing these spectral measures using the fact that the graphs are nimrep graphs. This method recovers the measures given in \cite{banica/bisch:2007} and will allow for an easy generalization to the case of $SU(3)$ and associated nimrep graphs.

A graph is called locally finite if each vertex is the start or endpoint for a finite number of edges. Let $\mathcal{G}$ be any locally finite bipartite graph, with a distinguished vertex labelled $\ast$ and bounded adjacency matrix $\Delta$ regarded as an operator on $\ell^2(\mathcal{G}^{(0)})$, where $\mathcal{G}^{(0)}$ denotes the vertices of $\mathcal{G}$. Let $A(\mathcal{G})_k$ be the algebra generated by pairs $(\eta_1,\eta_2)$ of paths from the distinguished vertex $\ast$ such that $r(\eta_1)=r(\eta_2)$ and $|\eta_1|=|\eta_2|=k$. Then $A(\mathcal{G}) = \overline{\bigcup_k A(\mathcal{G})_k}$ is called the path algebra of $\mathcal{G}$ (see \cite{evans/kawahigashi:1998} for more details).
Let $\varphi$ be a state on $C^{\ast}(\Delta)$. From (\ref{eqn:moments_selfadjoint}), we define the spectral measure of $\mathcal{G}$ to be the probability measure $\mu_{\Delta}$ on $\mathbb{R}$ given by $\int_{\mathbb{R}} \psi(x) \mathrm{d}\mu_{\Delta}(x) = \varphi(\psi(\Delta))$, for any continuous function $\psi:\mathbb{R} \rightarrow \mathbb{C}$, as in \cite{banica/bisch:2007}.

\subsection{Spectral measure for $A_{\infty,\infty}$}

We begin by looking at the fixed point algebra of $\bigotimes_{\mathbb{N}} M_2$ under the action of the group $\mathbb{T}$. Let $\rho$ be the fundamental representation of $SU(2)$, so that its restriction to $\mathbb{T}$ is given by
\begin{equation}  \label{eqn:restrict_rho_to_T}
(\rho|_{\mathbb{T}})(t) = \mathrm{diag}(t,\overline{t}),
\end{equation}
where $t \in \mathbb{T}$.

Let $\{ \chi_i \}_{i \in \mathbb{N}}$, $\{ \sigma_i \}_{i \in \mathbb{Z}}$ be the irreducible characters of $SU(2)$, $\mathbb{T}$ respectively, where $\chi_0$ is the trivial character of $SU(2)$, $\chi_1$ is the character of $\rho$, and $\sigma_i (z) = z^i$, $i \in \mathbb{Z}$. If $\sigma$ is the restriction of $\chi_1$ to $\mathbb{T}$, we have $\sigma = \sigma_1 + \sigma_{-1}$ (by (\ref{eqn:restrict_rho_to_T})), and $\sigma \sigma_i = \sigma_{i-1} + \sigma_{i+1}$, for any $i \in \mathbb{Z}$. Then the representation graph of $\mathbb{T}$ is identified with the doubly infinite graph $A_{\infty,\infty}$, illustrated in Figure \ref{fig:A_infty,infty}, whose vertices are labelled by the integers $\mathbb{Z}$ which correspond to the irreducible representations of $\mathbb{T}$, where we choose the distinguished vertex to be $\ast = 0$. The Bratteli diagram for the path algebra of the graph $A_{\infty,\infty}$ with initial vertex $\ast$ is given by Pascal's triangle. The dimension of the $0^{\mathrm{th}}$ level of the path algebra is 1, and we compute the dimensions of the matrix algebras corresponding to minimal central projections at the other levels according to the rule that for a vertex $(v,n)$ at level $n$ we take the sum of the dimensions at level $n-1$ corresponding to vertices $(v',n-1)$ for which there is an edge in the Bratteli diagram from $(v',n-1)$ to $(v,n)$. It is well-known that these numbers give the binomial coefficients, with the $j^{\mathrm{th}}$ vertex along level $m$ giving $C^m_j$, and we see that $\sigma^m = \sum_{j=0}^m C^m_j \sigma_{m-2j}$, where $C^m_j$ are the binomial coefficients.

\begin{figure}[tb]
\begin{center}
  \includegraphics[width=75mm]{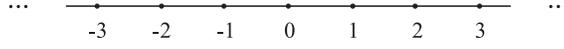}\\
 \caption{Doubly infinite graph $A_{\infty,\infty}$.} \label{fig:A_infty,infty}
\end{center}
\end{figure}

Recall that if $\{ \pi_i \}$ denote irreducible representations of a group $G$, and if $\pi = n_1 \pi_1 \oplus n_2 \pi_2 \oplus \cdots \;$ on a full matrix algebra, then the fixed point algebra of the action $\mathrm{Ad}(\pi)$ is isomorphic to $M = M_{n_1} \oplus M_{n_2} \oplus \cdots \;$, and the dimension of $M$ is given by the sum of the squares of the $n_i$. Then we see that $(\otimes^k M_2)^{\mathbb{T}} \cong \bigoplus_{j=0}^k M_{C^k_j}$, and $(\bigotimes_{\mathbb{N}} M_2)^{\mathbb{T}} \cong A(A_{\infty,\infty})$.
Hence $\mathrm{dim}(A(A_{\infty,\infty})_k) = \mathrm{dim}\left( \left(\otimes^k M_2 \right)^{\mathbb{T}} \right) = \sum_{j=0}^k (C^k_j)^2 = C^{2k}_k$. Counting the number $p_j$ of pairs of paths in $A(A_{\infty,\infty})_k$ which end at a vertex $k-2j$ of $A_{\infty,\infty}$ is the same as the dimension of the subalgebra of $A(A_{\infty,\infty})_k$ which corresponds to the vertex $k-2j$ at level $k$ of the Bratteli diagram for $A(A_{\infty,\infty})$, and hence $p_j$ is given by the binomial coefficient $p_j = C^k_j$.

We define an operator $w_Z$ on $\ell^2(\mathbb{Z})$ by $w_Z = s + s^{-1}$, where $s$ is the bilateral shift on $\ell^2(\mathbb{Z})$. Let $\Omega$ be the vector $(\delta_{i,0})_i$.
Then $w_Z$ is identified with the adjacency matrix $\Delta_{\infty, \infty}$ of $A_{\infty, \infty}$, where we regard the vector $\Omega$ as corresponding to the vertex $0$ of $A_{\infty, \infty}$, and the shifts $s$, $s^{-1}$ correspond to moving along an edge to the right, left respectively on $A_{\infty, \infty}$. Then $s^k \Omega$ corresponds to the vertex $k$ of $A_{\infty}$, $k \in \mathbb{Z}$, the identity $s^{-1} s = s s^{-1} = 1$ correspond to moving along an edge of $A_{\infty, \infty}$ and then back along the reverse edge, arriving back at the original vertex we started at. Applying $w_Z^n$, $n \geq 0$, to $\Omega$ gives a vector $v=(v_i)_{i \in \mathbb{Z}}$ in $\ell^2(A_{\infty, \infty})$, where $v_i$ gives the number of paths of length $n$ from the vertex 0 to the vertex $i$ of $A_{\infty, \infty}$.

The binomial coefficient $C^{2n}_n$ counts the number of `balanced' paths of length $2n$ on the integer lattice $\mathbb{Z}^2$ \cite{king/egecioglu:1999}, that is, paths of length $2n$ starting from the point $(0,0)$ and ending at the point $(2n,0)$ where each edge is a vector equal to a translation of the vectors $(0,0) \rightarrow (1,1)$ or $(0,0) \rightarrow (1,-1)$.

We define a state $\varphi$ on $C^{\ast}(w_Z)$ by $\varphi( \, \cdot \, ) = \langle \, \cdot \, \Omega, \Omega \rangle$.
The odd moments are all zero.
For the even moments we have
$$\varphi(w_Z^{2k}) = \varphi((s+s^{-1})^{2k}) = \sum_{j=0}^{2k} C^{2k}_j \varphi(s^{2k-2j}) = \sum_{j=0}^{2k} C^{2k}_j \delta_{j,k} = C^{2k}_k.$$

Suppose the operator $\Delta$ has norm $\leq 2$, so that the support of the spectral measure $\mu$ of $\Delta$ is contained in $[-2,2]$. There is a map $\Phi:\mathbb{T} \rightarrow [-2,2]$ given by
\begin{equation} \label{def:phi-SU(2)}
\Phi(u) = u + u^{-1},
\end{equation}
for $u \in \mathbb{T}$. Then any probability measure $\varepsilon$ on $\mathbb{T}$ produces a probability measure $\mu$ on $[-2,2]$ by
$$\int_{-2}^2 \psi(x) \mathrm{d}\mu(x) = \int_{\mathbb{T}} \psi(u + u^{-1}) \mathrm{d}\varepsilon(u),$$
for any continuous function $\psi:[-2,2] \rightarrow \mathbb{C}$.

The operator $\Delta_{\infty,\infty}$ has norm 2. Consider the measure $\varepsilon(u)$ given by $\mathrm{d}\varepsilon(u) = \mathrm{d}u$, where $\mathrm{d}u$ is the uniform Lebesgue measure on $\mathbb{T}$. Now $\int_{\mathbb{T}} u^m \mathrm{d}u = \delta_{m,0}$, hence $\int_{\mathbb{T}} (u + u^{-1})^m \mathrm{d}u = 0$ for $m$ odd, and
$$\int_{\mathbb{T}} (u + u^{-1})^{2k} \mathrm{d}u = \sum_{j=0}^{2k} C^{2k}_j \int_{\mathbb{T}} u^{2k-2j} \mathrm{d}u = C^{2k}_k = \varphi(w_Z^{2k}),$$
for $k \geq 0$ \cite[Theorem 2.2]{banica/bisch:2007}.
Now, we can write
$$\int_{\mathbb{T}} (u + u^{-1})^m \mathrm{d}u = \int_0^1 (e^{2 \pi i \theta} + e^{-2 \pi i \theta})^m \mathrm{d}\theta = 2 \int_0^{1/2} (e^{2 \pi i \theta} + e^{-2 \pi i \theta})^m \mathrm{d}\theta.$$
If we let $x = e^{2 \pi i \theta} + e^{-2 \pi i \theta} = 2\cos(2 \pi \theta)$, then $\mathrm{d}x/\mathrm{d}\theta = 2 \pi i (e^{2 \pi i \theta} - e^{-2 \pi i \theta}) = -4 \pi \sin(2 \pi \theta) = -2 \pi \sqrt{4 - x^2}$. Here the square root is always taken to be positive, since $\sin(2 \pi \theta) \geq 0$ in the range $0 \leq \theta \leq 1/2$. So
$$\int_{\mathbb{T}} (u + u^{-1})^m \mathrm{d}u = 2 \int_0^{1/2} (e^{2 \pi i \theta} + e^{-2 \pi i \theta})^n \mathrm{d}\theta = \frac{1}{\pi} \int_{-2}^2 x^m \frac{1}{\sqrt{4-x^2}} \; \mathrm{d}x.$$
Thus the spectral measure $\mu_{w_Z}$ of $w_Z$ (over $[-2,2]$) is given by $\mathrm{d}\mu_{w_Z}(x) = (\pi\sqrt{4-x^2})^{-1} \mathrm{d}x$.

Summarizing, we have the identifications
$$\mathrm{dim}(A(A_{\infty,\infty})_k) = \mathrm{dim}\left( \left(\otimes^k M_2 \right)^{\mathbb{T}} \right) = C^{2k}_k = \varphi(w_Z^{2k}) = \frac{1}{\pi} \int_{-2}^2 x^{2k} \frac{1}{\sqrt{4-x^2}} \; \mathrm{d}x.$$

\subsection{Spectral measure for $A_{\infty}$}

\begin{figure}[tb]
\begin{center}
  \includegraphics[width=130mm]{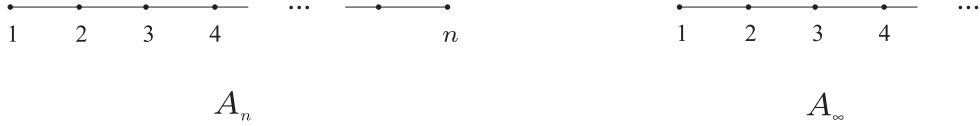}\\
 \caption{Dynkin diagrams $A_n$, $n=2,3,\ldots \;$, and $A_{\infty}$}\label{Fig:DynkinDiagramAn}
\end{center}
\end{figure}

We will now consider the fixed point algebra of $\bigotimes_{\mathbb{N}} M_2$ under the action of $SU(2)$. We have $\chi_1 \chi_i = \chi_{i-1} + \chi_{i+1}$, for $i = 0,1,2, \ldots \;$, where $\chi_{-1} = 0$. Then the representation graph of $SU(2)$ is identified with the infinite Dynkin diagram $A_{\infty}$ of Figure \ref{Fig:DynkinDiagramAn}, with distinguished vertex $\ast = 1$. Then $(\bigotimes_{\mathbb{N}} M_2)^{SU(2)} \cong A(A_{\infty})$.

We define an operator $w_N$ on $\ell^2(\mathbb{N})$ by $w_N = l + l^{\ast}$, where $l$ is the unilateral shift to the right on $\ell^2(\mathbb{N})$, and $\Omega$ by the vector $(\delta_{i,1})_i$. Then $w_N$ is identified with the adjacency matrix $\Delta_{\infty}$ of $A_{\infty}$, where we regard the vector $\Omega$ as corresponding to the vertex $\ast = 1$ of $A_{\infty}$, the creation operator $l$ as an edge to the right on $A_{\infty}$ and the annihilation operator $l^{\ast}$ as an edge to the left. For the graph $A_{\infty}$, $w_N^n \Omega =(v_i)_{i \in \mathbb{N}}$ in $\ell^2(A_{\infty})$, where $v_i$ gives the number of paths of length $n$ from the vertex 1 to the vertex $i$ of $A_{\infty}$.

Let $c_n$ be the $n^{\mathrm{th}}$ Catalan number which counts the number of Catalan (or Dyck) paths of length $2n$ in the sublattice $L$ of $\mathbb{Z}^2$ given by all points with non-negative co-ordinates. A Catalan path begins at the point $(0,0)$ and must end at the point $(2n,0)$, and is constructed from edges which are translations of the vectors $(0,0) \rightarrow (1,1)$ or $(0,0) \rightarrow (1,-1)$. The Catalan numbers $c_k$ are given explicitly by $c_k = C^{2k}_k/(k+1)$.

We define a state $\varphi$ on $C^{\ast}(w_N)$ by $\varphi( \, \cdot \, ) = \langle \, \cdot \, \Omega, \Omega \rangle$.
Once again, the odd moments are all zero. For the even moments we have $\varphi(w_N^{2k}) = c_k$, since the sequences in $l$, $l^{\ast}$ which contribute to the calculation of $\varphi(w_N^{2k})$ can be identified with the Catalan paths of length $2k$.
By \cite[Aside 5.1.1]{jones:1983}, the dimension of the $k^{\mathrm{th}}$ level of the path algebra for the infinite graph $A_{\infty}$ is given by $\mathrm{dim}(A(A_{\infty})_k) = c_k$. A connection with Catalan paths was also shown in \cite[Aside 4.1.4]{jones:1983}, since any ordered reduced word in the Temperley-Lieb algebra $\mathrm{alg}(1,e_1,\ldots,e_{k-1})$ is of the form
$$(e_{j_1} e_{j_1-1} \cdots e_{l_1})(e_{j_2} e_{j_2-1} \cdots e_{l_2}) \cdots (e_{j_p} e_{j_p-1} \cdots e_{l_p}),$$
where $j_p$ is the maximum index, $j_i \geq l_i$, $i=1,\ldots, p \;$, and $j_{i+1} > j_i$, $l_{i+1} > l_i$, $i=1,\ldots, p-1 \;$. In the generic case, when the Temperley-Lieb parameter $\delta \geq 2$, these words are linearly independent. Such an ordered reduced word corresponds to an increasing path on the integer lattice from $(0,0)$ to $(k,k)$ which does not go below the diagonal. Rotating any such path on the lattice by $\pi/4$, we obtain a path of length $2k$ corresponding to a Catalan path.
For $\delta < 2$, the ordered reduced words are linearly dependent, and we only have $\mathrm{dim}(A(A_{\infty})_k) \leq c_k$.

A self-adjoint bounded operator $a$ is called a semi-circular element with mean $\kappa \in \mathbb{R}$ and variance $r^2/4$ if its moments equal those of the semi-circular distribution centered at $\kappa$ and of radius $r > 0$, i.e. $a$ has the probability measure $\mu_a$ on $[\kappa-r,\kappa+r]$ given by
\begin{equation} \label{eqn:semicircular_measure}
\mathrm{d}\mu_a(t) = \frac{2}{\pi r^2} \sqrt{r^2 - (x-\kappa)^2} \mathrm{d}x.
\end{equation}
When $\kappa = 0$, $r=2$, this is equivalent to $a$ being an even variable with even moments given by the Catalan numbers:
$$\varphi(a^m) = \left\{ \begin{array}{ccl} c_k, & & \textrm{if } m=2k, \\ 0, & & \textrm{if } m \textrm{ odd,} \end{array} \right.$$
Thus the operator $w_N$ above is a semi-circular element.
We will reproduce a proof that the probability measure $\mu_{w_N}$ on $[-2,2]$ is given by $\mathrm{d}\mu_{w_N}(x) = (2 \pi)^{-1} \sqrt{4 - x^2} \mathrm{d}x$ in the next section. This is the spectral measure for $A_{\infty}$ given in \cite{dykema/voiculescu/nica:1992}.

Summarizing, we have the identifications
\begin{eqnarray*}
\mathrm{dim}(A(A_{\infty})_k) & = & \mathrm{dim}\left( \left(\otimes^k M_2 \right)^{SU(2)} \right) \;\; = \;\; c_k \;\; = \;\; \frac{1}{k+1}C^{2k}_k \\
& = & \varphi(w_N^{2k}) \;\; = \;\; \frac{1}{2\pi} \int_{-2}^2 x^{2k} \sqrt{4-x^2} \; \mathrm{d}x.
\end{eqnarray*}

\section{Spectral measures for the $ADE$ Dynkin diagrams via nimreps} \label{sect:spec_measure-ADE}

Let $\Delta_{\mathcal{G}}$ be the adjacency matrix of the finite (possibly affine) Dynkin diagram $\mathcal{G}$ with $s$ vertices. The $m^{\mathrm{th}}$ moment $\int x^m \mathrm{d}\mu(x)$ is given by $\langle \Delta_{\mathcal{G}}^m e_1, e_1 \rangle$, where $e_1$ is the basis vector in $\ell^2(\mathcal{G})$ corresponding to the distinguished vertex $\ast$ of $\mathcal{G}$. Note that we can in fact define many spectral measures for $\mathcal{G}$ by $\langle \Delta_{\mathcal{G}}^m e_j, e_j \rangle$, where the basis vector $e_j$ in $\ell^2(\mathcal{G})$ now corresponds to any fixed vertex $j$ of $\mathcal{G}$.

Let $\beta^j$ be the eigenvalues of $\mathcal{G}$, with corresponding eigenvectors $x^j$, $j=1,\ldots,s$. Now $\Delta_{\mathcal{G}} = \mathcal{U} \Lambda_{\mathcal{G}} \mathcal{U}^{\ast}$, where $\Lambda_{\mathcal{G}} = \mathrm{diag}(\beta^1, \beta^2, \ldots, \beta^s)$ is a diagonal matrix and $\mathcal{U} = (x^1, x^2, \ldots, x^s)$. Then $\Delta_{\mathcal{G}}^m = \mathcal{U} \Lambda_{\mathcal{G}}^m \mathcal{U}^{\ast}$, so that
\begin{eqnarray}
\int_{\mathbb{T}} \psi(u + u^{-1}) \mathrm{d}\varepsilon(u) & = & \langle \mathcal{U} \Lambda_{\mathcal{G}}^m \mathcal{U}^{\ast} e_1, e_1 \rangle \;\; = \;\; \langle \Lambda_{\mathcal{G}}^m \mathcal{U}^{\ast} e_1, \mathcal{U}^{\ast} e_1 \rangle \nonumber \\
& = & \sum_{j=1}^s (\beta^j)^m |y_i|^2, \label{eqn:moments_general_SU(2)graph}
\end{eqnarray}
where $y_i = x^i_1$ is the first entry of the eigenvector $x^i$.

For a Dynkin diagram $\mathcal{G}$ with Coxeter number $h$, its eigenvalues $\lambda^j$ are given by
\begin{equation} \label{eqn:evalues_Dynkin}
\lambda^j = 2 \cos (\pi m_j/h),
\end{equation}
with corresponding eigenvectors $(\psi^{m_j}_a)_{a \in \mathfrak{V}(\mathcal{G})}$, for the exponents $m_j$ of $\mathcal{G}$, $j=1,\ldots,s$. Then by (\ref{eqn:verlinde_formulaG}), equation (\ref{eqn:moments_general_SU(2)graph}) becomes
\begin{equation} \label{eqn:moments_Dynkin}
\int_{\mathbb{T}} \psi(u + u^{-1}) \mathrm{d}\varepsilon(u) = \sum_{j=1}^s (\lambda^j)^m |\psi^{m_j}_{\ast}|^2,
\end{equation}
where $\ast$ is the distinguished vertex of $\mathcal{G}$ with lowest Perron-Frobenius weight. Using (\ref{eqn:moments_Dynkin}) we can obtain the results for the spectral measures of the Dynkin diagrams given in \cite{banica/bisch:2007}. The advantage of this method is that it can be extended to the case of $SU(3)$ $\mathcal{ADE}$ graphs, which we will do in Section \ref{sect:spec_measureSU(3)ADE}, and also to subgroups of $SU(3)$, which we will do in the sequel \cite{evans/pugh:2009vi}.

\subsection{Dynkin diagrams $A_n$, $A_{\infty}$}

The eigenvalues $\lambda^j_n$ of $A_n$ are given by (\ref{eqn:evalues_Dynkin}) with corresponding eigenvectors $\psi^j_a = S_{a,j} = \sqrt{2/(n+1)} \; \sin(j a \pi /(n+1))$, where the exponents are $m_j = 1,2,\ldots,n$. The distinguished vertex $\ast$ of $A_n$ is the vertex 1 in Figure \ref{Fig:DynkinDiagramAn}. With $\widetilde{u} = e^{\pi i/(n+1)}$, we have $2\cos(j \pi /(n+1)) = \widetilde{u}^j + \widetilde{u}^{-j}$ and $\sin(j \pi /(n+1)) = \mathrm{Im}(\widetilde{u}^j)$. Note that $\mathrm{Im}(\widetilde{u}^j) = 0$ for $j=0,n+1$. Then
\begin{eqnarray}
\int_{\mathbb{T}} \psi(u + u^{-1}) \mathrm{d}\varepsilon(u) & = & \frac{2}{n+1} \sum_{j=1}^n \left(2\cos\left(\frac{j \pi}{n+1}\right)\right)^m \sin^2\left(\frac{j \pi}{n+1}\right) \label{eqn:momentsAn1} \\
& = & \frac{2}{n+1} \sum_{j=1}^n (\widetilde{u}^j + \widetilde{u}^{-j})^m \; \mathrm{Im}(\widetilde{u}^j)^2 \nonumber \\
& = & \frac{2}{2(n+1)} \sum_{j=0}^{2(n+1)} (\widetilde{u}^j + \widetilde{u}^{-j})^m \; \mathrm{Im}(\widetilde{u}^j)^2 \nonumber \\
& = & 2 \int_{\mathbb{T}} (u + u^{-1})^m \; \mathrm{Im}(u)^2 \, \mathrm{d}_{n+1}u \label{eqn:momentsAn2}
\end{eqnarray}
where $\mathrm{d}_{n+1}$ is the uniform measure on the $2(n+1)^{\mathrm{th}}$ roots of unity. Thus the spectral measure (over $\mathbb{T}$) for $A_n$ is $d\varepsilon(u) = 2 \mathrm{Im}(u)^2 \, \mathrm{d}_{n+1}u$.
This is the result given in \cite[Theorem 3.1]{banica/bisch:2007}

We again consider the infinite graph $A_{\infty}$, and note that the computation of the $m^{\mathrm{th}}$ moment is a finite problem, $\int x^m \mathrm{d}\mu_{w_N}(x) = \langle \Delta_{A_n}^m e_1, e_1 \rangle$, for $m < 2n$. Taking the limit in (\ref{eqn:momentsAn1}) as $n \rightarrow \infty$ (cf. the second proof of Theorem 1.1.5 in \cite{hiai/petz:2000}), we obtain a sum which is the approximation of an integral,
$$\int x^m \mathrm{d}\mu_{w_N}(x) = \frac{2}{\pi} \int_0^{\pi} (2 \cos t)^m \sin^2 t \, \mathrm{d}t = \frac{1}{2\pi} \int_{-2}^2 x^m \sqrt{4-x^2} \mathrm{d}x,$$
so that $\mathrm{d}\mu_{w_N}(x) = (2 \pi)^{-1} \sqrt{4-x^2} \mathrm{d}x$, and the operator $w_N$ is a semi-circular element.

Alternatively, if we take the limit as $n \rightarrow \infty$ in (\ref{eqn:momentsAn2}), we obtain
$$\int_{\mathbb{T}} \psi(u + u^{-1}) \mathrm{d}\varepsilon(u) = 2 \int_{\mathbb{T}} (u + u^{-1})^m \; \mathrm{Im}(u)^2 \mathrm{d}u,$$
where $\mathrm{d}u$ is the uniform measure over $\mathbb{T}$, as claimed in the previous section.

\subsection{Dynkin diagrams $D_n$}

\begin{figure}[tb]
\begin{center}
  \includegraphics[width=130mm]{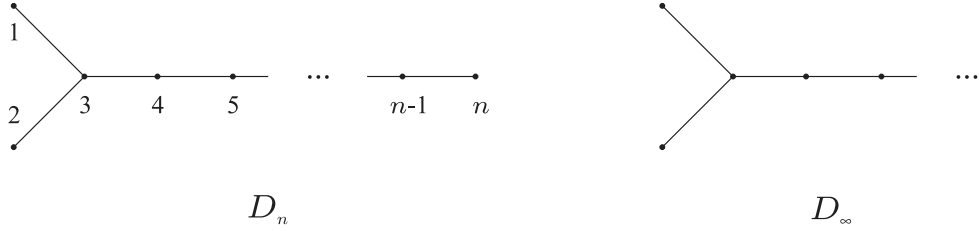}\\
 \caption{Dynkin diagrams $D_n$, $n=4,5,\ldots \;$, and $D_{\infty}$}\label{Fig:DynkinDiagramDn}
\end{center}
\end{figure}

For finite $n$, the distinguished vertex of the graph $D_n$ is the vertex $n$ in Figure \ref{Fig:DynkinDiagramDn}. The exponents $\mathrm{Exp}$ of $D_n$ are $1,3,5, \ldots, 2n-3, n-1$. For $n=2l$, the exponent $2l-1$ has multiplicity two, and we denote these exponents by $(2l-1,\pm)$. The eigenvectors of $D_{2l}$ are given by \cite[(B.6)]{behrend/pearce/petkova/zuber:2000} as:
\begin{eqnarray*}
& \psi^j_a = \sqrt{2} S_{2l+1-a,j}, \qquad \psi^j_1 = \psi^j_2 = \frac{1}{\sqrt{2}} S_{2l-1,j}, \qquad \psi^{(2l-1,\pm)}_a = S_{2l+1-a,2l-1}, & \\
& \psi^{(2l-1,\pm)}_{1+\epsilon} = \frac{1}{2} \left(S_{2l-1,2l-1} \pm (1-2\epsilon) \sqrt{(-1)^{l+1}}\right), &
\end{eqnarray*}
where $\epsilon = 0,1$, $a \neq 1,2$ and $j \in \mathrm{Exp}$, $j \neq 2l-1$. Using (\ref{eqn:moments_Dynkin}) and with $\widetilde{u} = e^{\pi i/(4l-2)}$,
\begin{eqnarray*}
\lefteqn{ \int_{\mathbb{T}} \psi(u + u^{-1}) \mathrm{d}\varepsilon(u) } \\
& = & \sum_{j \neq 2l-1} (2\cos(j \pi/(4l-2)))^m |\sqrt{2} S_{1,j}|^2 + 2(2\cos(j \pi/(4l-2)))^m |S_{1,j}|^2 \\
& = & \frac{4}{4l-2} \sum_{j \in \mathrm{Exp}} (2\cos(j \pi/(4l-2)))^m \sin^2(j \pi/(4l-2)) \\
& = & \frac{4}{4l-2} \sum_{j \in \mathrm{Exp}} (\widetilde{u}^j + \widetilde{u}^{-j})^m \; \mathrm{Im}(\widetilde{u}^j)^2 \\
& = & \frac{2}{4l-2} \sum_{j \in \{ 1,3,\ldots,8l-5 \}} (\widetilde{u}^j + \widetilde{u}^{-j})^m \; \mathrm{Im}(\widetilde{u}^j)^2 \;\; = \;\; 2 \int_{\mathbb{T}} (u + u^{-1})^m \; \mathrm{Im}(u)^2 \mathrm{d}_{4l-2}'u,
\end{eqnarray*}
where $\mathrm{d}_{4l-2}'$ is the uniform measure on the $(8l-4)^{\mathrm{th}}$ roots of unity of odd order.

For $D_{2l+1}$, the eigenvectors are given by \cite[(B.8)]{behrend/pearce/petkova/zuber:2000} as:
\begin{eqnarray*}
& \psi^j_a = (-1)^{\frac{j-1}{2}} \sqrt{2} S_{2l+2-a,j}, \qquad \psi^j_1 = \psi^j_2 = (-1)^{\frac{j-1}{2}} \frac{1}{\sqrt{2}} S_{2l,j} = \frac{1}{2\sqrt{l}}, & \\
& \psi^{2l}_a = 0, \qquad \psi^{2l}_1 = \frac{1}{\sqrt{2}}, \qquad \psi^{2l}_2 = -\frac{1}{\sqrt{2}}, &
\end{eqnarray*}
where $a \neq 1,2$ and $j \in \mathrm{Exp}\setminus\{2l\}$. Then, using (\ref{eqn:moments_Dynkin}) and with $\widetilde{u} = e^{\pi i/(4l)}$,
\begin{eqnarray*}
\lefteqn{ \int_{\mathbb{T}} \psi(u + u^{-1}) \mathrm{d}\varepsilon(u) \;\; = \;\; 2 \sum_{j \neq 2l} (2\cos(j \pi/4l))^m |S_{1,j}|^2 + 0} \\
& = & \frac{4}{4l} \sum_{j \in \{1,3,\ldots,4l-1\}} (2\cos(j \pi/4l))^m \sin^2(j \pi/4l) \\
& = & \frac{2}{4l} \sum_{j \in \{1,3,\ldots,8l-1\}} (\widetilde{u}^j + \widetilde{u}^{-j})^m \; \mathrm{Im}(\widetilde{u}^j)^2 \;\; = \;\; 2 \int_{\mathbb{T}} (u + u^{-1})^m \; \mathrm{Im}(u)^2 \mathrm{d}_{4l}'u.
\end{eqnarray*}
So the spectral measure $\varepsilon(u)$ on $\mathbb{T}$ for $D_n$ is given by $\mathrm{d}\varepsilon(u) = \alpha(u) \mathrm{d}_{2n-2}'u$, where
\begin{equation} \label{def:alpha(u)}
\alpha(u) = 2\mathrm{Im}(u)^2,
\end{equation}
which recovers the spectral measure given in \cite[Theorem 3.2]{banica/bisch:2007}.

Taking the limit of the graph $D_n$ as $n \rightarrow \infty$ with the vertex $n$ as the distinguished vertex, we just obtain the infinite graph $A_{\infty}$. In order to obtain the infinite graph $D_{\infty}$ we must set the distinguished vertex $\ast$ of $D_n$ to be the vertex 1 in Figure \ref{Fig:DynkinDiagramDn}. Then using (\ref{eqn:moments_Dynkin}), and taking the limit as $n \rightarrow \infty$, we obtain the spectral measure for $D_{\infty}$.

\subsection{Dynkin diagram $E_6$}

For $E_6$ the exponents are 1, 4, 5, 7, 8, 11. The eigenvectors for $E_6$ are given in \cite[(B.9)]{behrend/pearce/petkova/zuber:2000}. In particular,
$$\psi^1_1 = \psi^{11}_1 = \frac{1}{2} \sqrt{\frac{3-\sqrt{3}}{6}}, \qquad \psi^4_1 = \psi^8_1 = \frac{1}{2}, \qquad \psi^5_1 = \psi^7_1 = \frac{1}{2} \sqrt{\frac{3+\sqrt{3}}{6}}.$$
Then, by (\ref{eqn:moments_Dynkin}),
$$\int_{\mathbb{T}} \psi(u + u^{-1}) \mathrm{d}\varepsilon(u) = \sum_{j \in \mathrm{Exp}} |\psi^j_1|^2 (2\cos(j\pi/12))^m = \frac{1}{2} \sum_{p \in B_6} |\psi^p_1|^2 (2\cos(p\pi/12))^m,$$
where $B_6 = \{ 1,4,5,7,8,11,13,16,17,19,20,23 \}$, and for $j >12$ we define $\psi^j_1$ by $\psi^j_1 = \psi^{24-j}_1$. Then with $\widetilde{u} = e^{\pi i/12}$,
$$\int_{\mathbb{T}} \psi(u + u^{-1}) \mathrm{d}\varepsilon(u) = \frac{1}{24} \sum_{p \in B_6} 12 |\psi^p_1|^2 (\widetilde{u}^{p} + \widetilde{u}^{-p})^m.$$
Now for any $p \in B_6$, $\widetilde{u}^{p}$ is a $24^{\mathrm{th}}$ root of unity, but for $p = 4,8,16,20$, $\widetilde{u}^{p}$ is also a $6^{\mathrm{th}}$ root of unity. Since $|\psi^p_1|^2$ takes different values for different $p$, clearly we cannot write the above summation as an integral using the uniform measure over $24^{\mathrm{th}}$ roots of unity. However, with $\alpha$ as in (\ref{def:alpha(u)}), we have $\alpha(\widetilde{u}^{p}) = 12 |\psi^p_1|^2 - \alpha_p$, where $\alpha_p = 1/2$ for $p=1,5,7,11,13,17,19,23$ and $\alpha_p = 3/2$ for $p=4,8,16,20$.

By considering $a_p = \alpha(\widetilde{u}^{p}) + 1/2$, we can write
\begin{eqnarray*}
\lefteqn{ \int_{\mathbb{T}} \psi(u + u^{-1}) \mathrm{d}\varepsilon(u) \;\; = \;\; \frac{1}{24} \sum_{p \in B_6} a_p (\widetilde{u}^{p} + \widetilde{u}^{-p})^m } \\
& & - \frac{1}{24} \left( (\widetilde{u}^4 + \widetilde{u}^{-4})^m + (\widetilde{u}^8 + \widetilde{u}^{-8})^m + (\widetilde{u}^{16} + \widetilde{u}^{-16})^m + (\widetilde{u}^{20} + \widetilde{u}^{-20})^m \right).
\end{eqnarray*}
Since $\widetilde{u}^p$ is also a $6^{\mathrm{th}}$ root of unity for $p=4,8,16,20$, it may be possible to obtain the last four terms by considering an integral using the uniform measure on $6^{\mathrm{th}}$ roots of unity.
First, we consider the integral $\int (u+u^{-1})^m (2 \mathrm{Im}(u)^2 + 1/2) \mathrm{d}_{12} u$, where $\mathrm{d}_{12}$ is the uniform measure on the $24^{\mathrm{th}}$ roots of unity, to obtain the terms in the summation above, giving
\begin{eqnarray*}
\lefteqn{ \int_{\mathbb{T}} \psi(u + u^{-1}) \mathrm{d}\varepsilon(u) } \\
& = & \int_{\mathbb{T}} (u+u^{-1})^m (2 \mathrm{Im}(u)^2 + \textstyle \frac{1}{2} \displaystyle) \mathrm{d}_{12} u \; - \frac{1}{24} \sum_q a_q (\widetilde{u}^q + \widetilde{u}^{-q})^m \\
& & + \frac{1}{24} \left( (\widetilde{u}^4 + \widetilde{u}^{-4})^m + (\widetilde{u}^8 + \widetilde{u}^{-8})^m + (\widetilde{u}^{16} + \widetilde{u}^{-16})^m + (\widetilde{u}^{20} + \widetilde{u}^{-20})^m \right), \qquad
\end{eqnarray*}
where the summation is over $q \in \{ 2,3,6,9,10,12,14,15,18,21,22,24 \}$, that is, the integers $1 \leq q \leq 24$ such that $q \not \in B_6$. For these values of $q$, we have $a_2 = a_{10} = a_{14} = a_{22} = 1$, $a_3 = a_9 = a_{15} = a_{21} = 3/2$, $a_6 = a_{18} = 5/2$, and $a_{12} = a_{24} = 1/2$.
Using these values for $a_q$, we now isolate the terms involving the $12^{\mathrm{th}}$ roots of unity, giving
\begin{eqnarray*}
\lefteqn{ \int_{\mathbb{T}} \psi(u + u^{-1}) \mathrm{d}\varepsilon(u) } \\
& = & \int_{\mathbb{T}} (u+u^{-1})^m (2 \mathrm{Im}(u)^2 + \textstyle \frac{1}{2} \displaystyle) \mathrm{d}_{12} u - \frac{1}{24} \sum_{k=1}^{12} (\widetilde{u}^{2k} + \widetilde{u}^{-2k})^m \\
& & - \frac{1}{16} (\widetilde{u}^{3} + \widetilde{u}^{-3})^m + \frac{1}{12} (\widetilde{u}^{4} + \widetilde{u}^{-4})^m - \frac{1}{16} (\widetilde{u}^{6} + \widetilde{u}^{-6})^m + \frac{1}{12} (\widetilde{u}^{8} + \widetilde{u}^{-8})^m \\
& & - \frac{1}{16} (\widetilde{u}^{9} + \widetilde{u}^{-9})^m + \frac{1}{48} (\widetilde{u}^{12} + \widetilde{u}^{-12})^m - \frac{1}{16} (\widetilde{u}^{15} + \widetilde{u}^{-15})^m + \frac{1}{12} (\widetilde{u}^{16} + \widetilde{u}^{-16})^m \\
& & - \frac{1}{16} (\widetilde{u}^{18} + \widetilde{u}^{-18})^m + \frac{1}{12} (\widetilde{u}^{20} + \widetilde{u}^{-20})^m - \frac{1}{16} (\widetilde{u}^{21} + \widetilde{u}^{-21})^m + \frac{1}{48} (\widetilde{u}^{24} + \widetilde{u}^{-24})^m .
\end{eqnarray*}
Now $\sum_{k=1}^{12} (\widetilde{u}^{2k} + \widetilde{u}^{-2k})^m/12 = \int (u+u^{-1})^m  \mathrm{d}_6 u$. For the remaining terms, we notice that $\sum_{k=1}^{8} (\widetilde{u}^{3k} + \widetilde{u}^{-3k})^m/8 = \int (u+u^{-1})^m  \mathrm{d}_4 u$, giving
\begin{eqnarray*}
\lefteqn{ \int_{\mathbb{T}} \psi(u + u^{-1}) \mathrm{d}\varepsilon(u) } \\
& = & \int_{\mathbb{T}} (u+u^{-1})^m (2 \mathrm{Im}(u)^2 + \textstyle \frac{1}{2} \displaystyle) \mathrm{d}_{12} u - \frac{1}{2} \int_{\mathbb{T}} (u+u^{-1})^m \mathrm{d}_6 u - \frac{1}{2} \int_{\mathbb{T}} (u+u^{-1})^m \mathrm{d}_4 u \\
& & + \frac{1}{12} (\widetilde{u}^{4} + \widetilde{u}^{-4})^m + \frac{1}{12} (\widetilde{u}^{8} + \widetilde{u}^{-8})^m + \frac{1}{12} (\widetilde{u}^{12} + \widetilde{u}^{-12})^m \\
& & + \frac{1}{12} (\widetilde{u}^{16} + \widetilde{u}^{-16})^m + \frac{1}{12} (\widetilde{u}^{20} + \widetilde{u}^{-20})^m + \frac{1}{12} (\widetilde{u}^{24} + \widetilde{u}^{-24})^m.
\end{eqnarray*}
These last six terms are given by the integral $\int (u+u^{-1})^m  \mathrm{d}_3 u/2$ over $\mathbb{T}$.
Then the spectral measure $\varepsilon(u)$ (over $\mathbb{T}$) for $E_6$ is
$\mathrm{d}\varepsilon = \alpha \mathrm{d}_{12} + (\mathrm{d}_{12} - \mathrm{d}_6 - \mathrm{d}_4 + \mathrm{d}_3)/2$, which recovers the spectral measure given in \cite[Theorem 6.2]{banica/bisch:2007}.

\subsection{Dynkin diagrams $E_7$, $E_8$}

\begin{definition}
${}$ (\cite[Def. 7.1]{banica/bisch:2007})
A discrete measure supported by roots of unity is called \textbf{cyclotomic} if it is a linear combination of measures of type $\mathrm{d}_n$, $n \geq 1$, and $\alpha \mathrm{d}_n$, $n \geq 2$.
\end{definition}
Note that since $\mathrm{d}_n' = 2\mathrm{d}_{2n} - \mathrm{d}_n$, all the measures for the $A$ and $D$ diagrams, as well as for $E_6$, have been cyclotomic. However, Banica and Bisch \cite{banica/bisch:2007} proved that the spectral measures for $E_7$, $E_8$ are not cyclotomic. This can also be seen by our method using (\ref{eqn:moments_Dynkin}).

For $E_7$ the exponents are 1, 5, 7, 9, 11, 13, 17. The eigenvectors $\psi^j_1$ for $E_7$ are given by $\psi^j_1 = \surd (S_{1j} \sum_{i \in P} S_{ij})$, where $S$ is the $S$-matrix for $SU(2)_{16}$ and $P = \{ 1,9,17 \}$ \cite{behrend/pearce/petkova/zuber:2000}.
Then
$$\int_{\mathbb{T}} \psi(u + u^{-1}) \mathrm{d}\varepsilon(u) = \sum_{j \in \mathrm{Exp}} |\psi^j_1|^2 (2\cos(j\pi/18))^m = \frac{1}{2} \sum_{p \in B_7} |\psi^p_1|^2 (2\cos(p\pi/18))^m,$$
where $B_7 = \{ 1,5,7,9,11,13,17,19,23,25,27,29,31,35 \}$, and for $j >18$ we define $\psi^j_1$ by $\psi^j_1 = \psi^{36-j}_1$. Then with $\widetilde{u} = e^{\pi i/18}$,
\begin{equation} \label{eqn:sum_spec_measure_E7}
\int_{\mathbb{T}} \psi(u + u^{-1}) \mathrm{d}\varepsilon(u) = \frac{1}{36} \sum_{p \in B_7} 18 |\psi^p_1|^2 (\widetilde{u}^{p} + \widetilde{u}^{-p})^m.
\end{equation}
Now for any $p \in B_7$, $\widetilde{u}^{p}$ is a $36^{\mathrm{th}}$ root of unity, but not a root of unity of lower order, except for $p = 9,27$, in which case $\widetilde{u}^{p}$ is also a $4^{\mathrm{th}}$ root of unity. Since $|\psi^1_1|^2 \neq |\psi^5_1|^2$, clearly we cannot write the summation in (\ref{eqn:sum_spec_measure_E7}) as an integral using the uniform measure over $36^{\mathrm{th}}$ roots of unity. With $\alpha$ as in (\ref{def:alpha(u)}), and $\alpha_p = 18 |\psi^p_1|^2 - \alpha(\widetilde{u}^{p})$, we find that $\alpha_p = 0.4076$ for $p = 1,17,19,35$, $\alpha_p = 2.7057$ for $p = 5,13,23,31$, $\alpha_p = -0.1133$ for $p = 7,11,25,29$, and $\alpha_p = 4$ for $p = 9,27$.
Since $\alpha(\widetilde{u}^{p}) - 18 |\psi^p_1|^2$ also takes different values for certain $p \in B_7$, and for any $p \in B_7$, $\widetilde{u}^{p}$ is a $36^{\mathrm{th}}$ root of unity, but not a root of unity of lower order, the summation in (\ref{eqn:sum_spec_measure_E7}) cannot be written as an integral using the measure $\alpha \mathrm{d}_{18}$ either. So we see that the spectral measure for $E_7$ is not cyclotomic.

For $E_8$ the exponents are 1, 7, 11, 13, 17, 19, 23, 29. The eigenvectors $\psi^j_1$ for $E_8$ are given ny $\psi^j_1 = \surd (S_{1j} \sum_{i \in P} S_{ij})$, where $S$ is the $S$-matrix for $SU(2)_{28}$ and $P = \{ 1,11,19,29 \}$ \cite{behrend/pearce/petkova/zuber:2000}.
Then
\begin{equation} \label{eqn:sum_spec_measure_E8}
\int_{\mathbb{T}} \psi(u + u^{-1}) \mathrm{d}\varepsilon(u) = \sum_{j \in \mathrm{Exp}} |\psi^j_1|^2 (2\cos(j\pi/30))^m = \frac{1}{60} \sum_{p \in B_8} 30 |\psi^p_1|^2 (\widetilde{u}^{p} + \widetilde{u}^{-p})^m,
\end{equation}
where $\widetilde{u} = e^{\pi i/30}$, $B_8 = \{ 1,7,11,13,17,19,23,29,31,37,41,43,47,49,53,59 \}$, and for $j >30$ we define $\psi^j_1$ by $\psi^j_1 = \psi^{60-j}_1$.
With $\alpha_p = 30 |\psi^p_1|^2 - \alpha(\widetilde{u}^{p})$, we find that $\alpha_p = 0.4038$ for $p = 1,29,31,59$, $\alpha_p = 3.5135$ for $p = 7,23,37,53$, $\alpha_p = 2.0511$ for $p = 11,19,41,49$, and $\alpha_p = 4.5316$ for $p = 13,17,43,47$.
Now for all $p \in B_8$, $\widetilde{u}^{p}$ is a $60^{\mathrm{th}}$ root of unity, but not a root of unity of lower order. By similar considerations as in the case of $E_7$, we see that the summation in (\ref{eqn:sum_spec_measure_E8}) cannot be written as an integral using the uniform measure $\mathrm{d}_{30}$ or the measure $\alpha \mathrm{d}_{30}$ either. So we see that the spectral measure for $E_8$ is not cyclotomic.

However, in \cite{banica:2007}, Banica found explicit formulae for the spectral measures of $E_7$, $E_8$, using the densities $\alpha_j = \mathrm{Re}(1-u^{2j}) = 2 \mathrm{Im}(\widetilde{u}^{j})^2$, for $j=1,2,3$, where $\alpha = \alpha_1$ is the density in (\ref{def:alpha(u)}). A further simplification of the measures for these two graphs was obtained by considering the discrete measure $\mathrm{d}_n''=(3\mathrm{d}_{3n}'-\mathrm{d}_n')/2$, which is the uniform measure on the $12n^{\mathrm{th}}$ roots of unity of order $6k\pm1$. The support of the spectral measure over $\mathbb{T}$ for $E_6$, $E_7$, $E_8$ respectively basically coincides with the support of $\mathrm{d}_2''$, $\mathrm{d}_3''$, $\mathrm{d}_5''$ respectively, which can be easily seen from (\ref{eqn:moments_Dynkin}).

For $E_7$, (\ref{eqn:sum_spec_measure_E7}) gives that the spectral measure as a discrete weighted measure on the $36^{\mathrm{th}}$ roots of unity of order $6k\pm1$, plus the Dirac measure on the points $i$, $-i$ with weights $|\psi^9_1|^2 /2 = 1/6$. Now for $p \in B_7$,
\begin{eqnarray*}
|\psi^p_1|^2 & = & S_{1,p}^2 + S_{1,p} S_{9,p} + S_{1,p} S_{17,p} \;\; = \;\; 2 S_{1,p}^2 + S_{1,p} S_{9,p} \\
& = & \frac{1}{9} (2 \sin^2(p \pi/18) + \sin(p \pi/18)\sin(9p \pi/18)),
\end{eqnarray*}
whilst with $\widetilde{u} = e^{\pi i/18}$,
$$\alpha_2 (\widetilde{u}^p) = 2 \mathrm{Im}(\widetilde{u}^{2p})^2 = 2\sin^2(2p\pi/18) = 4\sin^2(p \pi/18) - 4\sin^4(p\pi/18).$$
Since $3\sin(\pi/18) - 4\sin^3(\pi/18) = \sin(3\pi/18) = 1/2$, we can write
$$\sin(\pi/18) = \frac{1}{2} \left( \frac{-1+i\sqrt{3}}{2} \right)^{2/3} + \frac{1}{2} \left( \frac{-1+i\sqrt{3}}{2} \right)^{1/3} \frac{-1-i\sqrt{3}}{2},$$
where the third root of $(-1+i\sqrt{3})/2$ takes its value in $\{ e^{i \theta} | \; 0 \leq \theta < 2 \pi/3 \}$. Using this expression for $\sin(\pi/18)$ one can find $\sin(j\pi/18)$ for all $j = 1,\ldots,18$. Then it is easy to check the identity $\sin(9p\pi/18) = 6\sin(p\pi/18) - 8\sin^3(p\pi/18)$ for $p \in B_7$, $p \neq 9,27$. Then
\begin{eqnarray*}
|\psi^p_1|^2 & = & \frac{1}{9} (2 \sin^2(p \pi/18) + \sin(p \pi/18)(6\sin(p\pi/18) - 8\sin^3(p\pi/18))) \\
& = & \frac{1}{9} (8\sin^2(p \pi/18) - 8\sin^4(p\pi/18))) \;\; = \;\; \frac{1}{9} \alpha_2 (\widetilde{u}^p),
\end{eqnarray*}
and from (\ref{eqn:sum_spec_measure_E7})
\begin{eqnarray*}
\int_{\mathbb{T}} \psi(u + u^{-1}) \mathrm{d}\varepsilon(u) & = & \frac{2}{3} \, \frac{1}{12} \sum_{\stackrel{p \in B_7}{\scriptscriptstyle{p \neq 9,27}}} \alpha_2 (\widetilde{u}^p) (\widetilde{u}^{p} + \widetilde{u}^{-p})^m + \frac{1}{6} \left( (i+i^{-1})^m + (-i+(-i)^{-1})^m \right) \\
& = & \frac{2}{3} \int_{\mathbb{T}} \psi(u + u^{-1}) \alpha_2(u) \mathrm{d}_3''u + \frac{1}{3} \int_{\mathbb{T}} \psi(u + u^{-1}) \mathrm{d}_1'u.
\end{eqnarray*}

Thus the spectral measure $\varepsilon(u)$ (over $\mathbb{T}$) for $E_7$ is $\mathrm{d}\varepsilon = (2 \alpha_2 \mathrm{d}_3'' + \mathrm{d}_1')/3$, which recovers the spectral measure for $E_7$ given in \cite[Theorem 8.7]{banica:2007}.

For $E_8$, (\ref{eqn:sum_spec_measure_E8}) gives that the spectral measure as a discrete weighted measure on the $60^{\mathrm{th}}$ roots of unity of order $6k\pm1$. However we need to remove the contribution given by $e^{2 \pi i p/30}$ for $p=5,25,35,55$, which are the $12^{\mathrm{th}}$ roots of unity of order $6k\pm1$. Now for $p \in B_8$,
\begin{eqnarray*}
\frac{1}{2} |\psi^p_1|^2 & = & \frac{1}{2} (S_{1,p}^2 + S_{1,p} S_{11,p} + S_{1,p} S_{19,p} + S_{1,p} S_{29,p}) \;\; = \;\; S_{1,p}^2 + S_{1,p} S_{11,p} \\
& = & \frac{1}{15} (\sin^2(p \pi/30) + \sin(p \pi/30)\sin(11p \pi/30)),
\end{eqnarray*}
whilst with $\widetilde{u} = e^{\pi i/30}$,
$$\alpha_1 (\widetilde{u}^p) + \alpha_3 (\widetilde{u}^p) = 2 \mathrm{Im}(\widetilde{u}^{p})^2 + 2 \mathrm{Im}(\widetilde{u}^{3p})^2 = 2(\sin^2(p\pi/30) + \sin^2(11p\pi/30).$$
Now $3\sin(\pi/30) - 4\sin^3(\pi/30) = \sin(3\pi/30) = (-1+\sqrt{5})/4$, so we can solve this cubic in $\sin(\pi/30)$ to write $\sin(\pi/30) = (-1-\sqrt{5} + \sqrt{6}\sqrt{5-\sqrt{5}})/8$. Using this expression for $\sin(\pi/30)$ one can find $\sin(j\pi/30)$ for all $j = 1,\ldots,30$. Then it is easy to check the identity $\sin^2(3p\pi/30) = \sin(p\pi/30)\sin(11p\pi/30)$ for $p \in B_8$. Then
$$|\psi^p_1|^2 = \frac{1}{15} (\sin^2(p \pi/30) + \sin^2(3p\pi/30)) = \frac{1}{30} (\alpha_1 (\widetilde{u}^p) + \alpha_3 (\widetilde{u}^p)).$$
For $p=5,25,35,55$, $\alpha_1 (\widetilde{u}^p) + \alpha_3 (\widetilde{u}^p) = 5/2$.
Then from (\ref{eqn:sum_spec_measure_E8})
\begin{eqnarray*}
\lefteqn{ \int_{\mathbb{T}} \psi(u + u^{-1}) \mathrm{d}\varepsilon(u) \;\; = \;\; \frac{2}{3} \, \frac{1}{20} \sum_{p \in B_8} (\alpha_1 (\widetilde{u}^p) + \alpha_3 (\widetilde{u}^p)) (\widetilde{u}^{p} + \widetilde{u}^{-p})^m } \\
& = & \frac{2}{3} \int_{\mathbb{T}} \psi(u + u^{-1}) (\alpha_1(u) + \alpha_3(u)) \mathrm{d}_5''u - \frac{2}{15} \int_{\mathbb{T}} \psi(u + u^{-1}) (\alpha_1(u) + \alpha_3(u)) \mathrm{d}_1''u \\
& = & \frac{2}{3} \int_{\mathbb{T}} \psi(u + u^{-1}) (\alpha_1(u) + \alpha_3(u)) \mathrm{d}_5''u - \frac{1}{3} \int_{\mathbb{T}} \psi(u + u^{-1}) \mathrm{d}_1''u.
\end{eqnarray*}

Thus the spectral measure $\varepsilon(u)$ (over $\mathbb{T}$) for $E_8$ is $\mathrm{d}\varepsilon = (2 (\alpha_1 + \alpha_3) \mathrm{d}_5'' - \mathrm{d}_1'')/3$, which recovers the spectral measure for $E_8$ given in \cite[Theorem 8.7]{banica:2007}.

\section{Spectral measures for the finite subgroups of $SU(2)$} \label{Sect:Spectral_Measures_for_subgroup_SU(2)}

\begin{figure}[tb]
\begin{center}
  \includegraphics[width=150mm]{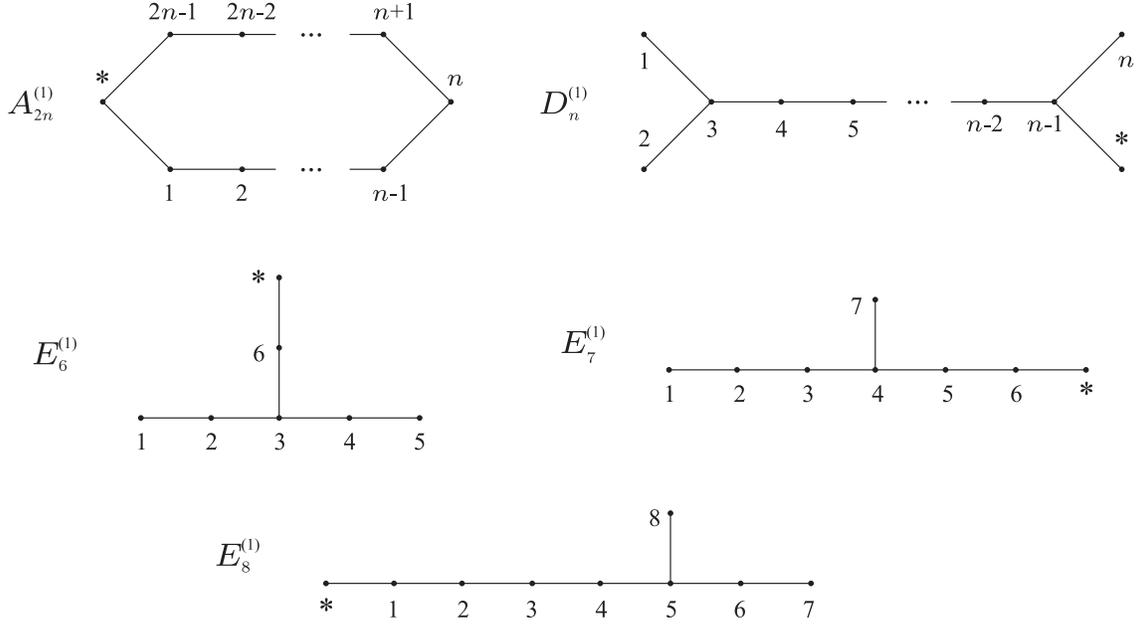}\\
 \caption{Affine Dynkin diagrams $A_{2n}^{(1)}$, $D_{2n}^{(1)}$, $E_6^{(1)}$, $E_7^{(1)}$ and $E_8^{(1)}$}\label{Fig:AffineDynkinADE}
\end{center}
\end{figure}

The McKay correspondence \cite{mckay:1980} associates to every finite subgroup $\Gamma$ of $SU(2)$ an affine Dynkin diagram $\mathcal{G}_{\Gamma}$ given by the fusion graph of the fundamental representation $\rho$ acting on the irreducible representations of $\Gamma$. These affine Dynkin diagrams are illustrated in Figure \ref{Fig:AffineDynkinADE}, where $\ast$ denotes the identity representation.
Hence there is associated to each finite subgroup of $SU(2)$ the corresponding (non-affine) $ADE$ Dynkin diagram $\mathcal{G}$, which is obtained from the affine diagram by deleting the vertex $\ast$ and all edges attached to it. This correspondence is shown in the following table. The second column indicates the type of the associated modular invariant.

\begin{center}
\begin{tabular}{|c|c|c|c|} \hline
Dynkin Diagram $\mathcal{G}$ & Type & Subgroup $\Gamma \subset SU(2)$
& $|\Gamma|$ \\
\hline\hline $A_l$ & I & cyclic, $\mathbb{Z}_{l+1}$ & $l+1$ \\
\hline $D_{2k}$ & I & binary dihedral, $BD_{2k} = Q_{2k-2}$ & $8k-8$ \\
\hline $D_{2k+1}$ & II & binary dihedral, $BD_{2k+1} = Q_{2k-1}$ & $8k-4$ \\
\hline $E_6$ & I & binary tetrahedral, $BT = BA_4$ & 24 \\
\hline $E_7$ & II & binary octahedral, $BO = BS_4$ & 48 \\
\hline $E_8$ & I & binary icosahedral, $BI = BA_5$ & 120 \\
\hline
\end{tabular}
\end{center}

It was shown in \cite{kawai:1989} that for any finite group $\Gamma$ the $S$-matrix, which simultaneously diagonalizes the representations of $\Gamma$, can be written in terms of the characters $\chi_j(\Gamma_i)$ of $\Gamma$ evaluated on the conjugacy classes $\Gamma_i$ of $\Gamma$, $S_{ij} = \sqrt{|\Gamma_j|} \chi_i(\Gamma_j)/\sqrt{|\Gamma|}$. Let $N_{\rho}$ be the fundamental representation matrix of the fusion rules of the irreducible characters of $\Gamma$. Then by the Verlinde formula (\ref{eqn:verlinde_formula}), the eigenvalues of $N_{\rho}$ are given by ratios of the $S$-matrix, $\sigma(N_{\rho}) = \{ S_{\rho,j}/S_{\rho,0} | j = 1,\ldots, p \}$, where $p$ is the number of conjugacy classes and $\rho$ is the fundamental representation of $G$. Now
$$\frac{\sqrt{|\Gamma_j|} \chi_{\rho}(\Gamma_j)/\sqrt{|\Gamma|}}{\sqrt{|\Gamma_j|} \chi_{\rho}(\Gamma_0)/\sqrt{|\Gamma|}} = \chi_{\rho}(\Gamma_j),$$
since $\chi_{\rho}(\Gamma_0) = 1$. Then any eigenvalue of $\Gamma$ can be written in the form $\chi_{\rho}(g) = \mathrm{Tr}(\rho(g))$, where $g$ is any element of $\Gamma_j$.

The elements $y_i$ in (\ref{eqn:moments_general_SU(2)graph}) are then given by $y_i = S_{0,j} = \sqrt{|\Gamma_j|} \chi_0(\Gamma_j)/\sqrt{|\Gamma|} = \sqrt{|\Gamma_j|}/\sqrt{|\Gamma|}$.
Then the $m^{\mathrm{th}}$ moment $\varsigma_m$ is given by
\begin{equation} \label{eqn:moments-subgroupSU(2)}
\varsigma_m \; = \; \int x^m \mathrm{d}\mu(x) \; = \; \sum_{j=1}^n \frac{|\Gamma_j|}{|\Gamma|} \chi_{\rho} (\Gamma_j)^m.
\end{equation}
We define an inverse $\Phi^{-1}:[-2,2] \rightarrow \mathbb{T}$ of the map $\Phi$ given in (\ref{def:phi-SU(2)}) by
\begin{equation} \label{def:inverse_phi-SU(2)}
\Phi^{-1}(x) = (x + i \sqrt{4-x^2})/2,
\end{equation}
for $x \in [-2,2]$. Then the spectral measure of $\Gamma$ (over $\mathbb{T}$) is given by
\begin{equation} \label{eqn:moments-subgroupSU(2)2}
\int_{\mathbb{T}} \psi(u + u^{-1}) \mathrm{d}\varepsilon(u) \; = \; \sum_{j=1}^n \frac{|\Gamma_j|}{|\Gamma|} (\Phi^{-1}(\chi_{\rho} (\Gamma_j)) + \overline{\Phi^{-1}(\chi_{\rho} (\Gamma_j))})^m.
\end{equation}

The generating series of the moments $G(q) = \sum_{m=0}^{\infty} \varsigma_m q^m = \int_{\mathbb{T}} (1-qu)^{-1} \mathrm{d}\varepsilon(u)$, is
\begin{equation} \label{eqn:Generating_series_momentsSU(2)}
G(q) = \sum_{m=0}^{\infty} \sum_{j=1}^n \frac{|\Gamma_j|}{|\Gamma|} \chi_{\rho} (\Gamma_j)^m q^m = \sum_{j=1}^n \frac{|\Gamma_j|}{|\Gamma|} \frac{1}{1 - q \chi_{\rho}(\Gamma_j)}.
\end{equation}

\subsection{Cyclic Group $\mathbb{Z}_{2n}$}

Suppose $\Gamma$ is the cyclic subgroup $\mathbb{Z}_{2n}$ of $SU(2)$, which has McKay graph $A_{2n}^{(1)}$. Then $|\Gamma| = 2n$, and each element of the group is a separate conjugacy class. Now $\chi_{\rho}(\Gamma_j) = \widetilde{u}^j + \widetilde{u}^{-j} \in [-2,2]$, where $\widetilde{u} = e^{\pi i/n}$, for each $j=1,\ldots,2n$. Then by (\ref{eqn:moments-subgroupSU(2)})
$$\int_{\mathbb{T}} \psi(u + u^{-1}) \mathrm{d}\varepsilon(u) = \sum_{j=1}^{2n} \frac{1}{2n} (\widetilde{u}^j + \widetilde{u}^{-j})^m = \int_{\mathbb{T}} (u + u^{-1})^m \; \mathrm{d}_n u.$$
Hence the spectral measure for $A_{2n}^{(1)}$ (over $\mathbb{T}$) is $\mathrm{d} \varepsilon(u) = \mathrm{d}_n u$, as in \cite[Theorem 2.1]{banica/bisch:2007}.

\subsection{Binary Dihedral Group $BD_n$}

Let $\Gamma$ be the binary dihedral group $BD_n = \langle \sigma, \tau | \tau^2 = \sigma^n = (\tau \sigma)^2 \rangle$, which has McKay graph $D_{n}^{(1)}$. Then $|\Gamma| = 4(n-2)$. The character table for $BD_n$ is given in Table \ref{table:Character_table-BDn}.
\begin{table}[bt]
\begin{center}
\begin{tabular}{|c||c|c|c|c|c|} \hline
$\Gamma_j$ & $\mathbf{1}$ & $(\tau \sigma)^2$ & $\sigma^j, \; j=1,\ldots,n-3$ & $\tau$ & $\tau \sigma$ \\
\hline $|\Gamma_j|$ & 1 & 1 & 2 & $n-2$ & $n-2$ \\
\hline $\chi_{\rho}(\Gamma_j) \in [-2,2]$ & 2 & $-2$ & $\xi^j + \xi^{-j}$ & 0 & 0 \\
\hline $e^{2 \pi i \theta} = \Phi^{-1}(\chi_{\rho}(\Gamma_j)) \in \mathbb{T}$ & 1 & $-1$ & $\xi^j$ & $i$ & $-i$ \\
\hline $\theta \in [0,1]$ & 0 & $\frac{n-2}{2(n-2)}$ & $\frac{j}{2(n-2)}$ & $\frac{n-2}{4(n-2)}$ & $\frac{3(n-2)}{4(n-2)}$ \\
\hline
\end{tabular} \\
\caption{Character table for $BD_n$. Here $\xi = e^{\pi i/(n-2)}$.} \label{table:Character_table-BDn}
\end{center}
\end{table}
Let $\widetilde{u} = e^{\pi i/2(n-2)}$ and $U(j) = (\widetilde{u}^j + \widetilde{u}^{-j})^m$. Then by (\ref{eqn:moments-subgroupSU(2)})
\begin{eqnarray*}
\lefteqn{ \int_{\mathbb{T}} \psi(u + u^{-1}) \mathrm{d}\varepsilon(u) } \\
& = & \frac{1}{4(n-2)} U(0) + \frac{1}{4(n-2)} U(n-2) + \sum_{j=1}^{n-3} \frac{2}{4(n-2)} \left( \frac{U(j) + U(2n-2-j)}{2} \right) \\
& & + \frac{n-2}{4(n-2)} U((n-2)/2) + \frac{n-2}{4(n-2)} U(3(n-2)/2) \\
& = & \sum_{j=0}^{2n-3} \frac{1}{4(n-2)} (\widetilde{u}^j + \widetilde{u}^{-j})^m + \frac{1}{4} \left( (\widetilde{u}^{(n-2)/2} + \widetilde{u}^{-(n-2)/2})^m + (\widetilde{u}^{3(n-2)/2} + \widetilde{u}^{-3(n-2)/2})^m \right) \\
& = & \frac{1}{2} \int_{\mathbb{T}} (u + u^{-1})^m \; \mathrm{d}_{n-2} u + \frac{1}{4} \int_{\mathbb{T}} (u + u^{-1})^m \; (\delta_i + \delta_{-i}),
\end{eqnarray*}
where $\delta_{\omega}$ is the Dirac measure at $\omega \in \mathbb{T}$.
Then the spectral measure for $D_{n}^{(1)}$ (over $\mathbb{T}$) is as given in \cite[Theorem 4.1]{banica/bisch:2007}:
$$\mathrm{d} \varepsilon(u) = \frac{1}{2} \mathrm{d}_{n-2} u + \frac{1}{4}(\delta_i + \delta_{-i}) = \frac{1}{2} \mathrm{d}_{n-2} u + \frac{1}{2} \mathrm{d}_1' u.$$

\subsection{Binary Tetrahedral Group $BT$}

Let $\Gamma$ be the binary tetrahedral group $BT$, which has McKay graph $E_6^{(1)}$. It has order 24, and is generated by $BD_4 = \langle \sigma, \tau \rangle$ and $\mu$:
$$\sigma = \left( \begin{array}{cc} i & 0 \\ 0 & -i \end{array} \right), \qquad \tau = \left( \begin{array}{cc} 0 & 1 \\ -1 & 0 \end{array} \right), \qquad \mu = \frac{1}{\sqrt{2}} \left( \begin{array}{cc} \varepsilon^7 & \varepsilon^7 \\ \varepsilon^5 & \varepsilon \end{array} \right),$$
where $\varepsilon = e^{2\pi i/ 8}$. The orders of the group elements $\sigma$, $\tau$, $\mu$ are 4, 4, 6 respectively.
The character table for $BT$ is given in Table \ref{table:Character_table-BT}.
\begin{table}[bt]
\begin{center}
\renewcommand{\arraystretch}{1.5}
\begin{tabular}{|c||c|c|c|c|c|c|c|} \hline
$\Gamma_j$ & $\mathbf{1}$ & $-\mathbf{1}$ & $\tau$ & $\mu$ & $\mu^2$ & $\mu^4$ & $\mu^5$ \\
\hline $|\Gamma_j|$ & 1 & 1 & 6 & 4 & 4 & 4 & 4 \\
\hline $\chi_{\rho}(\Gamma_j) \in [-2,2]$ & 2 & $-2$ & 0 & 1 & $-1$ & $-1$ & 1 \\
\hline $e^{2 \pi i \theta} = \Phi^{-1}(\chi_{\rho}(\Gamma_j)) \in \mathbb{T}$ & 1 & $-1$ & $i$ & $e^{\pi i/3}$ & $e^{2\pi i/3}$ & $e^{2\pi i/3}$ & $e^{\pi i/3}$\\
\hline $\theta \in [0,1]$ & 0 & $\frac{1}{2}$ & $\frac{1}{4}$ & $\frac{1}{6}$ & $\frac{1}{3}$ & $\frac{2}{3}$ & $\frac{5}{6}$ \\
\hline
\end{tabular} \\
\caption{Character table for the binary tetrahedral group $BT$.} \label{table:Character_table-BT}
\end{center}
\end{table}
Let $\widetilde{u} = e^{2\pi i/12}$ and $U(j) = (\widetilde{u}^j + \widetilde{u}^{-j})^m$. Then by (\ref{eqn:moments-subgroupSU(2)}), the integral $\int_{\mathbb{T}} \psi(u + u^{-1}) d\varepsilon(u)$ is equal to
$$\frac{1}{24} U(0) + \frac{1}{24} U(6) + \frac{6}{24} U(3) + \frac{4}{24} U(2) + \frac{4}{24} U(4) + \frac{4}{24} U(8) + \frac{4}{24} U(10).$$
For the $6^{\mathrm{th}}$ roots of unity we have $\alpha(e^{p \pi i/6}) - 1/2 = - 1/2$, $p=0,6$, and $\alpha(e^{p \pi i/6}) - 1/2 = 1$, $p=2,4,8,10$, where $\alpha$ is given in (\ref{def:alpha(u)}). Then since $U(3) = U(9)$:
\begin{eqnarray*}
\int_{\mathbb{T}} \psi(u + u^{-1}) \mathrm{d}\varepsilon(u) & = & \frac{3}{24} (U(0) + U(3) + U(6) + U(9)) \\
& & + \frac{1}{24}(-2U(0) + 4U(2) + 4U(2) -2U(6) + 4U(8) + 4U(10)) \\
& = & \frac{1}{2} \sum_{j=0}^{3} \frac{1}{4} (\widetilde{u}^{3j} + \widetilde{u}^{-3j})^m + \sum_{j=0}^{5} \frac{1}{6} (\alpha(\widetilde{u}^{2j}) - \textstyle \frac{1}{2} \displaystyle) (\widetilde{u}^{2j} + \widetilde{u}^{-2j})^m \\
& = & \frac{1}{2} \int_{\mathbb{T}} (u + u^{-1})^m \; \mathrm{d}_2 u + \int_{\mathbb{T}} (u + u^{-1})^m \; (\alpha(u) - \textstyle \frac{1}{2} \displaystyle) \mathrm{d}_3 u.
\end{eqnarray*}
Hence the spectral measure for $E_6^{(1)}$ (over $\mathbb{T}$) is $\mathrm{d} \varepsilon = (\alpha - 1/2) \mathrm{d}_3 + \mathrm{d}_2/2$, as given in \cite[Theorem 6.1]{banica/bisch:2007}.

\subsection{Binary Octahedral Group $BO$}

Let $\Gamma$ be the binary octahedral group $BO$, which has McKay graph $E_7^{(1)}$. It has order 48 and is generated by the binary tetrahedral group $BT$ and the element $\kappa$ of order 8 given by
$$\kappa = \left( \begin{array}{cc} \varepsilon & 0 \\ 0 & \varepsilon^7 \end{array} \right),$$
where again $\varepsilon = e^{2\pi i/ 8}$. Its McKay graph is $E_7^{(1)}$.
The character table for $BO$ is given in Table \ref{table:Character_table-BO}.
\begin{table}[bt]
\begin{center}
\renewcommand{\arraystretch}{1.5}
\begin{tabular}{|c||c|c|c|c|c|c|c|c|} \hline
$\Gamma_j$ & $\mathbf{1}$ & $-\mathbf{1}$ & $\mu$ & $\mu^2$ & $\tau$ & $\kappa$ & $\tau \kappa$ & $\kappa^3$ \\
\hline $|\Gamma_j|$ & 1 & 1 & 8 & 8 & 6 & 6 & 12 & 6 \\
\hline $\chi_{\rho}(\Gamma_j) \in [-2,2]$ & 2 & $-2$ & 1 & $-1$ & 0 & $\sqrt{2}$ & 0 & $-\sqrt{2}$ \\
\hline $e^{2 \pi i \theta} = \Phi^{-1}(\chi_{\rho}(\Gamma_j)) \in \mathbb{T}$ & 1 & $-1$ & $e^{\pi i/3}$ & $e^{2\pi i/3}$ & $i$ & $e^{\pi i/4}$ & $i$ & $e^{3\pi i/4}$ \\
\hline $\theta \in [0,1]$ & 0 & $\frac{1}{2}$ & $\frac{1}{6}$ & $\frac{1}{3}$ & $\frac{1}{4}$ & $\frac{1}{8}$ & $\frac{1}{4}$ & $\frac{3}{8}$ \\
\hline
\end{tabular} \\
\caption{Character table for the binary octahedral group $BO$.} \label{table:Character_table-BO}
\end{center}
\end{table}
Let $\widetilde{u} = e^{2\pi i/24}$ and $U(j) = (\widetilde{u}^j + \widetilde{u}^{-j})^m$. Then by (\ref{eqn:moments-subgroupSU(2)})
\begin{eqnarray*}
\lefteqn{ \int_{\mathbb{T}} \psi(u + u^{-1}) \mathrm{d}\varepsilon(u) } \\
& = & \frac{1}{48} U(0) + \frac{1}{48} U(12) + \frac{8}{48} U(4) + \frac{8}{48} U(8) + \frac{18}{48} U(6) + \frac{6}{48} U(3) + \frac{6}{48} U(9).
\end{eqnarray*}
For the $8^{\mathrm{th}}$ roots of unity we have $\alpha(e^{p \pi i/12}) - 1/2 = - 1/2$, for $p=0,12$, $\alpha(e^{p \pi i/12}) - 1/2 = 1/2$, for $p=3,9,15,21$, and $\alpha(e^{p \pi i/12}) - 1/2 = 3/2$, for $p=6,18$, where $\alpha$ is given in (\ref{def:alpha(u)}). Then since $U(j) = U(24-j)$, $j=1,\ldots,12$, we have
\begin{eqnarray*}
\lefteqn{ \int_{\mathbb{T}} \psi(u + u^{-1}) \mathrm{d}\varepsilon(u) } \\
& = & \frac{4}{48} (U(0) + U(4) + U(8) + U(12) + U(16) + U(20)) + \frac{1}{48} \Big(-3U(0) + 3U(3) \\
& & \qquad + 9U(6) + 3U(9) - 3U(12) + 3U(15) + 9 U(18) + 3U(21) \Big) \\
& = & \frac{1}{2} \sum_{j=0}^{5} \frac{1}{6} (\widetilde{u}^{4j} + \widetilde{u}^{-4j})^m + \sum_{j=0}^{7} \frac{1}{8} (\alpha(\widetilde{u}^{3j}) - \textstyle \frac{1}{2} \displaystyle) (\widetilde{u}^{3j} + \widetilde{u}^{-3j})^m \\
& = & \frac{1}{2} \int_{\mathbb{T}} (u + u^{-1})^m \; \mathrm{d}_3 u + \int_{\mathbb{T}} (u + u^{-1})^m \; (\alpha(u) - \textstyle \frac{1}{2} \displaystyle) \mathrm{d}_4 u.
\end{eqnarray*}
Hence the spectral measure for $E_7^{(1)}$ (over $\mathbb{T}$) is $\mathrm{d} \varepsilon = (\alpha - 1/2) \mathrm{d}_4 + \mathrm{d}_3/2$, as given in \cite[Theorem 6.1]{banica/bisch:2007}.

\subsection{Binary Icosahedral Group $BI$}

Let $\Gamma$ be the binary icosahedral group $BI$, which has McKay graph $E_8^{(1)}$. It has order 120, and is generated by $\sigma$, $\tau$:
$$\sigma = \left( \begin{array}{cc} - \varepsilon^3 & 0 \\ 0 & - \varepsilon^2 \end{array} \right), \qquad \tau = \frac{1}{\sqrt{5}} \left( \begin{array}{cc} \varepsilon^4 - \varepsilon & \varepsilon^2 - \varepsilon^3 \\ \varepsilon^2 - \varepsilon^3 & \varepsilon - \varepsilon^4 \end{array} \right),$$
where $\varepsilon = e^{2\pi i/ 5}$. The orders of $\sigma$, $\tau$ are 10, 4 respectively.
The character table for $BI$ is given in Table \ref{table:Character_table-BI}.
\begin{table}[bt]
\begin{center}
\renewcommand{\arraystretch}{1.5}
\begin{tabular}{|c||c|c|c|c|c|c|c|c|c|} \hline
$\Gamma_j$ & $\mathbf{1}$ & $-\mathbf{1}$ & $\sigma$ & $\sigma^2$ & $\sigma^3$ & $\sigma^4$ & $\tau$ & $\sigma^2 \tau$ & $\sigma^7 \tau$ \\
\hline $|\Gamma_j|$ & 1 & 1 & 12 & 12 & 12 & 12 & 30 & 20 & 20 \\
\hline $\chi_{\rho}(\Gamma_j) \in [-2,2]$ & 2 & $-2$ & $\mu^+$ & $-\mu^-$ & $\mu^-$ & $-\mu^+$ & 0 & $-1$ & 1 \\
\hline $e^{2 \pi i \theta} = \Phi^{-1}(\chi_{\rho}(\Gamma_j))$ & 1 & $-1$ & $e^{\pi i/5}$ & $e^{2\pi i/5}$ & $e^{3\pi i/5}$ & $e^{4\pi i/5}$ & $i$ & $e^{2\pi i/3}$ & $e^{\pi i/3}$ \\
\hline $\theta \in [0,1]$ & 0 & $\frac{1}{2}$ & $\frac{1}{10}$ & $\frac{1}{5}$ & $\frac{3}{10}$ & $\frac{2}{5}$ & $\frac{1}{4}$ & $\frac{1}{3}$ & $\frac{1}{6}$ \\
\hline
\end{tabular} \\
\caption{Character table for the binary icosahedral group $BI$. Here $\mu^{\pm} = (1 \pm \sqrt{5})/2$.} \label{table:Character_table-BI}
\end{center}
\end{table}
Let $\widetilde{u} = e^{2\pi i/60}$ and $U(j) = (\widetilde{u}^j + \widetilde{u}^{-j})^m$. By (\ref{eqn:moments-subgroupSU(2)})
\begin{eqnarray*}
\int_{\mathbb{T}} \psi(u + u^{-1}) \mathrm{d}\varepsilon(u) & = & \frac{1}{120} U(0) + \frac{1}{120} U(30) + \frac{12}{120} U(6) + \frac{12}{120} U(48) + \frac{12}{120} U(18) \\
& & + \frac{12}{120} U(36) + \frac{30}{120} U(15) + \frac{20}{120} U(20) + \frac{20}{120} U(10).
\end{eqnarray*}
For the $12^{\mathrm{th}}$ roots of unity we have $\alpha(e^{p \pi i/6}) - 1/2 = -1/2$, for $p= 0,6$, $\alpha(e^{p \pi i/6}) - 1/2 = 1$, for $p= 2,4,8,10$, $\alpha(e^{p \pi i/6}) - 1/2 = 3/2$, for $p= 3,9$, and $\alpha(e^{p \pi i/6}) - 1/2 = 0$, for $p= 1,5,7,11$, where $\alpha$ is given in (\ref{def:alpha(u)}). Then since $U(j) = U(60-j)$, $j=1,\ldots,30$, we have
\begin{eqnarray*}
\int_{\mathbb{T}} \psi(u + u^{-1}) \mathrm{d}\varepsilon(u) & = & \frac{6}{120} (U(0) + U(6) + U(12) + U(18) + U(24) + U(30) + U(36) \\
& & \qquad \quad + U(42) + U(48) + U(54)) \\
& & + \frac{1}{120}(-5U(0) + 10U(10) + 15U(15) + 10U(20) - 5U(30) \\
& & \qquad \quad + 10U(40) + 15U(45) + 10U(50)) \\
& = & \frac{1}{2} \sum_{j=0}^{9} \frac{1}{10} (\widetilde{u}^{6j} + \widetilde{u}^{-6j})^m + \sum_{j=0}^{11} \frac{1}{12} (\alpha(\widetilde{u}^{5j}) - \textstyle \frac{1}{2} \displaystyle) (\widetilde{u}^{5j} + \widetilde{u}^{-5j})^m \\
& = & \frac{1}{2} \int_{\mathbb{T}} (u + u^{-1})^m \; \mathrm{d}_5 u + \int_{\mathbb{T}} (u + u^{-1})^m \; (\alpha(u) - \textstyle \frac{1}{2} \displaystyle) \mathrm{d}_6 u.
\end{eqnarray*}
Hence the spectral measure for $E_8^{(1)}$ (over $\mathbb{T}$) is $\mathrm{d} \varepsilon = (\alpha - 1/2) \mathrm{d}_6 + \mathrm{d}_5/2$, as given in \cite[Theorem 6.1]{banica/bisch:2007}.

\section{Hilbert Series of dimensions of $ADE$ models.}

We now compare various polynomials related to $ADE$ models.

\subsection{$T$-Series} \label{Sect:T-series}

We begin first with the $T$-series of Banica and Bisch \cite{banica/bisch:2007}. Let $\mathcal{G}$ now be any bipartite graph with norm $\leq 2$, that is, its adjacency matrix $\Delta$ has norm $\leq 2$. These are the subgroups of $SU(2)$, with McKay graphs given by the affine Dynkin diagrams, and the modules and subgroups of $SU(2)_k$, which have McKay graphs given by the $ADE$ Dynkin diagrams.

Let $A(\mathcal{G})$ be the path algebra for $\mathcal{G}$, with initial vertex the distinguished vertex $\ast$ which has lowest Perron-Frobenius weight. The Hilbert series (also called Poincar\'{e} series in some literature)
\begin{equation} \label{eqn:f(z)}
f(z) = \sum_{k=0}^{\infty} \mathrm{dim}(A(\mathcal{G})_k) z^k
\end{equation}
of $\mathcal{G}$ is the generating function counting the numbers $l_{2k}$ of loops of length $2k$ on $\mathcal{G}$, from the vertex $\ast$ to itself, $f(z) = \sum_{k=0}^{\infty} l_{2k} z^k$. The Hilbert series $f$ measures the dimension of the algebra at level $k$ in the Bratteli diagram. If $\mathcal{G}$ is the principal graph of a subfactor $N \subset M$, the series $f$ measures the dimensions of the higher relative commutants, giving an invariant of the subfactor $N \subset M$. We define another function $\widehat{f}$ by
\begin{equation} \label{eqn:generalized-f(z)}
\widehat{f}(z) = \varphi \left( \left( \mathbf{1} - z^{\frac{1}{2}} \Delta \right)^{-1} \right).
\end{equation}
Then $\widehat{f}(z) = \varphi(\mathbf{1} + z^{1/2} \Delta + z \Delta^2 + z^{3/2} \Delta^3 + \ldots \; ) = \sum_{n=0}^{\infty} [\Delta^n]_{\ast,\ast} z^{n/2}$. Since $\mathcal{G}$ is bipartite, there are no paths of odd length from $\ast$ to $\ast$, and so $[\Delta^{2k+1}]_{\ast,\ast} = 0$ for $k=0,1,\ldots \;$. Then $\widehat{f}(z) = \sum_{k=0}^{\infty} [\Delta^{2k}]_{\ast,\ast} z^k = f(z)$. Then it is easily seen from (\ref{Stieltjes_transform}) and (\ref{eqn:generalized-f(z)}) that $f(z^2)$ is equal to the Stieltjes transform $\sigma(z)$ of $\mu_{\Delta}$.

Suppose $P$ is the ($A_1$-)planar algebra \cite{jones:planar} for a subfactor $N \subset M$ with Jones index $[M:N] < 4$ and principal graph $\mathcal{G}$. If $\mathrm{dim}(P_0^{\pm}) = 1$, the Hilbert series $f(z)$ is identical to the Hilbert series $\Phi_P(z)$ which gives the dimension of the planar algebra $P$:
$$\Phi_P(z) = \frac{1}{2} (\mathrm{dim}(P_0^+) + \mathrm{dim}(P_0^-)) + \sum_{j=1}^{\infty} \mathrm{dim}(P_j) z^j.$$
As a Temperley-Lieb module, $P$ decomposes into a sum of irreducible Temperley-Lieb modules, with the multiplicity of the irreducible module of lowest weight $k$ given by the non-negative integer $a_k$.
Jones \cite{jones:2001} then defined the series $\Theta$ by
$$\Theta_P(q) = \sum_{j=0}^{\infty} a_j q^j.$$
It was shown in \cite[Prop. 1.2]{banica/bisch:2007} that $\Theta(q^2) = 2G(q) + q^2 -1$, where $G(q)$ is the generating series of the moments of the spectral measure for $\mathcal{G}$, defined in Section \ref{Sect:Spectral_Measures_for_subgroup_SU(2)}. The series $\Theta(q)$ is essentially obtained from the Hilbert series $f(z)$ in (\ref{eqn:f(z)}) by a change of variables. More explicitly, in \cite{banica/bisch:2007}, $\Theta(q)$ is given in terms of $f(z)$ by
$$\Theta(q) = q + \frac{1-q}{1+q} f \left( \frac{q}{(1+q)^2} \right).$$

Banica and Bisch then introduced their $T$ series, which is defined for any Dynkin diagram (and affine Dynkin diagram) by
\begin{equation} \label{def:T-G}
T(q) = \frac{2 G(q^{1/2})-1}{1-q},
\end{equation}
in order to compute the spectral measures for the Dynkin diagrams (and affine Dynkin diagrams) of type $E$. In terms of the Hilbert series $f$, we have
$$T(q) = \frac{\Theta(q)-q}{1-q} = \frac{1}{1+q} f \left( \frac{q}{(1+q)^2} \right).$$

We can define a generalized $T$ series $\widetilde{T}_{ij}$ by
\begin{equation} \label{SU(N)_general_T-series}
\widetilde{T}(q) = \frac{1}{1+q} \widetilde{f} \left( \frac{q}{(1+q)^2} \right),
\end{equation}
where the matrix $\widetilde{f} (z) = \left( 1 - z^{\frac{1}{2}} \Delta_X \right)^{-1}$, and $[\widetilde{f} (z)]_{ij}$ counts paths from $i$ to $j$. Then $f(z) = \varphi(\widetilde{f} (z))$ and $T(q) = \varphi(\widetilde{T}(q))$.

The $T$ series for the $ADE$ Dynkin diagrams and their affine versions (except for $D_{n}^{(1)}$) were computed in \cite{banica/bisch:2007}.
These expressions can be easily derived from the spectral measures computed above for these graphs, since the $T$ series is additive with respect to the underlying measures; that is, if the measure $\varepsilon$ can be written as $\varepsilon = \alpha_1 \varepsilon_1 + \cdots + \alpha_s \varepsilon_s$ for some $s \in \mathbb{N}$, where $\sum_i \alpha_i = 1$, then the $T$ series $T_{\varepsilon}$ for $\varepsilon$ is $T_{\varepsilon} = \alpha_1 T_{\varepsilon_1} + \cdots + \alpha_s T_{\varepsilon_s}$.
The $T$ series for the measures $\mathrm{d}_n$, $\alpha \mathrm{d}_n$, $\mathrm{d}_n' = 2\mathrm{d}_{2n} - \mathrm{d}_n$, $\alpha \mathrm{d}_n'$ are easily computed from (\ref{def:T-G}) and using
$$\int_{\mathbb{T}} \frac{u^{-m}}{1-qu} \mathrm{d}_n u = \int_{\mathbb{T}} \sum_{j=0}^{\infty} q^j u^{j-m} \mathrm{d}_n u = \sum_{k=0}^{\infty} q^{2kn+r} = \frac{q^r}{1-q^{2n}},$$
where $m=2ln+r$ for $l \in \mathbb{Z}$, $r \in \{ 0,1,\ldots, 2n-1 \}$ (see \cite[Lemma 6.1]{banica/bisch:2007}).
Let $T^{\mathcal{G}}$ denote the $T$ series for the graph $\mathcal{G}$. Then the $T$ series are given by:
\begin{center}
\begin{minipage}[t]{14cm}
 \begin{minipage}[t]{7cm}
\begin{eqnarray*}
T^{A_n} & = & \frac{1-q^{n}}{1-q^{n+1}}, \\
T^{D_n} & = & \frac{1+q^{n-3}}{1+q^{n-2}}, \\
T^{E_6} & = & \frac{(1-q^6)(1-q^8)}{(1-q^3)(1-q^{12})}, \\
T^{E_7} & = & \frac{(1-q^9)(1-q^{12})}{(1-q^4)(1-q^{18})}, \\
T^{E_8} & = & \frac{(1-q^{10})(1-q^{15})(1-q^{18})}{(1-q^5)(1-q^9)(1-q^{30})},
\end{eqnarray*}
 \end{minipage}
 \hfill
 \begin{minipage}[t]{6cm}
\begin{eqnarray*}
T^{A_{2n}^{(1)}} & = & \frac{1+q^n}{(1-q)(1-q^n)}, \\
T^{D_n^{(1)}} & = & \frac{1+q^{n-1}}{(1-q^2)(1-q^{n-2})}, \\
T^{E_6^{(1)}} & = & \frac{1+q^6}{(1-q^3)(1-q^4)}, \\
T^{E_7^{(1)}} & = & \frac{1+q^9}{(1-q^4)(1-q^6)}, \\
T^{E_8^{(1)}} & = & \frac{1+q^{15}}{(1-q^6)(1-q^{10})}.
\end{eqnarray*}
 \end{minipage}
\end{minipage}
\end{center}

\subsection{Kostant Polynomial}

We now introduce a polynomial for finite subgroups of $SU(2)$ which is related to the $T$-series defined in Section \ref{Sect:T-series}. The precise relation between the two polynomials will be given later in Theorem \ref{Thm:SU(2)polynomials}.
For a subgroup $\Gamma \subset SU(2)$ and an irreducible representation $\gamma$ of $\Gamma$, the Kostant polynomial $F_{\gamma}$ counts the multiplicity of $\gamma$ in $(j)$, the $j+1$-dimensional irreducible representation of $SU(2)$ restricted to $\Gamma$. The Kostant polynomial $F_{\gamma}$ is given by
$$F_{\gamma}(t) = \sum_{j=0}^{\infty} \langle (j), \gamma, \rangle_{\Gamma} \; t^j,$$
where $\langle (j), \gamma \rangle_{\Gamma}$ counts the multiplicity of $\gamma$ in $(j)$. Let $F(t) = \sum_{j=0}^{\infty} t^j (j) = \sum_{\gamma} F_{\gamma}(t) \gamma$. Then we obtain the recursion formulae
\begin{eqnarray*}
F(t) \otimes (1) & = & \sum_{\gamma} F_{\gamma}(t) \gamma \otimes (1) \;\; = \;\; \sum_{j=0}^{\infty} t^j (j) \otimes (1) \\
& = & \sum_{j=0}^{\infty} t^j ((j-1) \oplus (j+1)) \;\; = \;\; (t^{-1} + t) F(t) - \frac{\mathrm{id}}{t},
\end{eqnarray*}
where $\mathrm{id}$ is the identity representation of $\Gamma$.
Evaluating this polynomial by taking its character on conjugation classes $\Gamma_i$ of $\Gamma$ we obtain \cite{itzykson:1989}:
\begin{equation} \label{eqn:SU(2)_Kostant-characters}
F_{\gamma}(t) = \sum_i \frac{|\Gamma_i|}{|\Gamma|}\frac{\chi_{\gamma}^{\ast} (\Gamma_i)}{1 - t \chi_{\rho} (\Gamma_i) + t^2}.
\end{equation}
The explicit result was worked out by Kostant in \cite{kostant:1984}, where he showed that the polynomials $F_{\gamma}(t)$ have the simple form
\begin{equation} \label{SU(2)_Kostant}
F_{\gamma}(t) = \frac{z_{\gamma}(t)}{(1-t^a)(1-t^b)},
\end{equation}
where $a,b$ are positive integers which satisfy $a+b=h+2$ and $ab=2 |\Gamma|$, where $h$ is the Coxeter number of the Dynkin diagram $\mathcal{G}$, and $z_{\gamma}(t)$ is now a finite polynomial. The values of $a,b$ are:

\begin{center}
\begin{tabular}{|c|c|c|} \hline
Dynkin Diagram $\mathcal{G}$ & $h$ & $a,b$ \\
\hline\hline $A_l$ & $l+1$ & 2, $l+1$ \\
\hline $D_l$ & $2l-2$ & 4, $2l-4$ \\
\hline $E_6$ & 12 & 6, 8 \\
\hline $E_7$ & 18 & 8, 12 \\
\hline $E_8$ & 30 & 12, 20 \\
\hline
\end{tabular}
\end{center}

The Kostant polynomial is related to subfactors realizing the $ADE$ modular invariants in \cite[$\S$3.3]{evans:2003}. Let $\ast$ label the trivial representation of $\Gamma$. By the argument of changing the $\iota$-vertex \cite{evans:2002} it may be assumed that the subfactor $N \subset M$ realizing the $ADE$ modular invariant has the $\iota$-vertex on the vertex which would join the extended vertex $\ast$ of the affine Dynkin diagram $\mathcal{G}_{\Gamma}$. For all $DE$ cases there is a natural bijection between (equivalence classes of) non-trivial irreducible representations of $\Gamma$ and $M$-$N$ sectors $[\iota \lambda_l]$, since the irreducible representations label the vertices of the $DE$ graph, as do the sectors $[\iota \lambda_l]$. Let $\rho$ denote the fundamental representation of $\Gamma$. Denoting the $M$-$N$ morphism associated to the irreducible representation $\gamma \neq \ast$ by $\overline{a}_{\gamma}$ (so $\iota = \overline{a}_{\rho}$), it was shown in \cite{evans:2003} that the polynomials $p_{\gamma}$ defined by
$$p_{\ast}(t) = 1 + q^{k+2}, \qquad \qquad
p_{\gamma}(t) = \sum_{i=0}^{k} \langle \overline{a}_{\gamma}, \iota \lambda_j \rangle t^{j+1},$$
are equal to the numerators $z_{\gamma}(t)$ of the Kostant polynomial $F_{\gamma}(t)$, and consequently $F_{\gamma}(t) = p_{\gamma}(t) / \Omega(t)$, where $\Omega(t) = (1 + t^2) p_{\ast}(t) - t p_{\rho}(t)$.
The Kostant polynomial $F_{\ast}(t)$ for the graphs $E_n$, $n=6,7,8$, is in fact just the $T$-series $T^{E_n^{(1)}}(t^2)$ of Section \ref{Sect:T-series}. This is because the generating series $G(q)$ of the moments of the spectral measures for $E_n^{(1)}$, $n=6,7,8$ is essentially equal to the Kostant polynomial for $E_n$, cf. (\ref{eqn:SU(2)_Kostant-characters}) and (\ref{eqn:Generating_series_momentsSU(2)}). More precisely, $F_{\ast}(t) = (1+t^2)^{-1} G \left( t/(1+t^2) \right) = T(t^2)$ (see also Theorem \ref{Thm:SU(2)polynomials} (iii)).

\subsection{Molien Series} \label{Sect:SU(2)-Molien}

Another related polynomial is the Molien series, which for subgroups of $SU(2)$ is in fact equal to the Kostant polynomial. Let $\Gamma$ be a finite subgroup of $SU(N)$ as above. For $i = 0,1,\ldots \;$, let $M_i$ be a representation of $\Gamma$ with $\textrm{dim} \; M_i < \infty$, and let $M = \bigoplus_{i=0}^{\infty} M_i$. With $\gamma$ an irreducible representation of $\Gamma$, the Molien series $P_{M, \gamma}$ of $M$ is defined in \cite{gomi/nakamura/shinoda:2004} by
$$P_{M, \gamma} (t) = \sum_{i=0}^{\infty} \langle M_i, \gamma \rangle_{\Gamma} \; t^i,$$
and counts the multiplicity $\langle M_i, \gamma \rangle_{\Gamma}$ of $\gamma$ in $M_i$.

Let $\overline{\mathbb{C}^N}$ denote the dual vector space of $\mathbb{C}^N$, and denote by $S = \bigoplus_k S^k(\overline{\mathbb{C}^N})$ the symmetric algebra of $\overline{\mathbb{C}^N}$ over $\mathbb{C}$, where $S^k(\overline{\mathbb{C}^N})$ is the $k^{\mathrm{th}}$ symmetric product of $\overline{\mathbb{C}^N}$. Let $\rho$ be the fundamental representation of $\Gamma$ and $\overline{\rho}$ its conjugate representation, let $\{ \rho_0 = \mathrm{id}, \rho_1 = \rho, \rho_2, \ldots, \rho_s \}$ be the irreducible representations of $\Gamma$ and $\chi_j$ be the character of $\rho_j$ for $j=0,1,\ldots,s$. Then we have Molien's formula for $P_{S, \gamma_j} (t)$ given as \cite{gomi/nakamura/shinoda:2004}:
$$P_{S, \rho_j} (t) = \frac{1}{|\Gamma|} \sum_{g \in \Gamma} \frac{\chi_j^{\ast} (g)}{\textrm{det}(1 - \overline{\rho} (g) t)}.$$

Let $R_k$ denote the sum of all the representations of $SU(N)$ which have Dynkin labels $\lambda_1, \lambda_2, \ldots , \lambda_{(N-1)}$ such that $\lambda_1 + \cdots + \lambda_{(N-1)} = k$, and $R = \bigoplus_{k=0}^{\infty} R_k$. Then in this notation, $P_{R, \gamma}$ recovers the Kostant polynomial $F_{\gamma}$, where $\gamma$ is an irreducible representation of $\Gamma$:
\begin{equation} \label{eqn:gen_Molien=Kostant}
P_{R, \gamma}(t) = \sum_{i=0}^{\infty} \langle R_i, \gamma \rangle_{\Gamma} \; t^i = F_{\gamma} (t,t,\ldots,t).
\end{equation}

Since there is only one Dynkin label $\lambda$ for any representation of $SU(2)$, $R_k = (k)$, the $(k+1)$-dimensional representation of $SU(2)$, for each $k$. Then by $(\ref{eqn:gen_Molien=Kostant})$ the Molien series $P_{R,\gamma}(t)$ for a subgroup $\Gamma \subset SU(2)$ is equal to the Kostant polynomial $F_{\gamma}(t)$. The $k^{\textrm{th}}$ symmetric product of $\overline{\mathbb{C}^2}$ gives the irreducible level $k$ representation, so that $R = S$ for $SU(2)$, and $P_{S,\gamma}(t) = F_{\gamma}(t)$.

\subsection{Hilbert Series of Pre-projective Algebras} \label{Sect:SU(2)-Hilbert_Series}

Finally, we introduce another related polynomial, the Hilbert series $H(t)$, which counts the dimensions of pre-projective algebras for the $ADE$ and affine Dynkin diagrams.
Let $\mathcal{G}$ be any (oriented or unoriented) graph, and let $\mathbb{C}\mathcal{G}$ be the algebra with basis given by the paths in $\mathcal{G}$, where paths may begin at any vertex of $\mathcal{G}$. Multiplication of two paths $a$, $b$ is given by concatenation of paths $a \cdot b$ (or simply $ab$), where $ab$ is defined to be zero if $r(a) \neq s(b)$.
Note that the algebra $\mathbb{C}\mathcal{G}$ is not the path algebra $A(\mathcal{G})$ for $\mathcal{G}$ in the usual operator algebraic meaning.
Let $[\mathbb{C}\mathcal{G}, \mathbb{C}\mathcal{G}]$ denote the subspace of $\mathbb{C}\mathcal{G}$ spanned by all commutators of the form $ab - ba$, for $a,b \in \mathbb{C}\mathcal{G}$. If $a,b$ are paths in $\mathbb{C}\mathcal{G}$ such that $r(a) = s(b)$ but $r(b) \neq s(a)$, then $a b - b a = a b$, so in the quotient $\mathbb{C}\mathcal{G} / [\mathbb{C}\mathcal{G}, \mathbb{C}\mathcal{G}]$ the path $a b$ will be zero. Then any non-cyclic path, i.e. any path $a$ such that $r(a) \neq s(a)$, will be zero in $\mathbb{C}\mathcal{G} / [\mathbb{C}\mathcal{G}, \mathbb{C}\mathcal{G}]$. If $a = a_1 a_2 \cdots a_k$ is a cyclic path in $\mathbb{C}\mathcal{G}$, then $a_1 a_2 \cdots a_k - a_k a_1 \cdots a_{k-1} = 0$ in $\mathbb{C}\mathcal{G} / [\mathbb{C}\mathcal{G}, \mathbb{C}\mathcal{G}]$, so $a_1 a_2 \cdots a_k$ is identified with $a_k a_1 \cdots a_{k-1}$. Similarly, $a = a_1 a_2 \cdots a_k$ is identified with every cyclic permutation of the edges $a_j$, $j=1,\ldots,k$. So the commutator quotient $\mathbb{C}\mathcal{G} / [\mathbb{C}\mathcal{G}, \mathbb{C}\mathcal{G}]$ may be identified, up to cyclic permutation of the arrows, with the vector space spanned by cyclic paths in $\mathcal{G}$.

The pre-projective algebra $\Pi$ of a finite unoriented graph $\mathcal{G}$ is defined as the quotient of $\mathbb{C} \mathcal{G}$ by the two-sided ideal generated by $\theta = \sum_{i,\sigma} \theta_i^{\sigma}$, where the summation is over all vertices $i$ and edges $\sigma$ of $\mathcal{G}$ such that $i$ is an endpoint for $\sigma$, and $\theta_i^{\sigma} \in \mathbb{C} \mathcal{G}$ is defined to be the loop of length two starting and ending at vertex $i$ formed by going along the edge $\sigma$ and back again. So the pre-projective algebra is the quotient algebra under relations $\theta$, and any closed loop of length 2 on $\mathcal{G}$ is identified with a linear combination of all the other closed loops of length 2 on $\mathcal{G}$ which have the same initial vertex. In the language of planar algebras for bipartite graphs (see \cite{jones:2000}), this is closely related to taking the (complement of the) kernel of the insertion operators given by the cups and caps.

For a graph $\mathcal{G}$ without any closed loops of length one, i.e. edges from a vertex to itself, the pre-projective algebra $\Pi$ has the following description as a quotient of a path algebra by a two-sided ideal generated by derivatives of a potential $\Phi$. We fix an orientation for the edges of $\mathcal{G}$, and form the double $\overline{\mathcal{G}}$ of $\mathcal{G}$, where for each (oriented) edge $\gamma$ we add the reverse edge $\widetilde{\gamma}$ which has $s(\widetilde{\gamma}) = r(\gamma)$, $r(\widetilde{\gamma}) = s(\gamma)$. We define a potential $\Phi$ by $\Phi = \sum_{\gamma} \gamma \widetilde{\gamma}$, where the summation is over all edges of $\mathcal{G}$. Let $\gamma_1 \gamma_2 \cdots \gamma_k$ be any closed loop of length $k$ in $\mathbb{C}\overline{\mathcal{G}} / [\mathbb{C}\overline{\mathcal{G}}, \mathbb{C}\overline{\mathcal{G}}]$, $k >1$. We define derivatives $\partial_i: \mathbb{C}\overline{\mathcal{G}} / [\mathbb{C}\overline{\mathcal{G}}, \mathbb{C}\overline{\mathcal{G}}] \rightarrow \mathbb{C}\overline{\mathcal{G}}$ for each vertex $i \in \mathfrak{V}_{\mathcal{G}}$ of $\mathcal{G}$ by $\partial_i (\gamma_1 \gamma_2 \cdots \gamma_k) = \sum_j \gamma_j \gamma_{j+1} \cdots \gamma_k \gamma_1 \cdots \gamma_{j-1}$, where the summation is over all $1 \leq j \leq k$ such that $s(\gamma_j) = i$. Then on paths $\gamma \widetilde{\gamma} \in \mathbb{C}\overline{\mathcal{G}} / [\mathbb{C}\overline{\mathcal{G}}, \mathbb{C}\overline{\mathcal{G}}]$, we have
$$\partial_i (\gamma \widetilde{\gamma}) = \left\{ \begin{array}{cl}
                                \gamma \widetilde{\gamma} & \textrm{ if } s(\gamma) = i, \\
                                \widetilde{\gamma} \gamma & \textrm{ if } r(\gamma) = i, \\
                                0 & \textrm{ otherwise.}
                                \end{array} \right.,$$
and $\Pi \cong \mathbb{C}\overline{\mathcal{G}} / (\partial_i \Phi: i \in \mathfrak{V}_{\mathcal{G}})$.
For any graph $\mathcal{G}$ and potential $\Phi$, Bocklandt \cite[Theorem 3.2]{bocklandt:2008} showed that if $A (\mathbb{C} \mathcal{G}, \Phi)$ is Calabi-Yau of dimension 2 then $A (\mathbb{C} \mathcal{G}, \Phi)$ is the pre-projective algebra of a non-Dynkin quiver.

We can define the Hilbert series for $A(\mathbb{C}\mathcal{G}, \Phi)$ as $H_A(t) = \sum_{k=0}^{\infty} H_{ji}^k t^k$, where the $H_{ji}^k$ are matrices which count the dimension of the subspace $\{ i \cdot a \cdot j | a \in A(\mathbb{C}\mathcal{G}, \Phi)_k \}$, where $A(\mathbb{C}\mathcal{G}, \Phi)_k$ is the subspace of $A(\mathbb{C}\mathcal{G}, \Phi)$ of all paths of length $k$, and $i$, $j$ are paths in $A(\mathbb{C}\mathcal{G}, \Phi)_0$, corresponding to vertices of $\mathcal{G}$.

Let $q \in \mathbb{C} \setminus \{ 0 \}$.
If $q = \pm 1$ or $q$ not a root of unity, the tensor category $\mathcal{C}_q$ of representations of the quantum group $SU(2)_q$ has a complete set $\{ L_s \}_{s = 0}^{\infty}$ of simple objects.
If $q$ is an $n^{\textrm{th}}$ root of unity, $\mathcal{C}_q$ is the semisimple subquotient of the category of representations of $SU(2)_q$. In this case, the set $\{ L_s \}_{s = 0}^{h(q)-2}$ is the complete set of simple objects of $\mathcal{C}_q$, where $L_s$ is the deformation of the $(s+1)$-dimensional representation of $SU(2)$, and $h(q)$ is $n$ when $n$ is odd and $n/2$ when $n$ is even,
satisfying:
\begin{equation} \label{eqn:C_q-FusionRules2}
L_r \otimes L_s \simeq  \bigoplus_{\stackrel{t=|r-s|}{\scriptscriptstyle{t \equiv r+s \mathrm{ mod } 2}}}^{k} L_t,
\end{equation}
where
$$k = \left\{ \begin{array}{cl}
                                r+s & \textrm{ if } r+s<h(q)-1, \\
                                2h(q)-4-r-s & \textrm{ if } r+s \geq h(q)-1. \end{array} \right.$$

Semisimple module categories over $\mathcal{C}_q$ where classified in \cite{etingof/ostrik:2004}. A semisimple $\mathcal{C}_q$-module category $\mathcal{D}$ is abelian, and is equivalent to the category of $I$-graded vector spaces $\mathcal{M}_I$, where $I$ are simple objects of $\mathcal{D}$. The structure of a $\mathcal{C}_q$ category on $\mathcal{M}_I$ is the same as a tensor functor $F$ from $\mathcal{C}_q$ to $\mathrm{Fun}(\mathcal{M}_I,\mathcal{M}_I) \cong \mathcal{M}_{I \times I}$, the category of additive functors from $\mathcal{M}_I$ to itself. When $q = \pm 1$ or $q$ is not a root of unity, by \cite[Theorem 2.5]{etingof/ostrik:2004}, such functors are classified by the following data:
\begin{itemize}
\item a collection of finite dimensional vector spaces $V_{ij}$, $i,j \in I$,
\item a collection of non-degenerate bilinear forms $E_{ij}:V_{ij} \otimes V_{ji} \rightarrow \mathbb{C}$, subject to the condition, $\sum_j \mathrm{Tr}(E_{ij} (E_{ji}^T)^{-1}) = -q-q^{-1}$, for each $i \in I$.
\end{itemize}
When $q$ is a root of unity there is an extra condition given in \cite{etingof/ostrik:2004}, due to the fact that $\mathcal{C}_q$ is now a quotient of the tensor category whose objects are $V^{\otimes m}$, $m \in \mathbb{N}$.

Let $\Delta$ be the matrix given by $\Delta_{i,j} = \mathrm{dim}V_{ij}$.
Quantum McKay correspondence gives a graph with adjacency matrix $\Delta$ and vertex set $I$. The free algebra $T$ in $\mathcal{C}_q$ generated by the self-dual object $V=L_1$ maps to the path algebra of the McKay graph under the functor $F:\mathcal{C}_q \rightarrow \mathcal{M}_{I \times I}$. Let $S$ be the quotient of $T$ by the two-sided ideal $J$ generated by the image of $\mathbf{1} = L_0$ under the map $\mathbf{1} \stackrel{\mathrm{coev}_V}{\rightarrow} V \otimes \overline{V} \stackrel{\mathrm{id}_V \otimes \phi^{-1}}{\rightarrow} V \otimes V$, where $\phi$ is any choice of isomorphism from $V$ to its conjugate representation $\overline{V}$. In the classical situation, $q = 1$, $S$ is the algebra of polynomials in two commuting variables. More generally, $S$ is called the $q$-symmetric algebra, or the algebra of functions on the quantum plane. The structure of these algebras is well known, see for example \cite{kassel:1995}.
Applying the functor $F$ to $S$ gives an algebra $\widetilde{\Pi} = F(S)$ which is the quotient of the path algebra with respect to the two-sided ideal $F(J)$. Then given any arbitrary connected graph $\mathcal{G}$, there exists a particular value of $q$ and choice of $\mathcal{C}_q$-module category $\mathcal{D}$ such that $\widetilde{\Pi}$ is equal to the pre-projective algebra $\Pi$ of $\mathcal{G}$ \cite[Lemma 2.2]{malkin/ostrik/vybornov:2006}.

When $q$ is not a root of unity, the $m^{\textrm{th}}$ graded component of the $q$-symmetric algebra $S$ is given by $S(m) = L_m$, for $m \in \mathbb{N}$, which satisfies
\begin{equation} \label{eqn:fusion_recursion_for_Hilbert}
L_1 \otimes L_m \simeq L_{m-1} \oplus L_{m+1}.
\end{equation}
Then summing (\ref{eqn:fusion_recursion_for_Hilbert}) over all $m \in \mathbb{N}$, with a grading $t^m$, gives $t L_1 \otimes S = t^2 S \oplus S \ominus L_0$. Applying the functor $F$ one obtains a recursion $t \Delta H(t) = H(t) + t^2 H(t) - 1$, where $\Delta$ is the adjacency matrix of the (quantum) McKay graph $\mathcal{G}$. Then we obtain the following result \cite[Theorem 2.3a]{malkin/ostrik/vybornov:2006}:
\begin{equation} \label{SU(2)_Hilbert(affine)}
H(t) = \frac{1}{1 - \Delta t + t^2}.
\end{equation}
For an $ADET$ graph $\mathcal{G}$, $q$ is an $n^{\textrm{th}}$ root of unity, and $h(q)=h$ is the Coxeter number of $\mathcal{G}$. The $m^{\textrm{th}}$ graded component is given by $S(m) = L_m$ for $0 \leq m \leq h-2$, and $S(m) = 0$ for $m \geq h-1$. Defining $\widehat{S} = S \ominus t^h(L_{h-2} \otimes S) \oplus t^{2h}(L_{h-2} \otimes L_{h-2} \otimes S) \ominus \cdots \;$, the fusion rules (\ref{eqn:C_q-FusionRules2}) give the recursion $L_1 \otimes \widehat{S}(m) \simeq \widehat{S}(m-1) \oplus \widehat{S}(m+1)$. Applying the functor $F$ gives $1 + t^h P + t \Delta H(t) = H(t) + t^2 H(t)$, where the matrix $P = F(L_{h-2})$. Then for the Dynkin diagrams (and the graph $\mathrm{Tad}_n$), there is a `correction' term in the numerator, so that \cite[Theorem 2.3b]{malkin/ostrik/vybornov:2006}:
\begin{equation}
H(t) = \frac{1 + Pt^h}{1 - \Delta t + t^2},
\end{equation}
where $P$ is a permutation corresponding to some involution of the vertices of the graph. Since $L_{h-2} \otimes L_{h-2} \simeq L_0$, $P^2 = F(L_{h-2} \otimes L_{h-2}) = F(\mathbf{1})$ so $P^2$ is the identity matrix. The matrix $P$ is an automorphism of the underlying graph \cite{malkin/ostrik/vybornov:2006}; for $A_n$, $D_{2n+1}$, $E_6$ it is the unique nontrivial involution, while for $D_{2n}$, $E_7$, $E_8$ (and $\mathrm{Tad}_n$) it is the identity matrix, i.e. the matrix $P$ corresponds to the Nakayama permutation $\pi$ for the $ADE$ graph \cite{erdmann/snashall:1998}. A Nakayama automorphism of $\Pi$ is an automorphism $\nu$ of edges for which there exists an element $\widehat{b}$ of the dual $\Pi^{\ast}$ of $\Pi$ such that $\widehat{b} a = \nu(a) \widehat{b}$ for all $a \in \Pi$. The Nakayama automorphism is related to the Nakayama permutation by $\nu(a) = \epsilon(a) \pi(a)$ for all edges $a$ of the Dynkin quiver, where $\epsilon(a) \in \{ \pm 1 \}$.

We now present the following result which relates these various polynomials:

\begin{theorem} \label{Thm:SU(2)polynomials}
Let $\Gamma$ be a finite subgroup of $SU(2)$ so that $\mathcal{G}_{\Gamma}$ is one of the affine Dynkin diagrams, with the vertices of $\mathcal{G}_{\Gamma}$ labelled by the irreducible representations $\gamma$ of $\Gamma$, with the distinguished vertex $\ast$ labelled by $\mathrm{id}$. Let $G(q)$ be the generating series of the moments for finite subgroups of $SU(2)$ in (\ref{eqn:Generating_series_momentsSU(2)}), $\widetilde{T}$ be the generalized T series defined in Section \ref{Sect:T-series}, and let $P_{\gamma}$, $F_{\gamma}$ be the Molien series, Kostant polynomial respectively of $\Gamma$. Then for the Hilbert series $H$ of $\mathcal{G}_{\Gamma}$ as in (\ref{SU(2)_Hilbert(affine)}), the following hold:
\begin{itemize}
\item[(i)] $\widetilde{T} (t^2) = H(t)$,

\item[(ii)] $H_{{\gamma},\mathrm{id}} (t) = P_{\gamma} (t) = F_{\gamma} (t)$,

\item[(iii)] $T(t^2) = H_{\mathrm{id},\mathrm{id}} (t) = P_{\mathrm{id}} (t) = F_{\mathrm{id}} (t) =  \frac{1}{1+t^2} G \left( \frac{t}{1+t^2} \right)$.
\end{itemize}
\end{theorem}
\emph{Proof:}
\begin{itemize}
\item[(i)]
From (\ref{SU(N)_general_T-series}) we have
\begin{eqnarray*}
\widetilde{T} (t^2) & = & \frac{1}{1 + t^2} \; \widetilde{f} \left( \frac{t^2}{(1 + t^2)^2} \right) = \frac{1}{1 + t^2} \cdot \frac{1}{1- t(1 + t^2)^{-1} \Delta} = \frac{1}{1 + t^2 - t \Delta} \\
& = & H(t).
\end{eqnarray*}

\item[(ii)]
By \cite[Cor. 2.4 (ii)]{gomi/nakamura/shinoda:2004}, for the symmetric algebra $S=S(\overline{\mathbb{C}^2})$, $P_{\gamma_j} = P_{S, \gamma_j}$ satisfies
$$\sum_{j=0}^{s} [\Delta_{\Gamma}]_{ij} P_{\gamma_j} (t) = (t + t^{-1}) P_{\gamma_i} (t) - t^{-1} \delta_{i,0},$$
where $\gamma_1, \ldots, \gamma_s$ are the irreducible representations associated with the vertices $1, \ldots, s$ of $\mathcal{G}_{\Gamma}$. Then
multiplying through by $t$ we obtain
$$\sum_{j=0}^s \left[ \mathbf{1} - \Delta_{\Gamma} t + \mathbf{1} t^2 \right]_{ij} P_{S, \gamma_j} (t) = \delta_{i,0}.$$
From (\ref{SU(2)_Hilbert(affine)}) we see that the matrix $\left( \mathbf{1} - \Delta_{\Gamma} t + \mathbf{1} t^2 \right)$ is invertible, and hence by the definition of matrix multiplication, we see that
$$P_{\gamma} (t) = \left[ \left( \mathbf{1} - \Delta_{\Gamma} t + \mathbf{1} t^2 \right)^{-1} \right]_{\gamma,\mathrm{id}},$$
which is the first equality. The second was shown in Section \ref{Sect:SU(2)-Molien}.

\item[(iii)]
The first equality follows from $T(q) = \varphi(\widetilde{T}(q))$, and the next two are immediate from $(ii)$. For the last equality, using (\ref{eqn:SU(2)_Kostant-characters}) we have
\begin{eqnarray*}
F_{\mathrm{id}}(t) & = & \sum_{j=1}^n \frac{|\Gamma_j|}{|\Gamma|} \frac{\chi_0^{\ast}(\Gamma_j)}{1 - t \chi_{\rho}(\Gamma_j) + t^2} \;\; = \;\; \frac{1}{1+t^2} \sum_{j=1}^n \frac{|\Gamma_j|}{|\Gamma|} \frac{1}{1 - \left( \frac{t}{1+t^2} \right) \chi_{\rho}(\Gamma_j)} \\
& = & \frac{1}{1+t^2} G \left( \frac{t}{1+t^2} \right).
\end{eqnarray*}
\end{itemize}

\section{$SU(3)$ Case}

We will now consider the case of $SU(3)$. We no longer have self-adjoint operators, but are in the more general setting of normal operators, whose moments are given by (\ref{eqn:moments_normal_operator}).
We will first consider the fixed point algebra of $\bigotimes_{\mathbb{N}} M_3$ under the action of the group $\mathbb{T}^2$ to obtain the spectral measure for the infinite graph which we call $\mathcal{A}^{(6\infty)}$. We will then generalize the method presented in Section \ref{sect:spec_measure-ADE} to the case of $SU(3)$ graphs.

\subsection{Spectral measure for $\mathcal{A}^{(6\infty)}$}

We first consider the fixed point algebra of $\bigotimes_{\mathbb{N}} M_3$ under the action of the group $\mathbb{T}^2$. Let $\rho$ be the fundamental representation of $SU(3)$, so that the restriction of $\rho$ to $\mathbb{T}^2$ is given by
\begin{equation} \label{eqn:restrict_rho_to_T2}
(\rho|_{\mathbb{T}^2})(\omega_1,\omega_2) = \left( \begin{array}{ccc} \omega_1 & 0 & 0 \\ 0 & \omega_2^{-1} & 0 \\ 0 & 0 & \omega_1^{-1}\omega_2 \end{array} \right),
\end{equation}
for $(\omega_1,\omega_2) \in \mathbb{T}^2$.

Let $\{ \chi_{(\lambda_1,\lambda_2)} \}_{\lambda_1,\lambda_2 \in \mathbb{N}}$, $\{ \sigma_{(\lambda_1,\lambda_2)} \}_{\lambda_1,\lambda_2 \in \mathbb{Z}}$ be the irreducible characters of $SU(3)$, $\mathbb{T}^2$ respectively, where if $\chi_{(\lambda_1,\lambda_2)}$ is the character of a representation $\pi$ then $\chi_{(\lambda_2,\lambda_1)}$ is the character of the conjugate representation $\overline{\pi}$ of $\pi$. The trivial character of $SU(3)$ is $\chi_{(0,0)}$, $\chi_{(1,0)}$ is the character of $\rho$, and $\sigma_{(\lambda_1,\lambda_2)}(p,q) = (p^{\lambda_1},q^{\lambda_2})$, for $\lambda_1,\lambda_2 \in \mathbb{Z}$. If $\sigma$ is the restriction of $\chi_{(1,0)}$ to $\mathbb{T}^2$, we have $\sigma = \sigma_{(1,0)} + \sigma_{(0,-1)} + \sigma_{(-1,1)}$ (by (\ref{eqn:restrict_rho_to_T2})), and $\sigma \sigma_{(\lambda_1,\lambda_2)} = \sigma_{(\lambda_1+1,\lambda_2)} + \sigma_{(\lambda_1,\lambda_2-1)} + \sigma_{(\lambda_1-1,\lambda_2+1)}$, for any $\lambda_1,\lambda_2 \in \mathbb{Z}$. So the representation graph of $\mathbb{T}^2$ is identified with the infinite graph $\mathcal{A}^{(6\infty)}$, illustrated in Figure \ref{fig:A^(infty,infty)}, whose vertices are labelled by pairs $(\lambda_1,\lambda_2) \in \mathbb{Z}^2$, and which has an edge from vertex $(\lambda_1,\lambda_2)$ to the vertices $(\lambda_1+1,\lambda_2)$, $(\lambda_1,\lambda_2-1)$ and $(\lambda_1-1,\lambda_2+1)$. The 6 in the notation $\mathcal{A}^{(6\infty)}$ is to indicate that for this graph we are taking six infinities, one in each of the directions of $\pm e_i$, $i=1,2,3$, for the vectors $e_i$ given by $e_1 = \Lambda_1$, $e_2 = \Lambda_2 - \Lambda_1$, $e_3 = - \Lambda_2$, where $\Lambda_1$, $\Lambda_2$ are the fundamental weights of $SU(3)$. We choose the distinguished vertex to be $\ast = (0,0)$. Hence $(\bigotimes_{\mathbb{N}} M_3)^{\mathbb{T}^2} \cong A(\mathcal{A}^{(6\infty)})$.

\begin{figure}[tb]
\begin{center}
  \includegraphics[width=55mm]{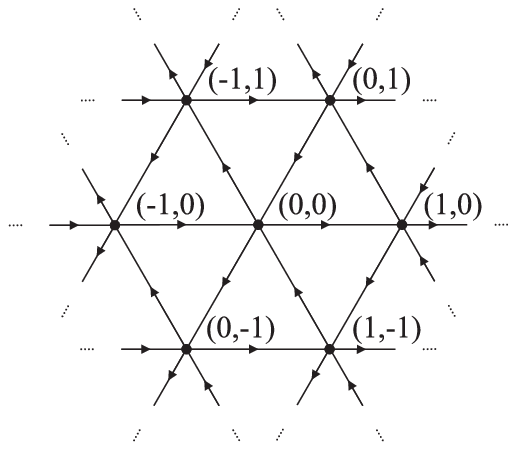}\\
 \caption{The infinite graph $\mathcal{A}^{(6\infty)}$.} \label{fig:A^(infty,infty)}
\end{center}
\end{figure}

We define a normal operator $v_Z$ in $\ell^2(\mathbb{Z}) \otimes \ell^2(\mathbb{Z})$ by $v_Z = s \otimes 1 + 1 \otimes s^{-1} + s^{-1} \otimes s$, where $s$ is again the bilateral shift on $\ell^2(\mathbb{Z})$. Let $\Omega \otimes \Omega$ be the vector $(\delta_{i,0})_i \otimes (\delta_{i,0})_i$. Then $v_Z$ is identified with the adjacency matrix $\Delta$ of $\mathcal{A}^{(6\infty)}$, where we regard the vector $\Omega \otimes \Omega$ as corresponding to the vertex $(0,0)$ of $\mathcal{A}^{(6\infty)}$, and the operators $s \otimes 1$, $s^{-1} \otimes s$, $1 \otimes s^{-1}$ as corresponding to an edge on $\mathcal{A}^{(6\infty)}$, in the direction of the vectors $e_1, e_2, e_3$ respectively. Then $(s^{\lambda_1} \otimes s^{-\lambda_2})(\Omega \otimes \Omega)$ corresponds to the vertex $(\lambda_1, \lambda_2)$ of $\mathcal{A}^{(6\infty)}$, for any $\lambda_1, \lambda_2 \in \mathbb{Z}$, and applying $v_Z^m v_Z^{\ast n} (\Omega \otimes \Omega)$ gives a vector $y=(y_{(\lambda_1,\lambda_2)})$ in $\ell^2(\mathcal{A}^{(6\infty)})$, where $y_{(\lambda_1,\lambda_2)}$ gives the number of paths of length $m+n$ from $(0,0)$ to the vertex $(\lambda_1,\lambda_2)$, where $m$ edges are on $\mathcal{A}^{(6\infty)}$ and $n$ edges are on the reverse graph $\widetilde{\mathcal{A}^{(6\infty)}}$. The relation $(1 \otimes s^{-1})(s^{-1} \otimes s)(s \otimes 1) = s^{-1} s \otimes s^{-1} s = 1 \otimes 1$ corresponds to the fact that traveling along edges in directions $e_1$ followed by $e_2$ and then $e_3$ forms a closed loop, and similarly for any permutations of $1 \otimes s^{-1}$, $s^{-1} \otimes s$, $s \otimes 1$.

Define a state $\varphi$ on $C^{\ast}(v_Z)$ by $\varphi( \, \cdot \, ) = \langle \, \cdot \, (\Omega \otimes \Omega), \Omega \otimes \Omega \rangle$.
When $m \not \equiv n \textrm{ mod } 3$ it is impossible for there to be a closed loop of length $m + n$ beginning and ending at the vertex $(0,0)$, with the first $m$ edges are on $\mathcal{A}^{(6\infty)}$ and the next $n$ edges are on the reverse graph $\widetilde{\mathcal{A}^{(6\infty)}}$. Hence $\varphi(v_Z^m v_Z^{\ast n}) = 0$ for $m \not \equiv n \textrm{ mod } 3$. We use the notation $(a,b,c)!$ to denote the multinomial coefficient $(a+b+c)!/(a!b!c!)$. For $m \equiv n \textrm{ mod } 3$, we have
\begin{eqnarray*}
\varphi(v_Z^m v_Z^{\ast n}) & = & \sum_{\stackrel{0 \leq k_1 + k_2 \leq m}{\scriptscriptstyle{0 \leq l_1 + l_2 \leq n}}} (k_1, k_2, m - k_1 - k_2)! (l_1, l_2, n - l_1 - l_2)! \; \varphi(s^{r_1} \otimes s^{r_2}) \\
& = & \sum_{\stackrel{0 \leq k_1 + k_2 \leq m}{\scriptscriptstyle{0 \leq l_1 + l_2 \leq n}}} (k_1, k_2, m - k_1 - k_2)! (l_1, l_2, n - l_1 - l_2)! \; \delta_{r_1, 0} \; \delta_{r_2, 0},
\end{eqnarray*}
where
\begin{equation} \label{eqn:r1,r2}
r_1 = 2k_1 + k_2 - 2l_1 - l_2 + n - m, \qquad r_2 = 2l_2 + l_1 - 2k_2 - k_1 + m - n.
\end{equation}
Then we get a non-zero contribution when $l_1 = k_1 + r$, $l_2 = k_2 + r$, where $n = m + 3r$, $r \in \mathbb{Z}$. So we obtain
\begin{equation} \label{eqn:momentsA(infty)6}
\varphi(v_Z^m v_Z^{\ast n}) = \sum_{k_1,k_2} (k_1, k_2, m - k_1 - k_2)! (k_1 + r, k_2 + r, m + r - k_1 - k_2)!
\end{equation}
where the summation is over all integers $k_1, k_2 \geq 0$ such that $\mathrm{max}(0,-r) \leq k_1, k_2 \leq \mathrm{min}(m,m+2r)$ and $k_1 + k_2 \leq \mathrm{min}(m,m+r)$.

\begin{proposition}
The dimension of the $m^{\mathrm{th}}$ level of the path algebra for the infinite graph $\mathcal{A}^{(6\infty)}$ is given by
$$\mathrm{dim}\left( \left(\otimes^m M_3 \right)^{\mathbb{T}^2} \right) = \mathrm{dim}(A(\mathcal{A}^{(6\infty)})_m) = \sum_{j=0}^m C^{2j}_j (C^m_j)^2.$$
\end{proposition}
\emph{Proof:}
When $m = n$ we have
\begin{eqnarray*}
\varphi(v_Z^m v_Z^{\ast m}) & = & \sum_{0 \leq k_1 + k_2 \leq m} ((k_1, k_2, m - k_1 - k_2)!)^2 \\
& = & \sum_{k_1 = 0}^m \sum_{k_2 = 0}^{m-k_1} \left( \frac{m!}{k_1! k_2! (m-k_1-k_2)!} \right)^2 \\
& = & \sum_{k_1 = 0}^m \left( \frac{m!}{k_1! (m-k_1)!} \right)^2 \sum_{k_2 = 0}^{m-k_1} \left( \frac{(m-k_1)!}{k_2! (m-k_1-k_2)!} \right)^2 \\
& = & \sum_{k_1 = 0}^m (C^m_{k_1})^2 \sum_{k_2 = 0}^{m-k_1} (C^{m-k_1}_{k_2})^2 \;\; = \;\; \sum_{k_1 = 0}^m (C^m_{k_1})^2 C^{2(m-k_1)}_{m-k_1}.
\end{eqnarray*}
\hfill
$\Box$

\begin{figure}[tb]
\begin{center}
  \includegraphics[width=70mm]{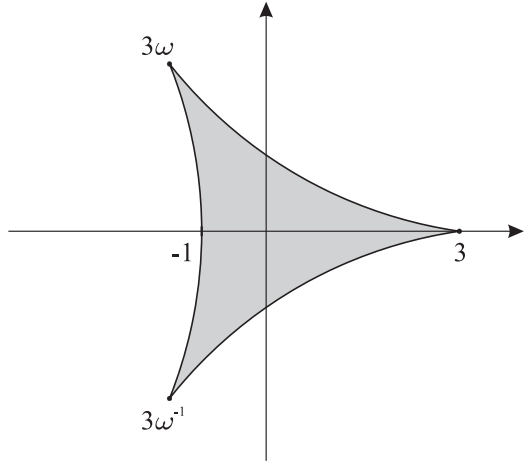}\\
  \caption{The set $\mathfrak{D}$, the closure of the interior of a deltoid.} \label{fig:hypocycloid-S}
\end{center}
\end{figure}

Since the spectrum $\sigma(s)$ of $s$ is $\mathbb{T}$, the spectrum $\sigma(v_Z)$ of $v_Z$ is $\mathfrak{D} = \{ \omega_1 + \omega_2^{-1} + \omega_1^{-1}\omega_2 | \; \omega_1,\omega_2 \in \mathbb{T} \}$, the closure of the interior of the three-cusp hypocycloid, called a deltoid, illustrated in Figure \ref{fig:hypocycloid-S}, where $\omega = e^{2\pi i/3}$.
Any point in $\mathfrak{D}$ can be parameterized by
\begin{equation} \label{eqn:parameterizeS}
x = r(2 \cos(2 \pi t) + \cos (4 \pi t)),  \qquad y = r(2 \sin(2 \pi t) - \sin (4 \pi t)),
\end{equation}
where $0 \leq r \leq 1$, $0 \leq t < 1$, with $r=1$ corresponding to the boundary of $\mathfrak{D}$.

Thus the support of the probability measure $\mu_{v_Z}$ is contained in $\mathfrak{D}$. There is a map $\Phi:\mathbb{T}^2 \rightarrow \mathfrak{D}$ from the torus to $\mathfrak{D}$ given by
\begin{equation} \label{Phi:T^2->S}
\Phi(\omega_1,\omega_2) = \omega_1 + \omega_2^{-1} + \omega_1^{-1} \omega_2,
\end{equation}
where $\omega_1, \omega_2 \in \mathbb{T}$.

Consider the permutation group $S_3$ as the subgroup of $GL(2,\mathbb{Z})$ generated by the matrices $T_2$, $T_3$, of orders 2, 3 respectively, given by
\begin{equation} \label{T1,T2}
T_2 = \left( \begin{array}{cc} 0 & -1 \\ -1 & 0 \end{array} \right), \qquad T_3 = \left( \begin{array}{cc} 0 & -1 \\ 1 & -1 \end{array} \right).
\end{equation}
The action of $S_3$ given by $T(\omega_1,\omega_2) = (\omega_1^{a_{11}} \omega_2^{a_{12}}, \omega_1^{a_{21}} \omega_2^{a_{22}})$, for $T = (a_{ij}) \in S_3$, leaves $\Phi(\omega_1,\omega_2)$ invariant, i.e.
\begin{eqnarray*}
& & \Phi(\omega_1,\omega_2) \;\; = \;\; \Phi(\omega_1^{-1}\omega_2,\omega_1^{-1}) \;\; = \;\; \Phi(\omega_2^{-1},\omega_1\omega_2^{-1}) \\
& = & \Phi(\omega_2^{-1},\omega_1^{-1}) \;\; = \;\; \Phi(\omega_1^{-1}\omega_2,\omega_2) \;\; = \;\; \Phi(\omega_1,\omega_1\omega_2^{-1}).
\end{eqnarray*}

Any $S_3$-invariant probability measure $\varepsilon$ on $\mathbb{T}^2$ produces a probability measure $\mu$ on $\mathfrak{D}$ by
$$\int_{\mathfrak{D}} \psi(z) \mathrm{d}\mu(z) = \int_{\mathbb{T}^2} \psi(\omega_1 + \omega_2^{-1} + \omega_1^{-1} \omega_2) \mathrm{d}\varepsilon(\omega_1,\omega_2),$$
for any continuous function $\psi:\mathfrak{D} \rightarrow \mathbb{C}$, where $\mathrm{d}\varepsilon(\omega_1,\omega_2) = \mathrm{d}\varepsilon(g(\omega_1,\omega_2))$ for all $g \in S_3$.

\begin{theorem}
The spectral measure $\varepsilon(\omega_1,\omega_2)$ (on $\mathbb{T}^2$) for the graph $\mathcal{A}^{(6\infty)}$ is given by the uniform Lebesgue measure $\mathrm{d}\varepsilon(\omega_1,\omega_2) = \mathrm{d}\omega_1 \; \mathrm{d}\omega_2$.
\end{theorem}
\emph{Proof:}
With this measure we have
\begin{eqnarray*}
\lefteqn{ \int_{\mathbb{T}^2} (\omega_1 + \omega_2^{-1} + \omega_1^{-1}\omega_2)^m (\omega_1^{-1} + \omega_2 + \omega_1\omega_2^{-1})^n \mathrm{d}\omega_1 \; \mathrm{d}\omega_2 } \\
& = & \sum_{\stackrel{0 \leq k_1 + k_2 \leq m}{\scriptscriptstyle{0 \leq l_1 + l_2 \leq n}}} \left( (k_1, k_2, m - k_1 - k_2)! (l_1, l_2, n - l_1 - l_2)! \; \int_{\mathbb{T}^2} \omega_1^{r_1} \omega_2^{r_2} \mathrm{d}\omega_1 \; \mathrm{d}\omega_2 \right) \\
& = & \sum_{\stackrel{0 \leq k_1 + k_2 \leq m}{\scriptscriptstyle{0 \leq l_1 + l_2 \leq n}}} (k_1, k_2, m - k_1 - k_2)! (l_1, l_2, n - l_1 - l_2)! \; \delta_{r_1, 0} \; \delta_{r_2, 0},
\end{eqnarray*}
where $r_1$, $r_2$ are as in (\ref{eqn:r1,r2}). This is equal to $\varphi(v_Z^m v_Z^{\ast n})$ given in (\ref{eqn:momentsA(infty)6}).
\hfill
$\Box$

\begin{figure}[tb]
\begin{center}
  \includegraphics[width=60mm]{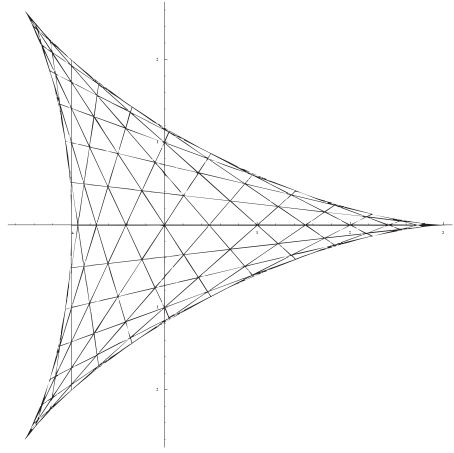}\\
 \caption{Mapping $\mathbb{T}^2$ onto the deltoid $\mathfrak{D}$.} \label{fig:poly-6}
\end{center}
\end{figure}

The quotient $\mathbb{T}^2/\mathbb{Z}_3$, where the $\mathbb{Z}_3$ action is given by left multiplication by $T_3$ is a two-sphere $\mathbb{S}^2$ with three singular points corresponding to  the points $(1,1)$, $(e^{2 \pi i/3}, e^{4 \pi i/3})$, $(e^{4 \pi i/3}, e^{2 \pi i/3})$ in $\mathbb{T}^2$ \cite{farsi/watling:1993}. Under the $\mathbb{Z}^2$ action given by left multiplication by $T_2$ on this two-sphere, we obtain a disc with three singular points, which is topologically equal to the deltoid $\mathfrak{D}$. The boundaries of the deltoid $\mathfrak{D}$ are given by the lines $\theta_1 = 1 - \theta_2$, $\theta_1 = 2 \theta_2$ and $2 \theta_1 = \theta_2$. The diagonal $\theta_1 = \theta_2$ in $\mathbb{T}^2$ is mapped to the real interval $[-1,3] \subset \mathfrak{D}$. The mapping of the `horizontal' lines on $\mathbb{T}^2$ between points $(e^{2 \pi i m/12}, e^{2 \pi i n/12})$ and $(e^{2 \pi i (m+1)/12}, e^{2 \pi i n/12})$, and the `vertical' lines on $\mathbb{T}^2$ between points $(e^{2 \pi i m/12}, e^{2 \pi i n/12})$ and $(e^{2 \pi i m/12}, e^{2 \pi i (n+1)/12})$, onto $\mathfrak{D}$, for $0 \leq m,n \leq 11$, is illustrated in Figure \ref{fig:poly-6}.

\begin{figure}[tb]
\begin{center}
  \includegraphics[width=55mm]{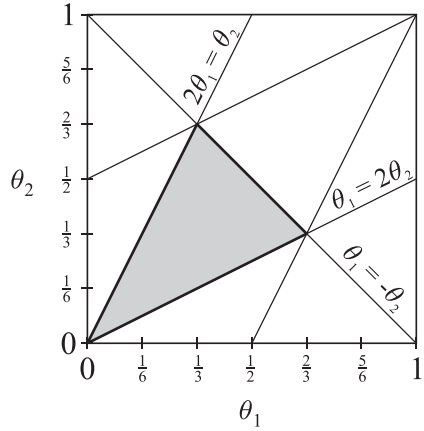}\\
 \caption{A fundamental domain $C$ of $\mathbb{T}^2/S_3$.} \label{fig:poly-4}
\end{center}
\end{figure}

Thus the quotient $\mathbb{T}^2/S_3$ is topologically equal to the deltoid $\mathfrak{D}$. A fundamental domain $C$ of $\mathbb{T}^2$ under the action of the group $S_3$ is illustrated in Figure \ref{fig:poly-4}, where the axes are labelled by the parameters $\theta_1$, $\theta_2$ in $(e^{2 \pi i \theta_1},e^{2 \pi i \theta_2}) \in \mathbb{T}^2$. The boundaries of $C$ map to the boundaries of the deltoid $\mathfrak{D}$. The torus $\mathbb{T}^2$ contains six copies of $C$.

\begin{figure}[tb]
\begin{center}
  \includegraphics[width=80mm]{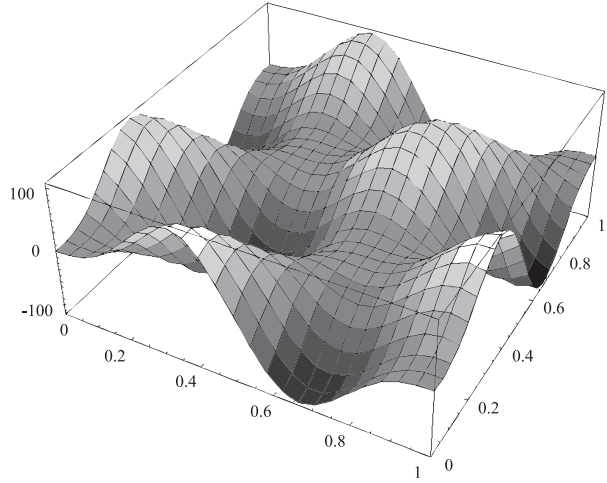}\\
 \caption{The Jacobian $J$.} \label{fig:poly-7}
\end{center}
\end{figure}

We will now determine the spectral measure $\mu_{v_Z}$ over $\mathfrak{D}$. Now
\begin{eqnarray*}
\lefteqn{ \int_{\mathbb{T}^2} (\omega_1 + \omega_2^{-1} + \omega_1^{-1}\omega_2)^m (\omega_1^{-1} + \omega_2 + \omega_1\omega_2^{-1})^n \mathrm{d}\omega_1 \; \mathrm{d}\omega_2 } \\
& = & 6 \int_{C} (\omega_1 + \omega_2^{-1} + \omega_1^{-1}\omega_2)^m (\omega_1^{-1} + \omega_2 + \omega_1\omega_2^{-1})^n \mathrm{d}\omega_1 \; \mathrm{d}\omega_2 \\
& = & 6 \int (e^{2 \pi i \theta_1} + e^{-2 \pi i \theta_2} + e^{2 \pi i (\theta_2 - \theta_1)})^m (e^{-2 \pi i \theta_1} + e^{2 \pi i \theta_2} + e^{2 \pi i (\theta_1 - \theta_2)})^n \mathrm{d}\theta_1 \; \mathrm{d}\theta_2,
\end{eqnarray*}
where the last integral is over the values of $\theta_1$, $\theta_2$ such that $(e^{2 \pi i \theta_1},e^{2 \pi i \theta_2}) \in C$.
Under the change of variable $z = e^{2 \pi i \theta_1} + e^{-2 \pi i \theta_2} + e^{2 \pi i (\theta_2 - \theta_1)}$, we have
\begin{eqnarray*}
x & := & \mathrm{Re}(z) \;\; = \;\; \cos(2 \pi \theta_1) + \cos(2 \pi \theta_2) + \cos(2 \pi (\theta_2 - \theta_1)), \\
y & := & \mathrm{Im}(z) \;\; = \;\; \sin(2 \pi \theta_1) - \sin(2 \pi \theta_2) + \sin(2 \pi (\theta_2 - \theta_1)).
\end{eqnarray*}
Then
\begin{eqnarray}
\lefteqn{ \int_{\mathbb{T}^2} (\omega_1 + \omega_2^{-1} + \omega_1^{-1}\omega_2)^m (\omega_1^{-1} + \omega_2 + \omega_1\omega_2^{-1})^n \mathrm{d}\omega_1 \; \mathrm{d}\omega_2 } \nonumber \\
& = & 6 \int_{\mathfrak{D}} (x+iy)^m (x+iy)^n |J^{-1}| \mathrm{d}x \; \mathrm{d}y, \label{eqn:integral-A(infty)6_over_S}
\end{eqnarray}
where the Jacobian $J = \mathrm{det}(\partial(x,y)/\partial(\theta_1,\theta_2))$ is the determinant of the Jacobian matrix.
We find that the Jacobian $J = J (\theta_1,\theta_2)$ is given by
\begin{equation} \label{eqn:J[theta]}
J (\theta_1,\theta_2) = 4 \pi^2 (\sin(2 \pi (\theta_1 + \theta_2)) - \sin(2 \pi (2\theta_1 - \theta_2)) - \sin(2 \pi (2\theta_2 - \theta_1))).
\end{equation}
The Jacobian is real and vanishes on the boundary of the deltoid $\mathfrak{D}$. For the values of $\theta_1$, $\theta_2$ such that $(e^{2 \pi i \theta_1},e^{2 \pi i \theta_2})$ are in the interior of the fundamental domain $C$ illustrated in Figure \ref{fig:poly-4}, the value of $J$ is always negative. In fact, restricting to any one of the fundamental domains shown in Figure \ref{fig:poly-4}, the sign of $J$ is constant. It is negative over three of the fundamental domains, and positive over the remaining three.
The Jacobian $J (\theta_1,\theta_2)$ is illustrated in Figure \ref{fig:poly-7}.
When evaluating $J$ at a point in $z \in \mathfrak{D}$, we pull back $z$ to $\mathbb{T}^2$. However, there are six possibilities for $(\omega_1,\omega_2) \in \mathbb{T}^2$ such that $\Phi(\omega_1,\omega_2) = z$, one in each of the fundamental domains of $\mathbb{T}^2$ in Figure \ref{fig:poly-4}. Thus over $\mathfrak{D}$, $J$ is only determined up to a sign.
To obtain a positive measure over $\mathfrak{D}$ we take the absolute value $|J|$ of the Jacobian in the integral (\ref{eqn:integral-A(infty)6_over_S}).

Writing $\omega_j = e^{2 \pi i \theta_j}$, $j=1,2$, $J$ is given in terms of $\omega_1,\omega_2 \in \mathbb{T}$ by,
\begin{eqnarray}
J (\omega_1,\omega_2) & = & 4 \pi^2 \mathrm{Im}(\omega_1 \omega_2 - \omega_1^2 \omega_2^{-1} - \omega_1^{-1} \omega_2^2) \nonumber \\
& = & - 2 \pi^2 i (\omega_1 \omega_2 - \omega_1^{-1} \omega_2^{-1} - \omega_1^2 \omega_2^{-1} + \omega_1^{-2} \omega_2 - \omega_1^{-1} \omega_2^2 + \omega_1 \omega_2^{-2}). \qquad \label{eqn:J[omega]}
\end{eqnarray}

Since
\begin{eqnarray*}
\lefteqn{ (\omega_1 \omega_2 - \omega_1^{-1} \omega_2^{-1} - \omega_1^2 \omega_2^{-1} + \omega_1^{-2} \omega_2 - \omega_1^{-1} \omega_2^2 + \omega_1 \omega_2^{-2})^2 } \\
& \qquad = & -6 + 2(\omega_1 \omega_2 + \omega_1^{-1} \omega_2^{-1} + \omega_1 \omega_2^{-2} + \omega_1^2 \omega_2^{-1} + \omega_1^{-1} \omega_2^2 + \omega_1^{-2} \omega_2) \\
& & - 2(\omega_1^3 + \omega_1^{-3} + \omega_2^3 + \omega_2^{-3} + \omega_1^3 \omega_2^{-3} + \omega_1^{-3} \omega_2^3) \\
& & + (\omega_1^2 \omega_2^2 + \omega_1^{-2} \omega_2^{-2} + \omega_1^2 \omega_2^{-4} + \omega_1^4 \omega_2^{-2} + \omega_1^{-2} \omega_2^4 + \omega_1^{-4} \omega_2^2),
\end{eqnarray*}
the square of the Jacobian is invariant under the action of $S_3$. Hence $J^2$ can be written in terms of $z$, $\overline{z}$, and we obtain $J (z,\overline{z})^2 = 4 \pi^4 (27 - 18z\overline{z} + 4z^3 + 4\overline{z}^3 -z^2 \overline{z}^2)$ for $z \in \mathfrak{D}$. Since $J$ is real, $J^2 \geq 0$.
We have the following expressions for the Jacobian $J$:
\begin{eqnarray*}
J (\theta_1,\theta_2) & = & 4 \pi^2 (\sin(2 \pi (\theta_1 + \theta_2)) - \sin(2 \pi (2\theta_1 - \theta_2)) - \sin(2 \pi (2\theta_2 - \theta_1))), \\
J (\omega_1,\omega_2) & = & - 2 \pi^2 i (\omega_1 \omega_2 - \omega_1^{-1} \omega_2^{-1} - \omega_1^2 \omega_2^{-1} + \omega_1^{-2} \omega_2 - \omega_1^{-1} \omega_2^2 + \omega_1 \omega_2^{-2}), \\
|J (z,\overline{z})| & = & 2 \pi^2 \sqrt{27 - 18z\overline{z} + 4z^3 + 4\overline{z}^3 -z^2 \overline{z}^2}, \\
|J (x,y)| & = & 2 \pi^2 \sqrt{27 - 18(x^2+y^2) + 8x(x^2-3y^2) -(x^2+y^2)^2}, \\
|J (r,t)| & = & 2 \pi^2 \sqrt{(1-r)((5+4\cos(6\pi t))^2 r^3 - 9(7+8\cos(6\pi t)) r^2 + 27r +27)},
\end{eqnarray*}
where $0 \leq \theta_1,\theta_2 < 1$, $\omega_1,\omega_2 \in \mathbb{T}$, $z = x+iy \in \mathfrak{D}$ and $0 \leq r \leq 1$, $0 \leq t < 1$. Here the expressions under the square root are always real and non-negative since $J^2$ is.
Consequently:

\begin{theorem}
The spectral measure $\mu_{v_Z}$ (over $\mathfrak{D}$) for the graph $\mathcal{A}^{(6\infty)}$ is
\begin{equation}
\mathrm{d}\mu_{v_Z}(z) = \frac{6}{|J|} \; \mathrm{d}z = \frac{3}{\pi^2 \sqrt{27 - 18z\overline{z} + 4z^3 + 4\overline{z}^3 -z^2 \overline{z}^2}} \; \mathrm{d}z.
\end{equation}
\end{theorem}

We thus have for the fixed point algebra under $\mathbb{T}^2$:
\begin{eqnarray*}
\mathrm{dim}(A(\mathcal{A}^{(6\infty)})_k) & = & \mathrm{dim}\left( \left(\otimes^k M_3 \right)^{\mathbb{T}^2} \right) \;\; = \;\; \sum_{j=0}^k C^{2j}_j (C^k_j)^2 \;\; = \;\; \varphi(|v_Z|^{2k}) \\
& = & \frac{3}{\pi^2} \int_{\mathfrak{D}} |z|^{2k} \frac{1}{\sqrt{27 - 18z\overline{z} + 4z^3 + 4\overline{z}^3 -z^2 \overline{z}^2}} \; \mathrm{d}z.
\end{eqnarray*}

\subsection{Spectral measure for $\mathcal{A}^{(\infty)}$} \label{Sect:spec_measure_for_A(infty)}

\begin{figure}[tb]
\begin{center}
  \includegraphics[width=70mm]{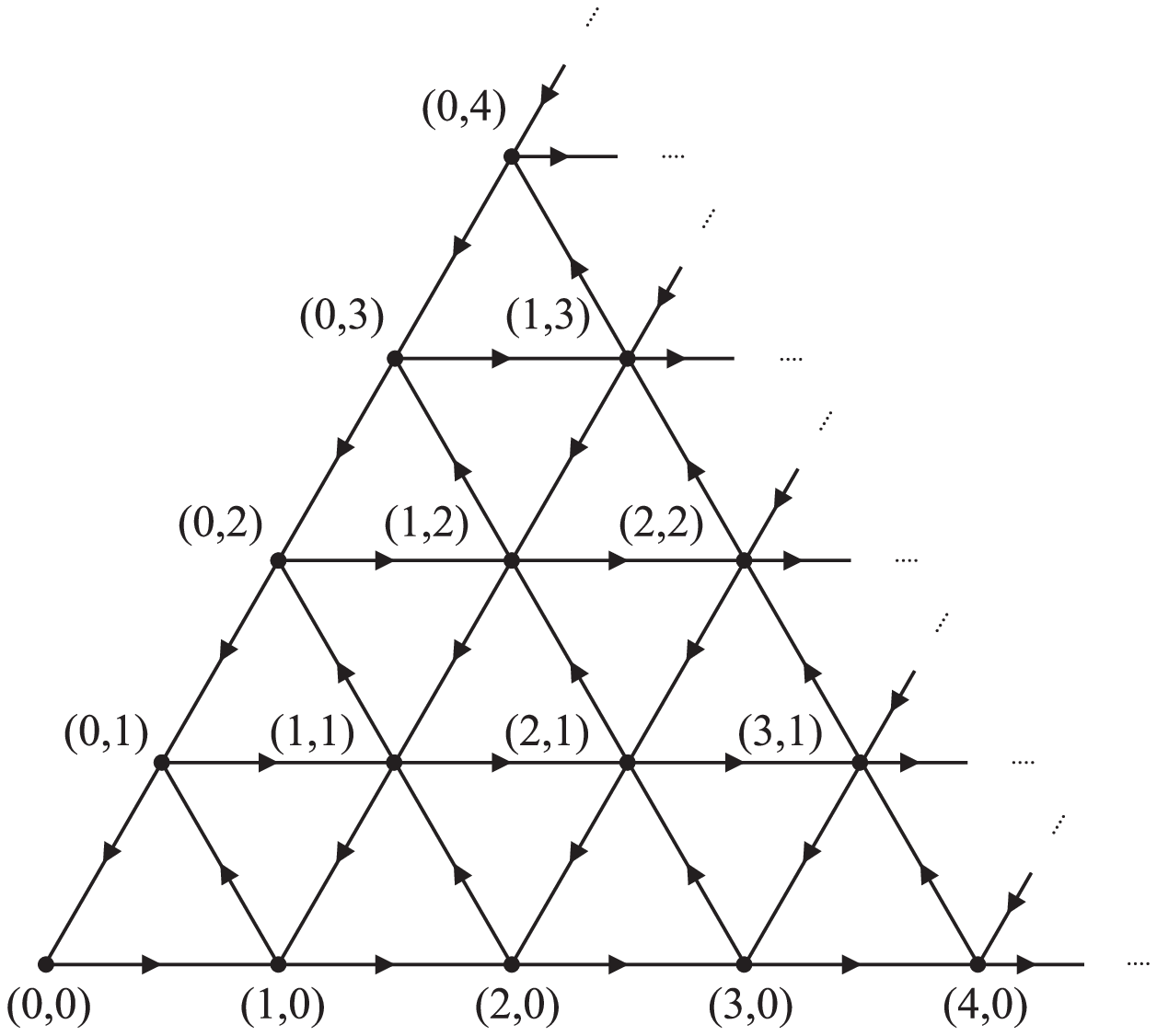}\\
 \caption{The infinite graph $\mathcal{A}^{(\infty)}$} \label{Fig:SU(3)-A(infty)}
\end{center}
\end{figure}

We now consider the fixed point algebra under the action of the group $SU(3)$. The characters of $SU(3)$ satisfy $\chi_{(1,0)} \chi_{(\lambda_1,\lambda_2)} = \chi_{(\lambda_1+1,\lambda_2)} + \chi_{(\lambda_1,\lambda_2-1)} + \chi_{(\lambda_1-1,\lambda_2+1)}$, for any $\lambda_1,\lambda_2 \geq 0$, where $\chi_{(\lambda,-1)} = 0$ for all $\lambda \geq 0$. So the representation graph of $SU(3)$ is identified with the infinite graph $\mathcal{A}^{(\infty)}$ illustrated in Figure \ref{Fig:SU(3)-A(infty)}, with distinguished vertex $\ast = (0,0)$. Hence $(\bigotimes_{\mathbb{Z}} M_3)^{SU(3)} \cong A(\mathcal{A}^{(\infty)})$.

We define a normal operator $v_N$ on $\ell^2(\mathbb{N}) \otimes \ell^2(\mathbb{N})$ by
\begin{equation} \label{def:v_N}
v_N = l \otimes 1 + 1 \otimes l^{\ast} + l^{\ast} \otimes l,
\end{equation}
where $l$ is again the unilateral shift on $\ell^2(\mathbb{N})$. If we regard the element $\Omega \otimes \Omega$ as corresponding to the apex vertex $(0,0)$, and the operators $l \otimes 1$, $l^{\ast} \otimes l$, $1 \otimes l^{\ast}$ as corresponding to the vectors $e_1, e_2, e_3$ on $\mathcal{A}^{(\infty)}$, then $(l^{\lambda_1} \otimes (l^{\ast})^{\lambda_2})(\Omega \otimes \Omega)$ corresponds to the vertex $(\lambda_1, \lambda_2)$ of $\mathcal{A}^{(\infty)}$, for $\lambda_1, \lambda_2 \geq 0$. We see that $v_N$ is identified with the adjacency matrix $\Delta_{\mathcal{A}}$ of $\mathcal{A}^{(\infty)}$, and $v_N^m v_N^{\ast n} (\Omega \otimes \Omega)$ gives a vector $y=(y_{(\lambda_1,\lambda_2)})$ in $\ell^2(\mathcal{A}^{(\infty)})$, where $y_{(\lambda_1,\lambda_2)}$ gives the number of paths of length $m+n$ from $(0,0)$ to the vertex $(\lambda_1,\lambda_2)$, where $m$ edges are on $\mathcal{A}^{(\infty)}$ and $n$ edges are on the reverse graph $\widetilde{\mathcal{A}^{(\infty)}}$. The relation $(l^{\ast} \otimes \, \cdot \, )(\Omega \otimes \, \cdot \,) = 0$ corresponds to the fact that there are no edges in the direction $-e_1$ from a vertex $(0,\lambda_2)$ on the boundary of $\mathcal{A}^{(\infty)}$, $\lambda_2 \geq 0$, and similarly $(\, \cdot \, \otimes l^{\ast})(\, \cdot \, \otimes \Omega) = 0$ corresponds to there being no edges in the direction $e_3$ from a vertex $(\lambda_1,0)$, $\lambda_1 \geq 0$. The relation $(1 \otimes l^{\ast})(l^{\ast} \otimes l)(l \otimes 1) = l^{\ast} l \otimes l^{\ast} l = 1 \otimes 1$ again corresponds to the fact that traveling along edges in directions $e_1$ followed by $e_2$ and then $e_3$ forms a closed loop, and similarly for any permutations of $1 \otimes l^{\ast}$, $l^{\ast} \otimes l$, $l \otimes 1$, but now the product will be 0 along one of the boundaries $\lambda_1 = 0$ or $\lambda_2 = 0$ for certain of the permutations, but 1 everywhere else.

The vector $\Omega \otimes \Omega$ is cyclic in $\ell^2(\mathbb{N}) \otimes \ell^2(\mathbb{N})$. We can show this by induction. Suppose any vector $l^{k_1} \Omega \otimes l^{k_2} \Omega \in \ell^2(\mathbb{N}) \otimes \ell^2(\mathbb{N})$, such that $k_1 + k_2 \leq p$, can be written as a linear combination of elements of the form $v_N^m v_N^{\ast n} (\Omega \otimes \Omega)$ where $m+n \leq p$. This is certainly true when $p=1$ since $v_N (\Omega \otimes \Omega) = (l \otimes 1 + 1 \otimes l^{\ast} + l^{\ast} \otimes l)(\Omega \otimes \Omega) = l \Omega \otimes \Omega$ and $v_N^{\ast} (\Omega \otimes \Omega) = \Omega \otimes l \Omega$. For $j=0,1,\ldots,p$, we have $v_N (l^{p-j} \Omega \otimes l^j \Omega) = l^{p-j+1} \Omega \otimes l^j \Omega + l^{p-j} \Omega \otimes l^{j-1} \Omega + l^{p-j-1} \Omega \otimes l^{j+1} \Omega$. Then $l^{p-j+1} \Omega \otimes l^j \Omega = v_N (l^{p-j} \Omega \otimes l^j \Omega) - l^{p-j} \Omega \otimes l^{j-1} \Omega - l^{p-j-1} \Omega \otimes l^{j+1} \Omega$, and $l^{p-j+1} \Omega \otimes l^j \Omega$, for $j=0,1,\ldots,p$, can be written as a linear combination of elements of the form $v_N^m v_N^{\ast n} (\Omega \otimes \Omega)$ where $m+n \leq p+1$. Since also $\Omega \otimes l^{p+1} \Omega = v_N^{\ast} (\Omega \otimes l^p \Omega) - l \Omega \otimes l^{p-1} \Omega$, then every $l^{k_1} \Omega \otimes l^{k_2} \Omega$, such that $k_1 + k_2 \leq p+1$, can be written as a linear combination of elements of the form $v_N^m v_N^{\ast n} (\Omega \otimes \Omega)$ where $m+n \leq p+1$. Then $\overline{C^{\ast}(v_N) (\Omega \otimes \Omega)} = \ell^2(\mathbb{N}) \otimes \ell^2(\mathbb{N})$.
We define a state $\varphi$ on $C^{\ast}(v_N)$ by $\varphi( \, \cdot \, ) = \langle \, \cdot \, (\Omega \otimes \Omega), \Omega \otimes \Omega \rangle$.
Since $C^{\ast}(v_N)$ is abelian and $\Omega \otimes \Omega$ is cyclic, it is the case that $\varphi$ is faithful.

The moments $\varphi(v_N^m v_N^{\ast n})$ are all zero if $m-n \not \equiv 0 \textrm{ mod } 3$, and for $m \equiv n \textrm{ mod } 3$ the moments $\varphi(v_N^m v_N^{\ast n})$ count the number of paths of length $m+n$ on the $SU(3)$ graph $\mathcal{A}^{(\infty)}$, starting from the apex vertex $(0,0)$, with the first $m$ edges on $\mathcal{A}^{(\infty)}$ and the other $n$ edges on the reverse graph $\widetilde{\mathcal{A}^{(\infty)}}$.
Let $A'(\mathcal{A}^{(\infty)})_{m,n}$ be the algebra generated by pairs $(\eta_1,\eta_2)$ of paths from $(0,0)$ such that $r(\eta_1)=r(\eta_2)$, $|\eta_1|=m$ and $|\eta_2|=n$. Then we define the general path algebra $A'(\mathcal{A}^{(\infty)})$ for the graph $\mathcal{A}^{(\infty)}$ to be $A'(\mathcal{A}^{(\infty)}) = \bigoplus_{m,n} A'(\mathcal{A}^{(\infty)})_{m,n}$. Then $\varphi(v_N^m v_N^{\ast n})$ gives the dimension of the $m,n^{\mathrm{th}}$ level $A'(\mathcal{A}^{(\infty)})_{m,n}$ of the general path algebra $A'(\mathcal{A}^{(\infty)})$.
In particular, $\varphi(v_N^m v_N^{\ast m})$ for $m=n$ gives the dimension of the $m^{\mathrm{th}}$ level of the path algebra for graph $\mathcal{A}^{(\infty)}$, i.e. $\varphi(v_N^m v_N^{\ast m}) = \mathrm{dim}(A(\mathcal{A}^{(\infty)})_m)$.

\begin{figure}[tb]
\begin{center}
  \includegraphics[width=140mm]{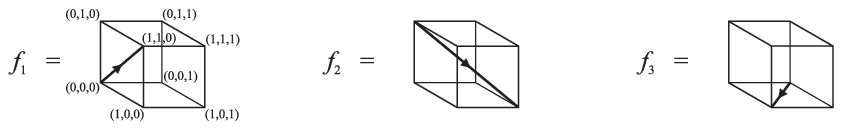}\\
 \caption{The vectors $f_i \in \mathbb{Z}^3$, $i=1,2,3$.} \label{fig:poly-5}
\end{center}
\end{figure}

The moments $\varphi(v_N^m v_N^{\ast n})$ have a realization in terms of a higher dimensional analogue of Catalan paths: Let $E = \{ f_1, f_2, f_3 \}$ be the set of vectors $f_1 = (1,1,0), f_2 = (1,-1,1), f_3 = (1,0,-1) \in \mathbb{Z}^3$, which are illustrated in Figure \ref{fig:poly-5}. These vectors correspond to the vectors $e_i$ above, $i=1,2,3$.

We define the conjugate $\overline{f}$ of a vector $f \in E$ by $\overline{(1,y,z)} = (1,-y,-z)$, and let $\overline{E} = \{ \overline{f}_1, \overline{f}_2, \overline{f}_3 \}$. Let $L$ be the sublattice of $\mathbb{Z}^3$ given by all points with non-negative co-ordinates. Then define $c_{m,n}$ to be the number of paths of length $m+n$ in $L$, starting from $(0,0,0)$ and ending at $(m+n,0,0)$, where $m$ edges are of the form of a vector from $E$ and $n$ edges are of the form of a vector from $\overline{E}$. Then $\varphi(v_N^m v_N^{\ast n}) = c_{m,n}$, and for $m=n$, $\varphi(v_N^m v_N^{\ast m}) = c_{m,m} = \mathrm{dim}(A(\mathcal{A}^{(\infty)})_m)$.

We now consider the probability measure $\mu_{v_N}$ on $\mathfrak{D}$ for the normal element $v_N$.
Since $\varphi$ is a faithful state, by \cite[Remark 2.3.2]{dykema/voiculescu/nica:1992} the support of $\mu_{v_N}$ is equal to the spectrum $\sigma(v_N)$ of $v_N$. Consider the exact sequence $0 \rightarrow \mathcal{K} \rightarrow C^{\ast}(v_N) \rightarrow C^{\ast}(v_N)/\mathcal{K} \rightarrow 0$, where $\mathcal{K} = \mathcal{K}(\ell^2(\mathbb{N}) \otimes \ell^2(\mathbb{N})) \subset B(\ell^2(\mathbb{N}) \otimes \ell^2(\mathbb{N}))$ are the compact operators. Let $\pi: B(\ell^2(\mathbb{N}) \otimes \ell^2(\mathbb{N})) \rightarrow B(\ell^2(\mathbb{N}) \otimes \ell^2(\mathbb{N}))/\mathcal{K}$ be the quotient map, then $\sigma(v_N) \supset \sigma(\pi(v_N))$. Now $\pi(v_N) = u \otimes 1 + 1 \otimes u^{\ast} + u^{\ast} \otimes u$ where $u$ is a unitary which has spectrum $\mathbb{T}$, so that the spectrum of $\pi(v_N)$ is given by $\sigma(\pi(v_N)) = \{ \omega_1 + \omega_2^{-1} + \omega_1^{-1}\omega_2 | \; \omega_1,\omega_2 \in \mathbb{T} \} = \mathfrak{D}$. Then $\sigma(v_N) \subset \mathfrak{D}$.

Consider the measure $\varepsilon(\omega_1,\omega_2)$ on $\mathbb{T}^2$ given by
\begin{eqnarray*}
\lefteqn{ \mathrm{d}\varepsilon(\omega_1,\omega_2) = \frac{1}{24 \pi^4} J(\omega_1,\omega_2)^2 \mathrm{d}\omega_1 \; \mathrm{d}\omega_2 } \\
& = & - \frac{1}{6}(\omega_1 \omega_2 + \omega_1 \omega_2^{-2} + \omega_1^{-2} \omega_2 - \omega_1^{-1} \omega_2^{-1} - \omega_1^2 \omega_2^{-1} - \omega_1^{-1} \omega_2^2)^2 \; \mathrm{d}\omega_1 \, \mathrm{d}\omega_2
\end{eqnarray*}
on $\mathbb{T}^2$, where $\mathrm{d}\omega_j$ is the uniform Lebesgue measure on $\mathbb{T}$, $j=1,2$.
We will prove in the next section that this is the spectral measure (over $\mathbb{T}^2$) of $v_N$, so that $\sigma(v_N) = \mathfrak{D}$.
With this measure we have
\begin{eqnarray*}
\lefteqn{ - \frac{1}{6} \int_{\mathbb{T}^2} (\omega_1 + \omega_2^{-1} + \omega_1^{-1}\omega_2)^m (\omega_1^{-1} + \omega_2 + \omega_1\omega_2^{-1})^n } \\
\lefteqn{ \qquad \qquad \times (\omega_1 \omega_2 + \omega_1 \omega_2^{-2} + \omega_1^{-2} \omega_2 - \omega_1^{-1} \omega_2^{-1} - \omega_1^2 \omega_2^{-1} - \omega_1^{-1} \omega_2^2)^2 \; \mathrm{d}\omega_1 \; \mathrm{d}\omega_2 } \\
& = & - \frac{1}{6} \sum_{\stackrel{0 \leq k_1 + k_2 \leq m}{\scriptscriptstyle{0 \leq l_1 + l_2 \leq n}}} \Bigg( (k_1, k_2, m - k_1 - k_2)! (l_1, l_2, n - l_1 - l_2)! \times \\
& & \int_{\mathbb{T}^2} \omega_1^{r_1} \omega_2^{r_2} (\omega_1 \omega_2 + \omega_1 \omega_2^{-2} + \omega_1^{-2} \omega_2 - \omega_1^{-1} \omega_2^{-1} - \omega_1^2 \omega_2^{-1} - \omega_1^{-1} \omega_2^2)^2 \; \mathrm{d}\omega_1 \, \mathrm{d}\omega_2 \Bigg) \\
& = & - \frac{1}{6} \sum_{a_1,a_2} \sum_{\stackrel{0 \leq k_1 + k_2 \leq m}{\scriptscriptstyle{0 \leq l_1 + l_2 \leq n}}} \Bigg( (k_1, k_2, m - k_1 - k_2)! (l_1, l_2, n - l_1 - l_2)! \; \gamma_{a_1,a_2} \\
& & \qquad \qquad \qquad \qquad \times \int_{\mathbb{T}^2} \omega_1^{r_1 + a_1} \omega_1^{r_2 + a_2} \; \mathrm{d}\omega_1 \; \mathrm{d}\omega_2 \Bigg),
\end{eqnarray*}
where $r_1$, $r_2$ are as in (\ref{eqn:r1,r2}), and the summation is over all integers $a_1$, $a_2$ such that $(a_1,a_2) \in \Upsilon = \{ (\lambda_1,\lambda_2)| \; \lambda_1 \equiv \lambda_2 \textrm{ mod } 3, |\lambda_1 + \lambda_2| \leq 4, |\lambda_1|+|\lambda_2| \leq 6 \}$. The set $\Upsilon$ is the set of all pairs $(a_1,a_2)$ of exponents of $\omega_1^{a_1} \omega_2^{a_2}$ that appear in the expansion of $(\omega_1 \omega_2 + \omega_1 \omega_2^{-2} + \omega_1^{-2} \omega_2 - \omega_1^{-1} \omega_2^{-1} - \omega_1^2 \omega_2^{-1} - \omega_1^{-1} \omega_2^2)^2$, and the integers $\gamma_{a_1,a_2}$ are the corresponding coefficients. Let $b_1 = (2a_1 + a_2)/3$ and $b_2 = (a_1 + 2a_2)/3$. The $m,n^{\mathrm{th}}$ moment for the measure $d\varepsilon(\omega_1,\omega_2)$ is zero if $m \not \equiv 0 \textrm{ mod } 3$, and for $n = m+3r$, $r \in \mathbb{Z}$, the $m,n^{\mathrm{th}}$ moment is given by
\begin{equation} \label{eqn:moments_using_J2uniform_measure}
- \frac{1}{6} \sum_{\stackrel{k_1,k_2}{\scriptscriptstyle{a_1,a_2}}} \gamma_{a_1,a_2} (k_1, k_2, m - k_1 - k_2)! \; (k_1 + r + b_1, k_2 + r - b_2, m + r - b_1 + b_2 - k_1 - k_2)!
\end{equation}
where the summation is over all $a_1, a_2 \in \mathbb{Z}$ such that $(a_1,a_2) \in \Upsilon$, and all non-negative integers $k_1$, $k_2$ such that
\begin{eqnarray}
\mathrm{max}(0,-r-b_1) & \leq & k_1 \leq \mathrm{min}(m,m+2r-b_1) \label{eqn:condition1_k1,k2} \\
\mathrm{max}(0,-r+b_2) & \leq & k_2 \leq \mathrm{min}(m,m+2r+b_2) \label{eqn:condition2_k1,k2} \\
k_1 + k_2 & \leq & \mathrm{min}(m,m+r-b_1+b_2). \label{eqn:condition3_k1,k2}
\end{eqnarray}

As in (\ref{eqn:integral-A(infty)6_over_S}), under the change of variables $\omega_1 + \omega_2^{-1} + \omega_1^{-1} \omega_2 = z$, the spectral measure $\mu_{v_N}(z)$ is given by
$$\mathrm{d}\mu_{v_N}(z) = \frac{6}{|J|} \frac{1}{24 \pi^4} J^2 \; \mathrm{d}z = \frac{1}{4\pi^4} |J| \; \mathrm{d}z.$$

We will have for the fixed point algebra under $SU(3)$:
\begin{eqnarray*}
\mathrm{dim}(A(\mathcal{A}^{(\infty)})_k) & = & \mathrm{dim}\left( \left(\otimes^k M_3 \right)^{SU(3)} \right) \;\; = \;\; \varphi(|v_N|^{2k}) \\
& = & \frac{1}{2\pi^2} \int_{\mathfrak{D}} |z|^{2k} \sqrt{27 - 18z\overline{z} + 4z^3 + 4\overline{z}^3 -z^2 \overline{z}^2} \; \mathrm{d}z.
\end{eqnarray*}

\section{Spectral measures for $\mathcal{ADE}$ graphs via nimreps} \label{sect:spec_measureSU(3)ADE}

Let $\Delta_{\mathcal{G}}$ be the adjacency matrix of a finite graph $\mathcal{G}$ with $s$ vertices, such that $\Delta_{\mathcal{G}}$ is normal. The $m,n^{\mathrm{th}}$ moment $\int z^m \overline{z}^n \mathrm{d}\mu(z)$ is given by $\langle \Delta_{\mathcal{G}}^m (\Delta_{\mathcal{G}}^{\ast})^n e_1, e_1 \rangle$, where $e_1$ is the basis vector in $\ell^2(\mathcal{G})$ corresponding to the distinguished vertex $\ast$ of $\mathcal{G}$.
For convenience we will use the notation
\begin{equation} \label{def:Rmn}
R_{m,n}(\omega_1,\omega_2) := (\omega_1 + \omega_2^{-1} + \omega_1^{-1}\omega_2)^m (\omega_1^{-1} + \omega_2 + \omega_1\omega_2^{-1})^n,
\end{equation}
so that $\int_{\mathbb{T}^2} R_{m,n}(\omega_1,\omega_2) \mathrm{d}\varepsilon(\omega_1,\omega_2) = \int z^m \overline{z}^n \mathrm{d}\mu(z) = \langle \Delta_{\mathcal{G}}^m (\Delta_{\mathcal{G}}^{\ast})^n e_1, e_1 \rangle$.

Let $\beta^j$ be the eigenvalues of $\mathcal{G}$, with corresponding eigenvectors $x^j$, $j=1,\ldots,s$. Then as for $SU(2)$, $\Delta_{\mathcal{G}}^m (\Delta_{\mathcal{G}}^{\ast})^n = \mathcal{U} \Lambda_{\mathcal{G}}^m (\Lambda_{\mathcal{G}}^{\ast})^n \mathcal{U}^{\ast}$, where $\Lambda_{\mathcal{G}}$ is the diagonal matrix $\Lambda_{\mathcal{G}} = \mathrm{diag}(\beta^1, \beta^2, \ldots, \beta^s)$ and $\mathcal{U} = (x^1, x^2, \ldots, x^s)$, so that
\begin{eqnarray}
\int_{\mathbb{T}^2} R_{m,n}(\omega_1,\omega_2) \mathrm{d}\varepsilon(\omega_1,\omega_2) & = & \langle \mathcal{U} \Lambda_{\mathcal{G}}^m (\Lambda_{\mathcal{G}}^{\ast})^n \mathcal{U}^{\ast} e_1, e_1 \rangle \;\; = \;\; \langle \Lambda_{\mathcal{G}}^m (\Lambda_{\mathcal{G}}^{\ast})^n \mathcal{U}^{\ast} e_1, \mathcal{U}^{\ast} e_1 \rangle \nonumber \\
& = & \sum_{j=1}^s (\beta^j)^m (\overline{\beta^j})^n |y_j|^2, \label{eqn:moments_general_SU(3)graph}
\end{eqnarray}
where $y_j = x^j_1$ is the first entry of the eigenvector $x^j$.

For a finite $\mathcal{ADE}$ graph $\mathcal{G}$ with Coxeter exponents $\mathrm{Exp}$, its eigenvalues $\beta^{(\lambda)}$ are ratios of the $S$-matrix given by $\beta^{(\lambda)} = S_{\rho \lambda}/S_{0 \lambda}$, for $\lambda \in \mathrm{Exp}$, with corresponding eigenvectors $(\psi^{\lambda}_a)_{a \in \mathfrak{V}(\mathcal{G})}$. Then (\ref{eqn:moments_general_SU(3)graph}) becomes
\begin{equation} \label{eqn:moments_SU(3)}
\int_{\mathbb{T}^2} R_{m,n}(\omega_1,\omega_2) \mathrm{d}\varepsilon(\omega_1,\omega_2) = \sum_{\lambda \in \mathrm{Exp}} (\beta^{(\lambda)})^m (\overline{\beta^{(\lambda)}})^n |\psi^{\lambda}_{\ast}|^2,
\end{equation}
where $\ast$ is the distinguished vertex of $\mathcal{G}$ with lowest Perron-Frobenius weight.

\subsection{Graphs $\mathcal{A}^{(l)}$, $l \leq \infty$.}

The distinguished vertex $\ast$ of the graph $\mathcal{A}^{(l)}$ is the apex vertex $(0,0)$. Its eigenvalues $\beta^{(\lambda)}$ are given by the ratio $S_{\rho \lambda}/S_{0 \lambda}$, with corresponding eigenvectors $\psi^{\lambda}_{\mu} = S_{\mu,\lambda}$, where the exponents of $\mathcal{A}^{(l)}$ are $\mathrm{Exp} = \{ (\lambda_1,\lambda_2) | \; 0 \leq \lambda_1,\lambda_2 \leq l-3; \; \lambda_1 + \lambda_2 \leq l-3 \}$, and the $S$-matrix for $SU(3)$ at level $k = l-3$ is given by \cite{gannon:1994}:
\begin{eqnarray*}
l \sqrt{3} i \; S_{\mu,\lambda} & = & e^{\xi (2 \lambda_1' \mu_1' + \lambda_1' \mu_2' + \lambda_2' \mu_1' + 2 \lambda_2' \mu_2')} + e^{\xi (\lambda_2' \mu_1' - \lambda_1' \mu_1' + 2 \lambda_1' \mu_2' - \lambda_2' \mu_2')} \\
& & + e^{\xi (\lambda_1' \mu_2' - \lambda_1' \mu_1' - 2 \lambda_2' \mu_1' - \lambda_2' \mu_2')} - e^{\xi (- 2 \lambda_1' \mu_2' - \lambda_1' \mu_1' - \lambda_2' \mu_2' - 2 \lambda_2' \mu_1')} \\
& & - e^{\xi (2 \lambda_1' \mu_1' + \lambda_1' \mu_2' + \lambda_2' \mu_1' - \lambda_2' \mu_2')} - e^{\xi (\lambda_1' \mu_2' - \lambda_1' \mu_1' + \lambda_2' \mu_1' + 2 \lambda_2' \mu_2')}
\end{eqnarray*}
where $\xi = - 2 \pi i /3l$, $\lambda = (\lambda_1,\lambda_2)$, $\mu = (\mu_1,\mu_2)$, and $\lambda_j' = \lambda_j + 1$, $\mu_j' = \mu_j + 1$, for $j = 1,2$.
Then setting $\mu = (0,0)$ we obtain
\begin{eqnarray}
\psi^{\lambda}_{\ast} & = & \frac{2}{l \sqrt{3}} \left( \sin(2\lambda_1' \pi/l) + \sin(2\lambda_2' \pi/l) - \sin(2(\lambda_1' + \lambda_2') \pi/l) \right) \label{eqn:PF-evector2} \\
& = & - \frac{1}{2 \sqrt{3} \pi^2 l} J \left( (\lambda_1 + 2\lambda_2 + 3)/3l, (2\lambda_1 + \lambda_2 + 3)/3l \right), \label{eqn:PF-evector=J}
\end{eqnarray}
where in (\ref{eqn:PF-evector=J}) $\theta_1 = (\lambda_1 + 2 \lambda_2 + 3)/3l$ and $\theta_2 = (2 \lambda_1 + \lambda_2 + 3)/3l$, so that $(\lambda_1 + 1)/l = 2 \theta_2 - \theta_1$ and $(\lambda_2 + 1)/l = 2 \theta_1 - \theta_2$.

Since the $S$-matrix is symmetric, we also have $\psi^{\lambda}_{\mu} = S_{\lambda,\mu}$, so that the Perron-Frobenius eigenvector $\psi^{(0,0)}$ has entries $\psi^{(0,0)}_{\lambda}$ given by (\ref{eqn:PF-evector2}). Since the $S$-matrix is unitary, the eigenvector $\psi^{(0,0)}$ has norm 1. Recall that the Perron-Frobenius eigenvector for $\mathcal{A}^{(l)}$ can also be written in the form \cite{di_francesco:1992}:
\begin{equation} \label{eqn:PF-evector1}
\phi^{(0,0)}_{\lambda} = \frac{ \sin((\lambda_1 +1)\pi/l) \sin((\lambda_2 +1)\pi/l) \sin((\lambda_1+\lambda_2 +2)\pi/l) }{ \sin^2(\pi/l) \sin(2\pi/l) },
\end{equation}
where $\phi^{(0,0)}$ has norm $> 1$. In fact, $\phi^{(0,0)}$ has norm $l \sqrt{3}(8 \sin(2 \pi/l) \sin^2(\pi/l))^{-1}$, so that $\psi^{(0,0)} = 8 \sin(2 \pi/l) \sin^2(\pi/l) \; \phi^{(0,0)} /l \sqrt{3}$. Then by (\ref{eqn:PF-evector=J}),
\begin{eqnarray*}
J(\theta_1,\theta_2) & = & - 2 \sqrt{3} \pi^2 l \; \psi^{(0,0)}_{(l(2\theta_2-\theta_1)-1,l(2\theta_1-\theta_2)-1)} \\
& = & - 2 \sqrt{3} \pi^2 l \; \frac{8}{l \sqrt{3}} \sin(2 \pi/l) \sin^2(\pi/l) \; \phi^{(0,0)}_{(l(2\theta_2-\theta_1)-1,l(2\theta_1-\theta_2)-1)} \\
& = & - 16 \pi^2 \sin((2\theta_2-\theta_1)\pi) \sin((2\theta_1-\theta_2)\pi) \sin((\theta_1+\theta_2)\pi),
\end{eqnarray*}
so that the Jacobian $J(\theta_1,\theta_2)$ can also be written as a product of sine functions.
From this form for $J$ we see that the expression for $J(\omega_1,\omega_2)$ in (\ref{eqn:J[omega]}) factorizes as
$$J (\omega_1,\omega_2) = - 2 \pi^2 i (u_1^{-1} u_2^2 - u_1 u_2^{-2})(u_1^2 u_2^{-1} - u_1^{-2} u_2)(u_1 u_2 - u_1^{-1} u_2^{-1}),$$
where $u_1 = \omega_1^{1/2}$ and $u_2 = \omega_2^{1/2}$ take their values in $\{ e^{i \theta} | \; 0 \leq \theta < \pi \}$.

We now compute the spectral measure for $\mathcal{A}^{(l)}$.
The exponents of $\mathcal{A}^{(l)}$ are all the vertices of $\mathcal{A}^{(l)}$, i.e. $\mathrm{Exp} = \{ (\lambda_1,\lambda_2) | \; \lambda_1,\lambda_2 \geq 0; \; \lambda_1 + \lambda_2 \leq l-3 \}$. Then summing over all $(\lambda_1,\lambda_2) \in \mathrm{Exp}$ corresponds to summing over all $(\theta_1, \theta_2) \in \{ (q_1/3l, q_2/3l) | \; q_1, q_2 = 0,1,\ldots,3l-1 \}$, such that $\theta_1 + \theta_2 \equiv 0 \textrm{ mod } 3$ and
\begin{eqnarray*}
2 \theta_2 - \theta_1 & = & (\lambda_1 + 1)/l \;\; \geq \;\; 1/l, \qquad \qquad
2 \theta_1 - \theta_2 \;\; = \;\; (\lambda_2 + 1)/l \;\; \geq \;\; 1/l, \\
\theta_1 + \theta_2 & = & (\lambda_1 + \lambda_2 + 2)/l \;\; \leq \;\; (l-1)/l \;\; = \;\; 1 - 1/l.
\end{eqnarray*}
Let $L_{(\theta_1,\theta_2)}$ be the set of all such $(\theta_1, \theta_2)$, and let $C_l$ be the set of all $(\omega_1, \omega_2) \in \mathbb{T}$, where $\omega_j = e^{2 \pi i \theta_j}$, $j=1,2$, such that $(\theta_1,\theta_2) \in L_{(\theta_1,\theta_2)}$. It is easy to check that $\beta^{(\lambda)} = \omega_1 + \omega_2^{-1} + \omega_1^{-1} \omega_2$.
Using (\ref{eqn:moments_SU(3)}),
\begin{eqnarray}
\lefteqn{\int_{\mathbb{T}^2} R_{m,n}(\omega_1,\omega_2) \mathrm{d}\varepsilon(\omega_1,\omega_2) } \nonumber \\
& = & \frac{1}{12 \pi^4 l^2} \sum_{\lambda \in \mathrm{Exp}} (\beta^{(\lambda)})^m (\overline{\beta^{(\lambda)}})^n J \left( (2 \lambda_1 + \lambda_2 + 3)/3l, (\lambda_1 + 2 \lambda_2 + 3)/3l \right)^2 \nonumber \\
& = & - \frac{1}{3 l^2} \sum_{(\omega_1,\omega_2) \in C_l} (\omega_1 + \omega_2^{-1} + \omega_1^{-1}\omega_2)^m (\omega_1^{-1} + \omega_2 + \omega_1\omega_2^{-1})^n \nonumber \\
& & \qquad \qquad \times (\omega_1 \omega_2 + \omega_1 \omega_2^{-2} + \omega_1^{-2} \omega_2 - \omega_1^{-1} \omega_2^{-1} - \omega_1^2 \omega_2^{-1} - \omega_1^{-1} \omega_2^2)^2. \label{eqn:SU(3)sum_for_measure-A(n)}
\end{eqnarray}

\begin{figure}[tb]
\begin{center}
  \includegraphics[width=55mm]{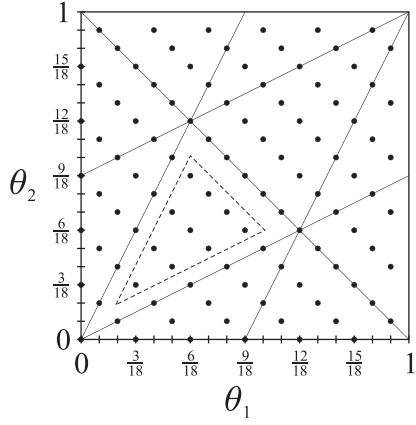}\\
 \caption{The points $(\theta_1,\theta_2)$ such that $(e^{2 \pi i \theta_1}, e^{2 \pi i \theta_2}) \in D_6$.} \label{fig:D6}
\end{center}
\end{figure}

If we let $C$ be the limit of $C_l$ as $l \rightarrow \infty$, then $C$ is a fundamental domain of $\mathbb{T}^2$ under the action of the group $S_3$, illustrated in Figure \ref{fig:poly-4}. Since $J=0$ along the boundary of $C$, which is mapped to the boundary of $\mathfrak{D}$ under the map $\Phi:\mathbb{T}^2 \rightarrow \mathfrak{D}$, we can take the summation in (\ref{eqn:SU(3)sum_for_measure-A(n)}) to include points on the boundary of $C$. Since $J^2$ is invariant under the action of $S_3$, we have
\begin{eqnarray}
\lefteqn{ \int_{\mathbb{T}^2} R_{m,n}(\omega_1,\omega_2) \mathrm{d}\varepsilon(\omega_1,\omega_2) } \nonumber \\
& = & - \frac{1}{6} \frac{1}{3 l^2} \sum_{(\omega_1,\omega_2) \in D_l} (\omega_1 + \omega_2^{-1} + \omega_1^{-1}\omega_2)^m (\omega_1^{-1} + \omega_2 + \omega_1\omega_2^{-1})^n \nonumber \\
& & \qquad \quad \times (\omega_1 \omega_2 + \omega_1 \omega_2^{-2} + \omega_1^{-2} \omega_2 - \omega_1^{-1} \omega_2^{-1} - \omega_1^2 \omega_2^{-1} - \omega_1^{-1} \omega_2^2)^2, \label{eqn:SU(3)sum_for_measure-A(n)2}
\end{eqnarray}
where
\begin{equation} \label{def:Dl}
D_l = \{ (e^{2 \pi i q_1/3l}, e^{2 \pi i q_2/3l}) \in \mathbb{T}^2 | \; q_1,q_2 = 0, 1, \ldots, 3l-1; q_1 + q_2 \equiv 0 \textrm{ mod } 3 \}
\end{equation}
is the image of $C_l$ under the action of $S_3$. We illustrate the points $(\theta_1,\theta_2)$ such that $(e^{2 \pi i \theta_1}, e^{2 \pi i \theta_2}) \in D_6$ in Figure \ref{fig:D6}. Notice that the points in the interior of the fundamental domain $C$, those enclosed by the dashed line, correspond to the vertices of the graph $\mathcal{A}^{(6)}$.

The number $\sharp^{(l)}_{\mathrm{int}}$ of such pairs in the interior of a fundamental domain $C$ can be seen to be equal to $n^{(l)} = (l-2)(l-1)/2$, where $n^{(l)}$ is the number of vertices of $\mathcal{A}^{(l)}$, whilst the number $\sharp^{(l)}_{\partial}$ of such pairs along the boundary of $C$ is $n^{(l+3)} - n^{(l)} = [(l+1)(l+2) - (l-2)(l-1)]/2 = 3l$. Then the total number of such pairs over the whole of $\mathbb{T}^2$ is $|D_l| = 6 \sharp^{(l)}_{\mathrm{int}} + 3 \sharp^{(l)}_{\partial} - 6$ since we count the interior of $C$ six times but only count its boundary three times. The vertices at the corners of the boundary of $C$ are overcounted twice each, hence the term $-6$. So $|D_l| = 3(l-2)(l-1) + 9l - 6 = 3l^2$, and we have
\begin{eqnarray*}
\lefteqn{ \int_{\mathbb{T}^2} R_{m,n}(\omega_1,\omega_2) \mathrm{d}\varepsilon(\omega_1,\omega_2) } \\
& = & - \frac{1}{6} \frac{1}{|D_l|} \sum_{(\omega_1,\omega_2) \in D_l} (\omega_1 + \omega_2^{-1} + \omega_1^{-1}\omega_2)^m (\omega_1^{-1} + \omega_2 + \omega_1\omega_2^{-1})^n \\
& & \qquad \qquad \qquad \times (\omega_1 \omega_2 + \omega_1 \omega_2^{-2} + \omega_1^{-2} \omega_2 - \omega_1^{-1} \omega_2^{-1} - \omega_1^2 \omega_2^{-1} - \omega_1^{-1} \omega_2^2)^2 \\
& = & - \frac{1}{6} \int_{\mathbb{T}^2} (\omega_1 + \omega_2^{-1} + \omega_1^{-1}\omega_2)^m (\omega_1^{-1} + \omega_2 + \omega_1\omega_2^{-1})^n \\
& & \qquad \times (\omega_1 \omega_2 + \omega_1 \omega_2^{-2} + \omega_1^{-2} \omega_2 - \omega_1^{-1} \omega_2^{-1} - \omega_1^2 \omega_2^{-1} - \omega_1^{-1} \omega_2^2)^2 \; \mathrm{d}^{(l)}(\omega_1,\omega_2),
\end{eqnarray*}
where $\mathrm{d}^{(l)}$ is the uniform measure over $D_l$.
Then we have proved the following:

\begin{theorem}
The spectral measure of $\mathcal{A}^{(l)}$ (over $\mathbb{T}^2$) is given by
\begin{equation}
\mathrm{d}\varepsilon(\omega_1,\omega_2) = \frac{1}{24 \pi^4} J(\omega_1,\omega_2)^2 \mathrm{d}^{(l)}(\omega_1,\omega_2).
\end{equation}
\end{theorem}

We can now easily deduce the spectral measure of $\mathcal{A}^{(\infty)}$ claimed in Section \ref{Sect:spec_measure_for_A(infty)}. Letting $l \rightarrow \infty$, the measure $\mathrm{d}^{(l)}(\omega_1,\omega_2)$ becomes the uniform Lebesgue measure $\mathrm{d}\omega_1 \; \mathrm{d}\omega_2$ on $\mathbb{T}^2$:

\begin{theorem} \label{thm:measureA(infty)}
The spectral measure of $\mathcal{A}^{(\infty)}$ (over $\mathbb{T}^2$) is
\begin{equation} \label{eqn:measure_overT2-A(infty)}
\mathrm{d}\varepsilon(\omega_1,\omega_2) = \frac{1}{24 \pi^4} J(\omega_1,\omega_2)^2 \mathrm{d}\omega_1 \; \mathrm{d}\omega_2,
\end{equation}
where $\mathrm{d}\omega$ is the uniform Lebesgue measure over $\mathbb{T}$. Over $\mathfrak{D}$, the spectral measure $\mu_{v_N}(z)$ of $\mathcal{A}^{(\infty)}$ is
\begin{equation}\label{eqn:measure_overS-A(infty)}
\mathrm{d}\mu_{v_N}(z) = \frac{1}{2 \pi^2} \sqrt{27 - 18z\overline{z} + 4z^3 + 4\overline{z}^3 -z^2 \overline{z}^2} \; \mathrm{d}z.
\end{equation}
\end{theorem}

\noindent \emph{Remark:}
For vertices $\nu$ of $\mathcal{A}^{(n)}$ we define polynomials $S_{\nu}(x,y)$ by $S_{(0,0)}(x,y) = 1$, and $x S_{\nu}(x,y) = \sum_{\mu} \Delta_{\mathcal{A}}(\nu, \mu) S_{\mu}(x,y)$ and $y S_{\nu}(x,y) = \sum_{\mu} \Delta_{\mathcal{A}}^T(\nu, \mu) S_{\mu}(x,y)$.
For concrete values of the first few $S_{\mu}(x,y)$ see \cite[p. 610]{evans/kawahigashi:1998}.
Gepner \cite{gepner:1991} proved that this is the measure required to make these polynomials $S_{\mu}(z,\overline{z})$ orthogonal, i.e.
$$\frac{1}{2 \pi^2} \int_{\mathbb{T}^2} S_{\mu}(z,\overline{z}) \overline{S_{\nu}(z,\overline{z})} \sqrt{27 - 18z\overline{z} + 4z^3 + 4\overline{z}^3 -z^2 \overline{z}^2} \; \mathrm{d}z = \delta_{\mu,\nu}.$$

Then in particular, it follows from Theorem \ref{thm:measureA(infty)} that the dimension of the $n^{\mathrm{th}}$ level of the path algebra for $\mathcal{A}^{(\infty)}$ is given by (\ref{eqn:moments_using_J2uniform_measure}) with $m = n$ (i.e. $r=0$), or equivalently by the integral $\int_{\mathfrak{D}} |z|^{2m} \mathrm{d}\mu_{v_N}(z)$ with measure given by (\ref{eqn:measure_overS-A(infty)}).

The dimension of the irreducible representation $\pi_{\lambda}$ of the Hecke algebra $H_n(q)$, labelled by a Young diagram $\lambda = (p_1, p_2, n-p_1-p_2)$ with at most 3 rows, is given by the determinantal formula (see e.g. \cite{sagan:2001}):
\begin{equation} \label{eqn:determinantal_formula}
\mathrm{dim}(\pi_{\lambda}) = n! \left| \begin{array}{ccc}
                                            1/p_1! & 1/(p_1+1)! & 1/(p_1+2)! \\
                                            1/(p_2-1)! & 1/p_2! & 1/(p_2+1)! \\
                                            1/(n-p_1-p_2-2)! &  1/(n-p_1-p_2-1)! &  1/(n-p_1-p_2)! \end{array} \right|,
\end{equation}
where $1/q!$ is understood to be zero if $q$ is negative. Computing the determinant in equation (\ref{eqn:determinantal_formula}), we can rewrite the right hand side as a sum of multinomial coefficients:
\begin{eqnarray}
\lefteqn{ \mathrm{dim}(\pi_{\lambda}) = (p_1,p_2,n-p_1-p_2)! - (p_1,p_2+1,n-p_1-p_2-1)! } \nonumber \\
& & \qquad + (p_1+1,p_2+1,n-p_1-p_2-2)! - (p_1+1,p_2-1,n-p_1-p_2)! \nonumber \\
& & \qquad + (p_1+2,p_2-1,n-p_1-p_2-1)! - (p_1+2,p_2,n-p_1-p_2-2)! \quad \label{def:f(n)p1,p2}
\end{eqnarray}

We can also obtain another formula for the dimension of $A(\mathcal{A}^{(\infty)})_n$. The number $c^{(n)}_{(\lambda_1,\lambda_2)}$ of paths of length $n$ on the graph $\mathcal{A}^{(\infty)}$ from the apex vertex $(0,0)$ to a vertex $(\lambda_1,\lambda_2)$ is given in \cite{di_francesco:1997} as
\begin{equation} \label{def:c(n)}
c^{(n)}_{(\lambda_1,\lambda_2)} = \frac{(\lambda_1+1)(\lambda_2+1)(\lambda_1 + \lambda_2 + 2) \; n!}{((n+2\lambda_1+\lambda_2+6)/3)!((n-\lambda_1+\lambda_2+3)/3)!((n-\lambda_1-2\lambda_2)/3)!}.
\end{equation}

Then we have the following:

\begin{lemma} \label{Lemma:A(infty)_identities}
Let $c^{(n)}_{(\lambda_1,\lambda_2)}$ be the number of paths of length $n$ from $(0,0)$ to the vertex $(\lambda_1,\lambda_2)$ on the graph $\mathcal{A}^{(\infty)}$, as given in (\ref{def:c(n)}), and let $A'(\mathcal{A}^{(\infty)})$ be the general path algebra defined in Section \ref{Sect:spec_measure_for_A(infty)}.
Then, for fixed integers $m,n < \infty$, the following are all equal:
\begin{itemize}
\item[(1)] $\mathrm{dim}(A'(\mathcal{A}^{(\infty)})_{m,n})$,
\item[(2)] $\frac{1}{2 \pi^2} \int_{\mathfrak{D}} z^m \overline{z}^n \sqrt{27 - 18z\overline{z} + 4z^3 + 4\overline{z}^3 -z^2 \overline{z}^2} \; \mathrm{d}z$,
\item[(3)] $\frac{1}{24\pi^4} \int_{\mathbb{T}^2} (\omega_1 + \omega_2^{-1} + \omega_1^{-1}\omega_2)^m (\omega_1^{-1} + \omega_2 + \omega_1\omega_2^{-1})^n J(\omega_1,\omega_2)^2 \mathrm{d}\omega_1 \; \mathrm{d}\omega_2$,
\item[(4)] $- \frac{1}{6} \sum \gamma_{a_1,a_2} (k_1, k_2, n - k_1 - k_2)! \; (k_1 + r + b_1, k_2 + r - b_2, m + r - b_1 + b_2 - k_1 - k_2)!$,
\item[(5)] $\sum c^{(m)}_{(\lambda_1,\lambda_2)} c^{(n)}_{(\lambda_1,\lambda_2)}$,
\end{itemize}
where in (4), $n = m+3r$, $r \in \mathbb{Z}$, $b_1 = (2a_1 + a_2)/3$, $b_2 = (a_1 + 2a_2)/3$ and the summation is over all $a_1, a_2 \in \mathbb{Z}$ such that $(a_1,a_2) \in \Upsilon$, and all non-negative integers $k_1$, $k_2$ which satisfy (\ref{eqn:condition1_k1,k2})-(\ref{eqn:condition3_k1,k2}). The summation in (5) is over all $0 \leq \lambda_1,\lambda_2 \leq \mathrm{min}(m,n)$ such that $\lambda_1 + \lambda_2 \leq \mathrm{min}(m,n)$ and $m \equiv n \equiv \lambda_1 + 2 \lambda_2 \textrm{ mod } 3$.
\end{lemma}
\emph{Proof:}
The identities (1) $=$ (2) $=$ (3) $=$ (4) were shown above. The identity (1) = (5) is trivial since the dimension of $A'(\mathcal{A}^{(\infty)})_{m,n}$ is equal to the number of pairs of paths (with lengths $m$, $n$ respectively) which begin at $(0,0)$ and end at the same vertex of $\mathcal{A}^{(\infty)}$.
\hfill
$\Box$

\begin{corollary}
Let $f^{(n)}_{p_1,p_2}$ be the sum of multinomial coefficients given by (\ref{def:f(n)p1,p2}). Then, in particular, for fixed $n < \infty$, the following are all equal:
\begin{itemize}
\item[(1)] $\mathrm{dim}((\bigotimes^n M_3)^{SU(3)})$,
\item[(2)] $\frac{1}{2 \pi^2} \int_{\mathfrak{D}} |z|^{2n} \sqrt{27 - 18z\overline{z} + 4z^3 + 4\overline{z}^3 -z^2 \overline{z}^2} \; \mathrm{d}z$,
\item[(3)] $\frac{1}{24\pi^4} \int_{\mathbb{T}^2} |\omega_1 + \omega_2^{-1} + \omega_1^{-1}\omega_2|^{2n} \; J(\omega_1,\omega_2)^2 \, \mathrm{d}\omega_1 \; \mathrm{d}\omega_2$,
\item[(4)] $- \frac{1}{6} \sum \gamma_{a_1,a_2} (k_1, k_2, n - k_1 - k_2)! \; (k_1 + b_1, k_2 - b_2, n - b_1 + b_2 - k_1 - k_2)!$,
\item[(5)] $\sum f^{(n)}_{p_1,p_2}$,
\item[(6)] $\sum (c^{(n)}_{(\lambda_1,\lambda_2)})^2$,
\end{itemize}
where in (4), $b_1 = (2a_1 + a_2)/3$, $b_2 = (a_1 + 2a_2)/3$ and the summation is over all $a_1, a_2 \in \mathbb{Z}$ such that $(a_1,a_2) \in \Upsilon$, and all non-negative integers $k_1$, $k_2$ which satisfy (\ref{eqn:condition1_k1,k2})-(\ref{eqn:condition3_k1,k2}). The summation in (5) is over all $0 \leq p_2 \leq p_1 \leq n$ such that $n-p_1 \leq 2p_2$, whilst the summation in (6) is over all $0 \leq \lambda_1,\lambda_2 \leq n$ such that $\lambda_1 + \lambda_2 \leq n$ and $n \equiv \lambda_1 + 2 \lambda_2 \textrm{ mod } 3$.
\end{corollary}
\emph{Proof:}
The identities (1) $=$ (2) $=$ (3) $=$ (4) $=$ (6) follow from Lemma \ref{Lemma:A(infty)_identities}. The identity (1) = (5) follows from (\ref{def:f(n)p1,p2}) and the fact that $(\bigotimes^n M_3)^{SU(3)} = A(\mathcal{A}^{(\infty)})_n = \bigoplus_{\lambda} \pi_{\lambda}(H_n(q))$, where the summation is again over all Young diagrams $\lambda$ with $n$ boxes.

\subsection{Graphs $\mathcal{D}^{(n)}$, $n \equiv 0 \textrm{ mod } 3$.}

The exponents of $\mathcal{D}^{(3k)}$, for integers $k \geq 2$, are the 0-coloured vertices of $\mathcal{A}^{(3k)}$, i.e. $\mathrm{Exp} = \{ (\lambda_1,\lambda_2)| \; \lambda_1,\lambda_2 \geq 0; \; \lambda_1 + \lambda_2 \leq 3k-3; \; \lambda_1 - \lambda_2 \equiv 0 \textrm{ mod } 3 \}$, where the exponent $(k-1,k-1)$ has multiplicity three.

For $\mathcal{D}^{(3k)}$ we have $|\psi_{\ast}^{\lambda}| = \sqrt{3} S_{(0,0),\lambda}$ for all $\lambda \in \mathrm{Exp}$ except for $\lambda = (k-1,k-1)$. For this exponent however the eigenvalue $\beta^{(k,k)} = 0$, so that this term does not contribute in (\ref{eqn:moments_SU(3)}). Then for $\lambda \neq (k-1,k-1)$, the weight $|\psi_{\ast}^{\lambda}|$ is given by $|\psi_{\ast}^{\lambda}| = J \left( (\lambda_1 + 2\lambda_2 + 3)/3l, (2\lambda_1 + \lambda_2 + 3)/3l \right)/ 6k \pi^2$.

Since the exponents for $\mathcal{D}^{(3k)}$ are all of colour zero, under the above identification between $\lambda_1$, $\lambda_2$ and $\theta_1$, $\theta_2$, the exponents $\lambda$ correspond to all pairs $(\theta_1,\theta_2)$ such that $\theta_1 - \theta_2 \equiv 0 \textrm{ mod } 3$ and $(e^{2 \pi i \theta_1},e^{2 \pi i \theta_2}) \in C$. These pairs $(\theta_1,\theta_2)$ are thus in fact all of the form $(p_1/3k,p_2/3k)$, for $p_1,p_2 \in \{ 1,2,\ldots, 3k-1 \}$. Under the action of $S_3$, these pairs are mapped to the all the points $(q_1,q_2) \in [0,1]^2$ such that $e^{2 \pi i q_j}$ is a $3k^{\mathrm{th}}$ root of unity, for $j=1,2$, except for the points $(q_1,q_2)$ which parameterize the boundary of $\mathfrak{D}$.
However, we can again use the fact that the Jacobian is zero at the points which parameterize the boundary of $\mathfrak{D}$.

Then by (\ref{eqn:moments_SU(3)}) we have
\begin{eqnarray*}
\lefteqn{ \int_{\mathbb{T}^2} R_{m,n}(\omega_1,\omega_2) \mathrm{d}\varepsilon(\omega_1,\omega_2) } \\
& = & \frac{1}{4 \pi^4} \frac{1}{(3k)^2} \sum_{\lambda \in \mathrm{Exp}} (\beta^{(\lambda)})^m (\overline{\beta^{(\lambda)}})^n J \left( (\lambda_1 + 2\lambda_2 + 3)/3l, (2\lambda_1 + \lambda_2 + 3)/3l \right)^2 \\
& = & \frac{1}{24 \pi^4} \frac{1}{(3k)^2} \sum_{\theta_1,\theta_2} (\beta^{(\lambda)})^m (\overline{\beta^{(\lambda)}})^n J(\theta_1,\theta_2)^2.
\end{eqnarray*}
The last summation is over $(\theta_1,\theta_2) \in \{ (p_1/3k,p_2/3k) | \; p_1,p_2 = 1,\ldots, 3k-1 \}$.
Then we have obtained the following result:

\begin{theorem} \label{Thm:spec_measure-D(3k)}
The spectral measure of $\mathcal{D}^{(3k)}$, $k \geq 2$, (over $\mathbb{T}^2$) is
\begin{equation}
\mathrm{d}\varepsilon(\omega_1,\omega_2) = \frac{1}{24 \pi^4} J(\omega_1,\omega_2)^2 \; \mathrm{d}_{3k/2}\omega_1 \; \mathrm{d}_{3k/2}\omega_2,
\end{equation}
where $\mathrm{d}_{3k/2}$ is the uniform measure over the $3k^{\mathrm{th}}$ roots of unity.
\end{theorem}

For the limit as $k \rightarrow \infty$ we simply recover the measure (\ref{eqn:measure_overT2-A(infty)}) for $\mathcal{A}^{(\infty)}$.
This is due to the fact that taking the limit of the graph $\mathcal{D}^{(3k)}$ as $k \rightarrow \infty$ with the vertex $\ast = (0,0)$ as the distinguished vertex, we just obtain the infinite graph $\mathcal{A}^{(\infty)}$. In order to obtain the infinite graph $\mathcal{D}^{(\infty)}$ we must set the distinguished vertex $\ast$ of $\mathcal{D}^{(3k)}$ to be one of the triplicated vertices $(k-1,k-1)_i$, $i=1,2,3$, which come from the fixed vertex $(k-1,k-1)$ of $\mathcal{A}^{(3k)}$ under the $\mathbb{Z}_3$ action. Then using (\ref{eqn:moments_SU(3)}), and taking the limit as $k \rightarrow \infty$, we would obtain the spectral measure for $\mathcal{D}^{(\infty)}$.

\subsection{Graphs $\mathcal{A}^{(l)\ast}$, $l \leq \infty$.}

\begin{figure}[tb]
\begin{center}
  \includegraphics[width=120mm]{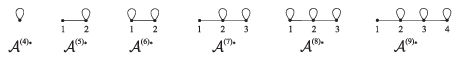}\\
 \caption{$\mathcal{A}^{(n)\ast}$ for $n = 4,5,6,7,8,9$}\label{fig:A(star)-graphs}
\end{center}
\end{figure}

The exponents of $\mathcal{A}^{(l)\ast}$ are $\mathrm{Exp} = \{ (j,j)| \; j=0,1,\ldots, \lfloor (l-3)/2 \rfloor \}$. From \cite{gaberdiel/gannon:2002} its eigenvectors are $\psi^{\lambda}_a = 2 \sqrt{l^{-1}} \sin(2 \pi a (\lambda_1 + 1)/l)$, where $\lambda = (\lambda_1, \lambda_2) \in \mathrm{Exp}$ and $a = 1,2,\ldots,\lfloor (l-1)/2 \rfloor$, as in Figure \ref{fig:A(star)-graphs}. Then
$$\int_{\mathbb{T}^2} R_{m,n}(\omega_1,\omega_2) \mathrm{d}\varepsilon(\omega_1,\omega_2) = \frac{4}{l} \sum_{j=0}^{\lfloor (l-3)/2 \rfloor} (\beta^{(j,j)})^m (\overline{\beta^{(j,j)}})^n \sin^2( 2 \pi (j+1)/l ).$$

Since all the eigenvalues $\beta^{(j,j)}$ of $\mathcal{A}^{(l)\ast}$ are real, there is a map $\Phi_1: \mathbb{T} \rightarrow \mathfrak{D}$ given by $\Phi_1 (u) = u + u^{-1} + 1$ so that the eigenvalues are given by $\Phi_1 (e^{2 \pi i (j+1)/l}) \in [-1,3]$ for $j=0,1,\ldots,\lfloor (l-3)/2 \rfloor \}$. Then the spectral measure of $\mathcal{A}^{(l)\ast}$ can be written as a measure over $\mathbb{T}$. Then with $\widetilde{u} = e^{2 \pi i/l}$, we have
$$\int_{\mathbb{T}} (u + u^{-1})^{m+n} \mathrm{d}\varepsilon(u) = \frac{4}{l} \sum_{j=1}^{\lfloor (l-1)/2 \rfloor} (\widetilde{u}^j + \widetilde{u}^{-j} + 1)^{m+n} \; \sin(\widetilde{u}^j)^2.$$
For all $l$, $\sin(\widetilde{u}^0) = 0$, and $\sin(\widetilde{u}^j) = \sin(\widetilde{u}^l-j)$, for $l=1,2,\ldots,\lfloor (l-1)/2 \rfloor$. If $l$ is even, we also must consider when $j=l/2$. In this case $\sin(\widetilde{u}^{l/2}) = 0$. Then we can write
\begin{eqnarray}
\int_{\mathbb{T}} (u + u^{-1})^{m+n}  \mathrm{d}\varepsilon(u) & = & \frac{2}{l} \sum_{j=0}^l (\widetilde{u}^j + \widetilde{u}^{-j} + 1)^{m+n} \; \sin^2(\widetilde{u}^j) \label{eqn:moment-A(n)star} \\
& = & 2 \int_{\mathbb{T}} (u + u^{-1} + 1)^{m+n} \; \sin^2(u)\mathrm{d}_{l/2}u, \nonumber
\end{eqnarray}
where $\mathrm{d}_p$ is the uniform measure over the $2p^{\mathrm{th}}$ roots of unity. Then we have:

\begin{theorem} \label{Thm:spec_measure-A(n)star}
The spectral measure of $\mathcal{A}^{(l)\ast}$, $l < \infty$, (over $\mathbb{T}$) is
\begin{equation}
\mathrm{d}\varepsilon(u) = \alpha(u) \mathrm{d}_{l/2}u,
\end{equation}
where $\mathrm{d}_{l/2}u$ is the uniform measure over $l^{\mathrm{th}}$ roots of unity, and $\alpha(u) = 2\mathrm{Im}(u)^2$.
\end{theorem}

Since $(u + u^{-1} + 1)^l = \sum_{i=0}^l C^l_i (u + u^{-1})^i$, for even $l = 2k$ we can express the $m,n^{\mathrm{th}}$ moment as a linear combination of the moments of the Dynkin diagram $A_{k-1}$:
$$\int_{\mathbb{T}} (u + u^{-1})^{m+n}  \mathrm{d}\varepsilon(u) = \sum_{j=0}^{m+n} C^{m+n}_j \int_{\mathbb{T}} (u + u^{-1})^j \; 2\mathrm{Im}(u)^2 \mathrm{d}_{l/2}u = \sum_{j=0}^{m+n} C^{m+n}_j \varsigma^j,$$
where $\varsigma^j$ is the $j^{\mathrm{th}}$ moment of $A_{k-1}$.
When $l \rightarrow \infty$, the $j^{\mathrm{th}}$ moment $\varsigma^j$ of $A_{\infty}$ is given by the Catalan number $c_{j/2}$ when $j$ is even, and 0 when $j$ is odd. Then for the infinite graph $\mathcal{A}^{(\infty)\ast}$,
$$\int_{\mathbb{T}} (u + u^{-1})^{m+n}  \mathrm{d}\varepsilon(u) = \sum_{k=0}^{\lfloor (m+n)/2 \rfloor} C^{m+n}_{2k} c_k.$$

In fact, the spectral measure for $\mathcal{A}^{(\infty)\ast}$ has semicircle distribution: Letting $l \rightarrow \infty$ in (\ref{eqn:moment-A(n)star}), we have the approximation of an integral
$$\lim_{l \rightarrow \infty} \frac{2}{l} \sum_{j=0}^l (\widetilde{u}^j + \widetilde{u}^{-j} + 1)^{m+n} \; \sin^2(\widetilde{u}^j) = 2 \int_0^1 (e^{2 \pi i \theta} + e^{-2 \pi i \theta} + 1)^m \; \sin^2 (2 \pi \theta) \mathrm{d}\theta.$$
Making the change of variable $x = e^{2 \pi i \theta} + e^{-2 \pi i \theta} + 1 = 2 \cos(2 \pi \theta) +1$, we have $2 \sin(2 \pi \theta) = \sqrt{4 - (x-1)^2}$, and $\mathrm{d}x/\mathrm{d}\theta = - 4 \pi \sin(2 \pi \theta) = -2 \pi \sqrt{4 - (x-1)^2}$. Then
\begin{eqnarray*}
\int x^m \mathrm{d}\mu(x) & = & 2 \int_0^1 (e^{2 \pi i \theta} + e^{-2 \pi i \theta} + 1)^m \; \sin^2 (2 \pi \theta) \mathrm{d}\theta \\
& = & 4 \int_0^{\frac{1}{2}} (e^{2 \pi i \theta} + e^{-2 \pi i \theta} + 1)^m \; \sin^2 (2 \pi \theta) \mathrm{d}\theta \\
& = & \frac{-4}{8 \pi} \int_{3}^{-1} x^m \sqrt{4 - (x-1)^2} \mathrm{d}x \;\; = \;\; \frac{1}{2 \pi} \int_{-1}^3 x^m \sqrt{4 - (x-1)^2} \mathrm{d}x,
\end{eqnarray*}
which is the semicircle law centered at 1 with radius 2. Then the spectral measure $\mu$ (over $[-1,3]$) for the infinite graph $\mathcal{A}^{(\infty)\ast}$ has semicircle distribution with mean 1 and variance 1, i.e. $\mathrm{d}\mu(x) = \sqrt{4 - (x-1)^2} \mathrm{d}x$.
The graph $\mathcal{A}^{(2l)\ast}$ has adjacency matrix $\Delta^{(2l)\ast} = \Delta_{l-1} + \mathbf{1}$, where $\Delta_{l}$ is the adjacency matrix of the Dynkin diagram $A_{l}$. Hence the spectral measure for $\mathcal{A}^{(2l)\ast}$ is the spectral measure for $A_{l-1}$ but with a shift by one.

\subsection{Graph $\mathcal{E}^{(8)}$}

The spectral measures for the graphs $\mathcal{A}^{(l)}$, $\mathcal{D}^{(3k)}$ are measures of type $\mathrm{d}_{p/2} \times \mathrm{d}_{p/2}$, $J^2 \mathrm{d}_{p/2} \times \mathrm{d}_{p/2}$, $\mathrm{d}^{(p)}$ or $J^2 \mathrm{d}^{(p)}$, for $p \in \mathbb{N}$.
However, we will now show that the spectral measure for $\mathcal{E}^{(8)}$ is not a linear combination of measures of these types.
The exponents of $\mathcal{E}^{(8)}$ are
$$\mathrm{Exp} = \{ (0,0), (5,0), (0,5), (2,2), (2,1), (1,2), (3,0), (2,3), (0,2), (0,3), (3,2), (2,0) \}.$$
Let $\omega = e^{2 \pi i/3}$ and $A$ be the automorphism of order 3 on the vertices of $\mathcal{A}^{(8)}$ given by $A(\mu_1,\mu_2) = (5-\mu_1-\mu_2,\mu_1)$. For the eigenvalues $\beta^{(\lambda)}$, $\beta^{(A(\lambda))} = \omega \beta^{(\lambda)}$ and $\beta^{(A^2(\lambda))} = \overline{\omega} \beta^{(\lambda)}$, the corresponding eigenvectors are $(v^{\lambda},v^{\lambda},v^{\lambda})$, $(v^{\lambda}, \omega v^{\lambda}, \overline{\omega} v^{\lambda})$ and $(v^{\lambda}, \overline{\omega} v^{\lambda}, \omega v^{\lambda})$ respectively, where the row vectors $v^{\lambda}$ are given in \cite[Table 17.3]{di_francesco/mathieu/senechal:1997} (We normalize the eigenvectors so that $||\psi^{\lambda}|| = 1$). Hence $\psi^{\lambda}_{\ast} = \psi^{A(\lambda)}_{\ast} = \psi^{A^2(\lambda)}_{\ast}$ for $\lambda \in \mathrm{Exp}$.
With $\theta_1 = (\lambda_1 + 2\lambda_2 + 3)/24$, $\theta_2 = (2\lambda_1 + \lambda_2 + 3)/24$, we have
\begin{center}
\renewcommand{\arraystretch}{1.5}
\begin{tabular}{|c|c|c|c|} \hline
$\lambda \in \mathrm{Exp}$ & $(\theta_1,\theta_2) \in [0,1]^2$ & $|\psi^{\lambda}_{\ast}|^2$ & $\frac{1}{16\pi^4} J(\theta_1,\theta_2)^2$ \\
\hline $(0,0)$, $(5,0)$, $(0,5)$ & $\left(\frac{1}{8},\frac{1}{8}\right)$, $\left(\frac{1}{3},\frac{13}{24}\right)$, $\left(\frac{13}{24},\frac{1}{3}\right)$ & $\frac{2-\sqrt{2}}{24}$ & $3-2\sqrt{2}$ \\
\hline $(2,2)$, $(2,1)$, $(1,2)$ & $\left(\frac{3}{8},\frac{3}{8}\right)$, $\left(\frac{7}{24},\frac{1}{3}\right)$, $\left(\frac{1}{3},\frac{7}{24}\right)$ & $\frac{2+\sqrt{2}}{24}$ & $3+2\sqrt{2}$ \\
\hline $(3,0)$, $(2,3)$, $(0,2)$ & $\left(\frac{1}{4},\frac{3}{8}\right)$, $\left(\frac{11}{24},\frac{5}{12}\right)$, $\left(\frac{7}{24},\frac{5}{24}\right)$ & $\frac{1}{12}$ & 2 \\
\hline $(0,3)$, $(3,2)$, $(2,0)$ & $\left(\frac{3}{8},\frac{1}{4}\right)$, $\left(\frac{5}{12},\frac{11}{24}\right)$, $\left(\frac{5}{24},\frac{7}{24}\right)$ & $\frac{1}{12}$ & 2 \\
\hline
\end{tabular}
\end{center}

From (\ref{eqn:moments_SU(3)}),
\begin{equation} \label{eqn:moments-E(8)}
\int_{\mathbb{T}^2} R_{m,n}(\omega_1,\omega_2) \mathrm{d}\varepsilon(\omega_1,\omega_2) = \frac{1}{6} \sum_{g \in S_3} \sum_{\lambda \in \mathrm{Exp}} (\beta^{(g(\lambda))})^m (\overline{\beta^{(g(\lambda))}})^n |\psi^{g(\lambda)}_{\ast}|^2.
\end{equation}

\begin{figure}[bt]
\begin{center}
  \includegraphics[width=55mm]{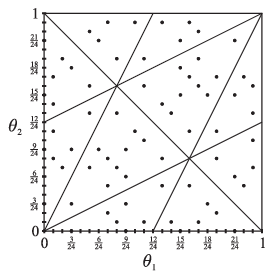}\\
 \caption{The points $(\theta_1,\theta_2) \in \{ g(\lambda) | \; \lambda \in \mathrm{Exp}, g \in S_3 \}$ for $\mathcal{E}^{(8)}$.} \label{fig:poly-11}
\end{center}
\end{figure}

Now the pairs $(\theta_1,\theta_2)$ given by $g(\lambda)$ for $\lambda \in \mathrm{Exp}$, $g \in S_3$, are illustrated in Figure \ref{fig:poly-11}.
Consider the pairs $(\theta_1,\theta_2) = (7/24,8/24), (8/24,13/24), (10/24,11/24)$. For each of these, $(\omega_1,\omega_2) = (e^{2 \pi i \theta_1}, e^{2 \pi i \theta_2}) \in \mathbb{T}^2$ can only be obtained in the integral in (\ref{eqn:moments-E(8)}) from either the product measure $\mathrm{d}_{12} \times \mathrm{d}_{12}$ on pairs of $24^{\mathrm{th}}$ roots of unity, or the uniform measure $\mathrm{d}^{(8)}$ on the elements of $D_8$ ($(7/24,8/24)$, $(8/24,13/24)$, $(10/24,11/24)$ are each in $D_8$, but none are in $D_{k}$ for any integer $k < 8$). Since these points $(\omega_1,\omega_2)$ cannot be obtained independently of each other, we must find a linear combination $\varepsilon' = c_1 \varepsilon_1 + c_2 J^2 \varepsilon_2$ of measures, where $\varepsilon_j$ must be either $\mathrm{d}_{12} \times \mathrm{d}_{12}$ or $\mathrm{d}^{(8)}$ for $j=1,2$ (it doesn't matter at this stage which of the two measures we take $\varepsilon_j$ to be), such that the weight $\varepsilon'(e^{2 \pi i \theta_1}, e^{2 \pi i \theta_2})$ is $(2+\sqrt{2})/24$ for $(\theta_1,\theta_2) = (7/24,8/24)$, $(2-\sqrt{2})/24$ for $(\theta_1,\theta_2) = (8/24,13/24)$ and $1/12$ for $(\theta_1,\theta_2) = (10/24,11/24)$. Suppose for now that $\varepsilon_1 = \varepsilon_2$. Then we must find solutions $c_1,c_2 \in \mathbb{C}$ such that
\begin{equation} \label{eqns:forE(8)}
c_1 + (3-2\sqrt{2})c_2 = \frac{2-\sqrt{2}}{24}, \qquad c_1 + (3+2\sqrt{2})c_2 = \frac{2+\sqrt{2}}{24}, \qquad c_1 + 2c_2 = \frac{1}{12}.
\end{equation}
Solving the first two equations we obtain $c_1 = c_2 = 1/48$. However, substituting for these values into the third equation we get $1/48 + 2/48 = 1/16 \neq 1/12$, hence no solution exists to the equations (\ref{eqns:forE(8)}), and hence the spectral measure for $\mathcal{E}^{(8)}$ is not a linear combination of measures of type $\mathrm{d}_{p/2} \times \mathrm{d}_{p/2}$, $J^2 \mathrm{d}_{p/2} \times \mathrm{d}_{p/2}$, $\mathrm{d}^{(p)}$ or $J^2 \mathrm{d}^{(p)}$, for $p \in \mathbb{N}$.

\subsection{Graph $\mathcal{E}_1^{(12)}$}

We will now show that the spectral measure for $\mathcal{E}_1^{(12)}$ is also not a linear combination of measures of type $\mathrm{d}_{p/2} \times \mathrm{d}_{p/2}$, $J^2 \mathrm{d}_{p/2} \times \mathrm{d}_{p/2}$, $\mathrm{d}^{(p)}$ or $J^2 \mathrm{d}^{(p)}$, for $p \in \mathbb{N}$.
The exponents of $\mathcal{E}_1^{(12)}$ are
$$\mathrm{Exp} = \{ (0,0), (9,0), (0,9), (4,4), (4,1), (1,4), \textrm{ and twice } (2,2), (5,2), (2,5) \}.$$
Computing the first entries of the eigenvectors, we have
\begin{eqnarray*}
|\psi^{(0,0)}_{\ast}|^2 \;\; = \;\; |\psi^{(9,0)}_{\ast}|^2 & = & |\psi^{(0,9)}_{\ast}|^2 \;\; = \;\; (2-\sqrt{3})/36, \\
|\psi^{(4,4)}_{\ast}|^2 \;\; = \;\; |\psi^{(4,1)}_{\ast}|^2 & = & |\psi^{(1,4)}_{\ast}|^2 \;\; = \;\; (2+\sqrt{3})/36,
\end{eqnarray*}
whilst for the repeated eigenvalues, for the exponents with multiplicity two which we will label by $(\lambda_1,\lambda_2)_1$, $(\lambda_1,\lambda_2)_2$, we have
$$|\psi^{(2,2)_1}_{\ast}|^2 + |\psi^{(2,2)_2}_{\ast}|^2 \;\; = \;\; |\psi^{(5,2)_1}_{\ast}|^2 + |\psi^{(5,2)_2}_{\ast}|^2 \;\; = \;\; |\psi^{(2,5)_1}_{\ast}|^2 + |\psi^{(2,5)_2}_{\ast}|^2 \;\; = \;\; 2/9.$$

With $\theta_1 = (\lambda_1 + 2\lambda_2 + 3)/24$, $\theta_2 = (2\lambda_1 + \lambda_2 + 3)/24$, we have
\begin{center}
\renewcommand{\arraystretch}{1.5}
\begin{tabular}{|c|c|c|} \hline
$\lambda \in \mathrm{Exp}$ & $(\theta_1,\theta_2) \in [0,1]^2$ & $\frac{1}{16\pi^4} J(\theta_1,\theta_2)^2$ \\
\hline $(0,0)$, $(9,0)$, $(0,9)$ & $\left(\frac{1}{12},\frac{1}{12}\right)$, $\left(\frac{7}{12},\frac{1}{3}\right)$, $\left(\frac{1}{3},\frac{7}{12}\right)$ & $\frac{7-4\sqrt{3}}{4}$ \\
\hline $(4,4)$, $(4,1)$, $(1,4)$ & $\left(\frac{5}{12},\frac{5}{12}\right)$, $\left(\frac{1}{3},\frac{1}{4}\right)$, $\left(\frac{1}{4},\frac{1}{3}\right)$ & $\frac{7+4\sqrt{3}}{4}$ \\
\hline $(2,2)$, $(5,2)$, $(2,5)$ & $\left(\frac{1}{4},\frac{1}{4}\right)$, $\left(\frac{5}{12},\frac{1}{3}\right)$, $\left(\frac{1}{3},\frac{5}{12}\right)$ & 4 \\
\hline
\end{tabular}
\end{center}

Again, from (\ref{eqn:moments_SU(3)}),
\begin{equation} \label{eqn:moments-E1(12)}
\int_{\mathbb{T}^2} R_{m,n}(\omega_1,\omega_2) \mathrm{d}\varepsilon(\omega_1,\omega_2) = \frac{1}{6} \sum_{g \in S_3} \sum_{\lambda \in \mathrm{Exp}} (\beta^{(g(\lambda))})^m (\overline{\beta^{(g(\lambda))}})^n |\psi^{g(\lambda)}_{\ast}|^2.
\end{equation}

\begin{figure}[bt]
\begin{center}
  \includegraphics[width=55mm]{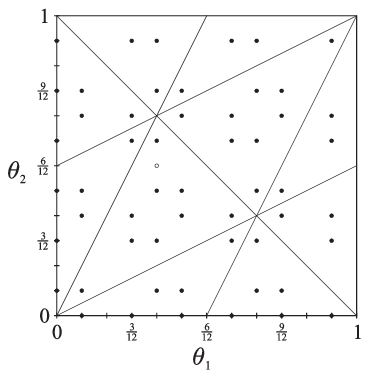}\\
 \caption{The points $(\theta_1,\theta_2) \in \{ g(\lambda) | \; \lambda \in \mathrm{Exp}, g \in S_3 \}$ for $\mathcal{E}_1^{(12)}$. The white circle indicates is the point $(5/12,6/12)$.} \label{fig:poly-12}
\end{center}
\end{figure}

We illustrate the pairs $(\theta_1,\theta_2)$ given by $g(\lambda)$ for $\lambda \in \mathrm{Exp}$, $g \in S_3$, in Figure \ref{fig:poly-12}.
Consider the pairs $(\theta_1,\theta_2) = (4/12,7/12), (3/12,5/12)$. For both of these, $(\omega_1,\omega_2) = (e^{2 \pi i \theta_1}, e^{2 \pi i \theta_2}) \in \mathbb{T}^2$ can only be obtained in the integral in (\ref{eqn:moments-E1(12)}) by using either the product measure $\mathrm{d}_{6} \times \mathrm{d}_{6}$ or the measure $\mathrm{d}^{(4)}$ ($(4/12,7/12)$, $(3/12,5/12)$ are both in $D_{4}$, but neither are in $D_{k}$ for any integer $k < 4$). With either of these measures, we will also obtain the point $(e^{2 \pi i 5/12}, e^{2 \pi i 6/12})$ in the integral (\ref{eqn:moments-E1(12)}). The corresponding pair $(\theta_1,\theta_2)$ is indicated by the white circle in Figure \ref{fig:poly-12}. The point $(e^{2 \pi i 5/12}, e^{2 \pi i 6/12})$ can also only obtained by using the measures $\mathrm{d}_{6} \times \mathrm{d}_{6}$ or $\mathrm{d}^{(4)}$.
Since these points $(\omega_1,\omega_2)$ cannot be obtained independently of each other, we must find a linear combination $\varepsilon' = c_1 \varepsilon_1 + c_2 J^2 \varepsilon_2$ of measures, where $\varepsilon_j$ must be either $\mathrm{d}_{6} \times \mathrm{d}_{6}$ or $\mathrm{d}^{(4)}$ for $j=1,2$, such that the weight $\varepsilon'(e^{2 \pi i \theta_1}, e^{2 \pi i \theta_2})$ is $(2-\sqrt{3})/36$, $(2+\sqrt{3})/36$, 0 for $(\theta_1,\theta_2) = (4/12,7/12), (3/12,5/12), (5/12,6/12)$ respectively. Suppose for now that $\varepsilon_1 = \varepsilon_2$ (again, it doesn't matter at this stage which of the two measures we take $\varepsilon_1$, $\varepsilon_2$ to be). Then since $J(5/12,6/12)^2 = 3/4$, we must find solutions $c_1,c_2 \in \mathbb{C}$ such that
\begin{equation} \label{eqns:forE1(12)}
c_1 + \frac{7-4\sqrt{3}}{4}c_2 = \frac{2-\sqrt{3}}{36}, \qquad c_1 + \frac{7+4\sqrt{3}}{4}c_2 = \frac{2+\sqrt{3}}{36}, \qquad c_1 + \frac{3}{4} c_2 = 0.
\end{equation}
However, no solution exists to the equations (\ref{eqns:forE1(12)}), and so the spectral measure for $\mathcal{E}_1^{(12)}$ is not a linear combination of measures of type $\mathrm{d}_{p/2} \times \mathrm{d}_{p/2}$, $J^2 \mathrm{d}_{p/2} \times \mathrm{d}_{p/2}$, $\mathrm{d}^{(p)}$ or $J^2 \mathrm{d}^{(p)}$, for $p \in \mathbb{N}$.

\section{Hilbert Series of $q$-deformations of CY-Algebras of Dimension 3} \label{sect:Hilbert-SU(3)}

We will now introduce the Calabi-Yau and $q$-deformed Calabi-Yau algebras of dimension 3, which are the $SU(3)$ generalizations of the pre-projective algebras of Section \ref{Sect:SU(2)-Hilbert_Series}. For certain $\mathcal{ADE}$ graphs we will also compute the Hilbert series of the $q$-deformed CY-algebras of dimension 3.

Let $\mathcal{G}$ be an oriented graph, and $\mathbb{C}\mathcal{G}$, $[\mathbb{C}\mathcal{G}, \mathbb{C}\mathcal{G}]$ be as in Section \ref{Sect:SU(2)-Hilbert_Series}. We define a derivation $\partial_a : \mathbb{C}\mathcal{G} / [\mathbb{C}\mathcal{G}, \mathbb{C}\mathcal{G}] \rightarrow \mathbb{C}\mathcal{G}$ by
$$\partial_{a} (a_1 \cdots a_n) = \sum_{j} a_{j+1} \cdots a_n a_1 \cdots a_{j-1},$$
where the summation is over all indices $j$ such that $a_j = a$.
Then for a potential $\Phi \in \mathbb{C}\mathcal{G} / [\mathbb{C}\mathcal{G}, \mathbb{C}\mathcal{G}]$, which is some linear combination of cyclic paths in $\mathcal{G}$, we define the algebra
$$A(\mathbb{C}\mathcal{G}, \Phi) \cong \mathbb{C}\mathcal{G} / \{ \partial_a \Phi \},$$
which is the quotient of the path algebra by the two-sided ideal generated by the elements $\partial_a \Phi \in \mathbb{C}\mathcal{G}$, for all edges $a$ of $\mathcal{G}$.
We define the Hilbert series $H_A(t)$ as in Section \ref{Sect:SU(2)-Hilbert_Series}.

If $A (\mathbb{C} \mathcal{G}, \Phi)$ is a Calabi-Yau algebra of dimensions $d \geq 3$ and $\textrm{deg} \; \Phi = d$, then \cite[Theorem 4.6]{bocklandt:2008}
\begin{equation} \label{eqn:H(t)-CYd}
H_A(t) = \frac{1}{1 - \Delta_{\mathcal{G}} t + \Delta_{\mathcal{G}}^T t^{d-1} - t^d}.
\end{equation}

Let $\Gamma$ be a subgroup of $SU(3)$. We do not concern ourselves here with the computation of the spectral measure of $\Gamma$, reserving that for a future publication \cite{evans/pugh:2009vi}. However, we make the following observation. Let $\Phi:\mathbb{T}^2 \rightarrow \mathfrak{D}$ be the map defined in (\ref{Phi:T^2->S}) and suppose we wish to compute `inverse' maps $\Phi^{-1}:\mathfrak{D} \rightarrow \mathbb{T}^2$ such that $\Phi \circ \Phi^{-1} = \mathrm{id}$, as we did for $SU(2)$ in (\ref{def:inverse_phi-SU(2)}).
For $z \in \mathfrak{D}$, we can write $z = \omega_1 + \omega_2^{-1} + \omega_1^{-1} \omega_2$ and $\overline{z} = \omega_1^{-1} + \omega_2 + \omega_1 \omega_2^{-1}$. Multiplying the first equation through by $\omega_1$, we obtain $z \omega_1 = \omega_1^2 + \omega_1 \omega_2^{-1} + \omega_2$. Then we need to find solutions $\omega_1$ to the cubic equation
\begin{equation} \label{eqn:cubic_w1}
\omega_1^3 - z \omega_1^2 + \overline{z} \omega_1 - 1 = 0.
\end{equation}
Similarly, we need to find solutions $\omega_2$ to the cubic equation $\omega_2^3 - \overline{z} \omega_2^2 + z \omega_2 - 1 = 0$, hence the three solutions for $\omega_2$ are given by the complex conjugate of the three solutions for $\omega_1$.
Solving (\ref{eqn:cubic_w1}) we obtain solutions $\omega^{(k)}$, $k = 0,1,2$, given by
$$\omega^{(k)} = (z + 2^{-1/3} \epsilon_k P + 2^{1/3} \overline{\epsilon_k} (z^2-3\overline{z}) P^{-1})/3,$$
where $\epsilon_k = e^{2 \pi i k /3}$, $2^{1/3}$ takes a real value, and $P$ is the cube root $P = (27 - 9z\overline{z} + 2z^3 + 3 \sqrt{3} \sqrt{27 - 18z\overline{z} + 4z^2 + 4\overline{z}^3 - z^2\overline{z}^2})^{1/3}$ such that $P \in \{ r e^{i \theta} | \; 0 \leq \theta < 2 \pi/3 \}$. For the roots of a cubic equation that it does not matter whether the square root in $P$ is taken to be positive or negative. We notice that the Jacobian $J$ appears in the expression for $P$ as the discriminant of the cubic equation (\ref{eqn:cubic_w1}).

We now consider the Hilbert series for $\Gamma$. For the McKay graph $\mathcal{G}_{\Gamma}$ one can define a cell system $W$ as in \cite{ginzburg:2006}, where $W(\triangle_{ijk})$ is a complex number for every triangle $\triangle_{ijk}$ on $\mathcal{G}_{\Gamma}$ whose vertices are labelled by the irreducible representations $i$, $j$, $k$ of $\Gamma$.
We introduce the following potential
$$\Phi_{\Gamma} = \sum_{\triangle_{ijk} \in \mathcal{G}_{\Gamma}} W(\triangle_{ijk}) \cdot \triangle_{ijk} \quad \in \mathbb{C} \mathcal{G}_{\Gamma} / [\mathbb{C} \mathcal{G}_{\Gamma}, \mathbb{C} \mathcal{G}_{\Gamma}].$$
Then dividing out $\mathbb{C} \mathcal{G}_{\Gamma}$ by the ideal generated by $\delta_a \Phi_{\Gamma}$ for all edges $a$ of $\mathcal{G}_{\Gamma}$, by \cite[Theorem 4.4.6]{ginzburg:2006}, $A(\mathbb{C}\mathcal{G}_{\Gamma}, \Phi_{\Gamma})$ is a Calabi-Yau algebra of dimension 3, and the Hilbert series is given by (\ref{eqn:H(t)-CYd}).

\begin{theorem}
Let $\Gamma$ be a finite subgroup of $SU(3)$,$\{ \rho_0 = \mathrm{id}, \rho_1 = \rho, \rho_2, \ldots, \rho_s \}$ its irreducible representations and $\mathcal{G}_{\Gamma}$ its McKay graph. Then if $P_{S, \rho_j}$ is the Molien series of the symmetric algebra $S$ of $\overline{\mathbb{C}^N}$, and $H(t)$ is the Hilbert series of $A(\mathbb{C}\mathcal{G}_{\Gamma}, \Phi_{\Gamma})$,
$$H_{\rho_j, 1_0} (t) = P_{S, \rho_j} (t).$$
\end{theorem}
\emph{Proof:}
Let $\Gamma$ be a subgroup of $SU(3)$ with irreducible representations $\rho_j$, $j=1, \ldots \; , s$, where $\rho_0 = \mathrm{id}$ is the identity representation and $\rho_1 = \rho$ the fundamental representation. The fundamental matrices $\Delta_{\Gamma}$, $\Delta_{\Gamma}^T$ defined by $\rho \otimes \rho_i = \sum_{j=0}^s (\Delta_{\Gamma})_{i,j} \rho_j$, $\overline{\rho} \otimes \rho_i = \sum_{j=0}^s (\Delta_{\Gamma}^T)_{i,j} \rho_j$, satisfy, by \cite[Cor. 2.4(i)]{gomi/nakamura/shinoda:2004},
$$\sum_{j=0}^s \left( -(\Delta_{\Gamma})_{\rho_i, \rho_j} t + (\Delta_{\Gamma}^T)_{\rho_i, \rho_j} t^2 \right) P_{S, \rho_j} (t) = - (1-t^3) P_{S, \rho_i} (t) + \delta_{i,0},$$
so we have
\begin{eqnarray*}
\sum_{j=0}^s \left(\textbf{1}_{\rho_i, \rho_j} -(\Delta_{\Gamma})_{\rho_i, \rho_j} t + (\Delta_{\Gamma}^T)_{\rho_i, \rho_j} t^2 - \textbf{1}_{\rho_i, \rho_j} t^3 \right) P_{S, \rho_j} (t) & = & \delta_{i,0} \\
\sum_{j=0}^s \left(\textbf{1} -(\Delta_{\Gamma}) t + (\Delta_{\Gamma}^T) t^2 - \textbf{1} t^3 \right)_{\rho_i, \rho_j} P_{S, \rho_j} (t) & = & \delta_{i,0}.
\end{eqnarray*}
Then $\left( P_{S,\rho_j} (t) \right)_{\rho_j}$ is given by the first column of the inverse of the invertible matrix $\left( \textbf{1} -(\Delta_{\Gamma}) t + (\Delta_{\Gamma}^T) t^2 - \textbf{1} t^3 \right)$, that is,
$$\hspace{30mm} P_{S, \rho_j} (t) = \left( \left( \textbf{1} -(\Delta_{\Gamma}) t + (\Delta_{\Gamma}^T) t^2 - \textbf{1} t^3 \right)^{-1} \right)_{\rho_j, \rho_0} = H_{\rho_j, \rho_0}. \hspace{25mm} \Box$$

For the $\mathcal{ADE}$ graphs, we define a potential $\Phi$ by
$$\Phi = \sum_{i,j,k} W(\triangle_{ijk}) \cdot \triangle_{ijk} \quad \in \mathbb{C} \mathcal{G} / [\mathbb{C} \mathcal{G}, \mathbb{C} \mathcal{G}],$$
where the Ocneanu cells $W(\triangle_{ijk})$ are computed in \cite{evans/pugh:2009i}.
The Hilbert series for the $q$-deformed $A(\mathbb{C}\mathcal{G}, \Phi)$ is given by

\begin{equation} \label{eqn:Hilbert_Series-SU(3)ADE}
H_{\mathcal{G}} (t) = \frac{1 - P t^h}{1 - \Delta_{\mathcal{G}} t + \Delta_{\mathcal{G}}^T t^2 - t^3},
\end{equation}
where $P$ is the permutation matrix corresponding to a $\mathbb{Z}/3\mathbb{Z}$ symmetry of the graph, and $h$ is the Coxeter number of $\mathcal{G}$.

The permutation matrix $P$ is an automorphism of the underlying graph, which is the identity for $\mathcal{D}^{(n)}$, $\mathcal{A}^{(n)\ast}$, $n \geq 5$, $\mathcal{E}^{(8)\ast}$, $\mathcal{E}_l^{(12)}$, $l=1,2,4,5$, and $\mathcal{E}^{(24)}$. For the remaining graphs, let $V$ be the permutation matrix corresponding to the clockwise rotation of the graph by $2 \pi /3$. Then
$$ P = \left\{
\begin{array}{cl} V & \mbox{ for } \quad \mathcal{A}^{(n)}, n \geq 4, \mbox{ and } \mathcal{E}^{(8)}, \\
                  V^{2n} & \mbox{ for } \quad \mathcal{D}^{(n) \ast}, n \geq 5.
\end{array} \right.$$

The numerator and denominator in (\ref{eqn:Hilbert_Series-SU(3)ADE}) commute. To see this note that $Q \Delta_{\mathcal{G}} = \Delta_{\mathcal{G}} Q$ and $Q \Delta_{\mathcal{G}}^T = \Delta_{\mathcal{G}}^T Q$, since $Q$ is a permutation matrix which corresponds to a symmetry of the graph $\mathcal{G}$.
The proof of (\ref{eqn:Hilbert_Series-SU(3)ADE}) will appear in \cite{evans/pugh:2009vi}.

In the $SU(2)$ case, the permutation matrices $P$ appearing in the numerator of $H_A(t)$ corresponded to the Nakayama permutation of the Dynkin diagram. The above claim then raises the question of the relation between the automorphisms which appear in the numerators of the expressions for $H_A(t)$ with Nakayama's automorphisms.

\paragraph{Acknowledgements}

This paper is based on work in \cite{pugh:2008}. The first author was partially supported by the EU-NCG network in Non-Commutative Geometry MRTN-CT-2006-031962, and the second author was supported by a scholarship from the School of Mathematics, Cardiff University.


\begin{thebibliography}{99}

\bibitem{banica/bisch:2007}
T. Banica and D. Bisch, Spectral measures of small index principal graphs, Comm. Math. Phys. {\bf 269} (2007), 259--281.

\bibitem{banica:2007}
T. Banica, Cyclotomic expansion of exceptional spectral measures, 2007. arXiv:0712.2524v2 [math.QA].

\bibitem{behrend/pearce/petkova/zuber:2000}
R. E. Behrend, P. A. Pearce, V. B. Petkova and J.-B. Zuber, Boundary conditions in rational conformal field theories, Nuclear Phys. B {\bf 579} (2000), 707--773.

\bibitem{bion-nadal:1991}
J. Bion-Nadal, An example of a subfactor of the hyperfinite {${\rm II}\sb 1$} factor whose principal graph invariant is the {C}oxeter graph {$E\sb 6$}, in {\it Current topics in operator algebras ({N}ara, 1990)}, 104--113, World Sci. Publ., River Edge, NJ, 1991.

\bibitem{bockenhauer/evans:1998}
J. B{\"o}ckenhauer and D. E. Evans, Modular invariants, graphs and {$\alpha$}-induction for nets of subfactors, Comm. Math. Phys., {I} {\bf 197} (1998), 361--386. {II} {\bf 200} (1999), 57--103. {III} {\bf 205} (1999), 183--228.

\bibitem{bockenhauer/evans:2000}
J. B{\"o}ckenhauer and D. E. Evans, Modular invariants from subfactors: {T}ype {I} coupling matrices and intermediate subfactors, Comm. Math. Phys. {\bf 213} (2000), 267--289.

\bibitem{bockenhauer/evans/kawahigashi:1999}
J. B{\"o}ckenhauer, D. E. Evans and Y. Kawahigashi, On {$\alpha$}-induction, chiral generators and modular invariants for subfactors, Comm. Math. Phys. {\bf 208} (1999), 429--487.

\bibitem{bockenhauer/evans/kawahigashi:2000}
J. B{\"o}ckenhauer, D. E. Evans and Y. Kawahigashi, Chiral structure of modular invariants for subfactors, Comm. Math. Phys. {\bf 210} (2000), 733--784.

\bibitem{bocklandt:2008}
R. Bocklandt, Graded {C}alabi {Y}au algebras of dimension 3, J. Pure Appl. Algebra {\bf 212} (2008), 14--32.

\bibitem{brenner/butler/king:2002}
S. Brenner, M. C. R. Butler and A. D. King, Periodic algebras which are almost {K}oszul, Algebr. Represent. Theory {\bf 5} (2002), 331--367.

\bibitem{cappelli/itzykson/zuber:1987ii}
A. Cappelli, C. Itzykson, C. and J.-B. Zuber, The {${\rm A}$}-{${\rm D}$}-{${\rm E}$} classification of minimal and {$A\sp {(1)}\sb 1$} conformal invariant theories, Comm. Math. Phys. {\bf 113} (1987), 1--26.

\bibitem{di_francesco:1992}
P. Di Francesco, Integrable lattice models, graphs and modular invariant conformal field theories, Internat. J. Modern Phys. A {\bf 7} (1992), 407--500.

\bibitem{di_francesco:1997}
P. Di Francesco, {${\rm SU}(N)$} meander determinants, J. Math. Phys. {\bf 38} (1997), 5905--5943.

\bibitem{di_francesco/mathieu/senechal:1997}
P. Di Francesco, P. Mathieu and D. S{\'e}n{\'e}chal, {\it Conformal field theory}, Graduate Texts in Contemporary Physics, Springer-Verlag, New York, 1997.

\bibitem{di_francesco/zuber:1990}
P. Di Francesco and J.-B. Zuber, {${\rm SU}(N)$} lattice integrable models associated with graphs, Nuclear Phys. B {\bf 338} (1990), 602--646.

\bibitem{king/egecioglu:1999}
{\"O}. E{\u{g}}ecio{\u{g}}lu and A. King, Random walks and {C}atalan factorization, in {\it Proceedings of the {T}hirtieth {S}outheastern {I}nternational {C}onference on {C}ombinatorics, {G}raph {T}heory, and {C}omputing ({B}oca {R}aton, {FL}, 1999)}, Congr. Numer. {\bf 138}, 129--140, 1999.

\bibitem{erdmann/snashall:1998}
K. Erdmannand N. Snashall, Preprojective algebras of {D}ynkin type, periodicity and the second {H}ochschild cohomology, in {\it Algebras and modules, {II} ({G}eiranger, 1996)}, CMS Conf. Proc. {\bf 24}, 183--193, Amer. Math. Soc., Providence, RI, 1998.

\bibitem{etingof/ostrik:2004}
P. Etingof and V. Ostrik, Module categories over representations of {${\rm SL}\sb q(2)$} and graphs, Math. Res. Lett. {\bf 11} (2004), 103--114.

\bibitem{evans:2002}
D. E. Evans, Fusion rules of modular invariants, Rev. Math. Phys. {\bf 14} (2002), 709--731.

\bibitem{evans:2003}
D. E. Evans, Critical phenomena, modular invariants and operator algebras, in {\it Operator algebras and mathematical physics (Constan\c ta, 2001)}, 89--113, Theta, Bucharest, 2003.

\bibitem{evans/kawahigashi:1998}
D. E. Evans and Y. Kawahigashi, {\it Quantum symmetries on operator algebras}, Oxford Mathematical Monographs. The Clarendon Press Oxford University Press, New York, 1998. Oxford Science Publications.

\bibitem{evans/pugh:2009i}
D. E. Evans and M. Pugh, Ocneanu Cells and {B}oltzmann Weights for the $SU(3)$ $\mathcal{ADE}$ Graphs, M\"{u}nster J. Math. (to appear). arxiv:0906.4307

\bibitem{evans/pugh:2009ii}
D. E. Evans and M. Pugh, {$SU(3)$}-{G}oodman-de la {H}arpe-{J}ones subfactors and the realisation of {$SU(3)$} modular invariants, Rev. Math. Phys. (to appear). arxiv:0906.4252 (math.OA)

\bibitem{evans/pugh:2009iii}
D. E. Evans and M. Pugh, {$A_2$}-Planar Algebras {I}. Preprint, arxiv:0906.4225 (math.OA).

\bibitem{evans/pugh:2009iv}
D. E. Evans and M. Pugh, {$A_2$}-Planar Algebras {II}: Planar Modules. Preprint, arxiv:0906.4311 (math.OA).

\bibitem{evans/pugh:2009vi}
D. E. Evans and M. Pugh, Spectral Measures and Generating Series for Nimrep Graphs in Subfactor Theory {II}. In preparation.

\bibitem{farsi/watling:1993}
C. Farsi and N. Watling, Cubic algebras, J. Operator Theory {\bf 30} (1993), 243--266.

\bibitem{gaberdiel/gannon:2002}
M. R. Gaberdiel and T. Gannon, Boundary states for {WZW} models, Nuclear Phys. B {\bf 639} (2002), 471--501.

\bibitem{gannon:1994}
T. Gannon, The classification of affine {${\rm SU}(3)$} modular invariant partition functions, Comm. Math. Phys. {\bf 161} (1994), 233--263.

\bibitem{gepner:1991}
D. Gepner, Fusion rings and geometry, Comm. Math. Phys. {\bf 141} (1991), 381--411.

\bibitem{ginzburg:2006}
V. Ginzburg, Calabi-{Y}au algebras, 2006. arXiv:math/0612139 [math.AG].

\bibitem{gomi/nakamura/shinoda:2004}
Y. Gomi, I. Nakamura and K.-I. Shinoda, Coinvariant algebras of finite subgroups of {${\rm SL}(3,{\bf C})$}, Canad. J. Math. {\bf 56} (2004), 495--528.

\bibitem{goodman/de_la_harpe/jones:1989}
F. M. Goodman, P. de la Harpe and V. F. R. Jones, {\it Coxeter graphs and towers of algebras}, MSRI Publications, 14, Springer-Verlag, New York, 1989.

\bibitem{hiai/petz:2000}
F. Hiai and D. Petz, {\it The semicircle law, free random variables and entropy}, Mathematical Surveys and Monographs {\bf 77}, Amer. Math. Soc., Providence, RI, 2000.

\bibitem{itzykson:1989}
C. Itzykson, From the harmonic oscillator to the {$A$}-{$D$}-{$E$} classification of conformal models, in {\it Integrable systems in quantum field theory and statistical mechanics}, Adv. Stud. Pure Math. {\bf 19}, 287--346, Academic Press, Boston, MA, 1989.

\bibitem{izumi:1991}
M. Izumi, Application of fusion rules to classification of subfactors, Publ. Res. Inst. Math. Sci. {\bf 27} (1991), 953--994.

\bibitem{izumi:1994}
M. Izumi, On flatness of the {C}oxeter graph {$E\sb 8$}, Pacific J. Math. {\bf 166} (1994), 305--327.

\bibitem{jones:1983}
V. F. R. Jones, Index for subfactors, Invent. Math. {\bf 72} (1983), 1--25.

\bibitem{jones:planar}
V. F. R. Jones, Planar algebras. {I}, New Zealand J. Math. (to appear).

\bibitem{jones:2000}
V. F. R. Jones, The planar algebra of a bipartite graph, in {\it Knots in Hellas '98 (Delphi)}, Ser. Knots Everything {\bf 24}, 94--117, World Sci. Publ., River Edge, NJ, 2000.

\bibitem{jones:2001}
V. F. R. Jones, The annular structure of subfactors, in {\it Essays on geometry and related topics, {V}ol. 1, 2}, Monogr. Enseign. Math. {\bf 38}, 401--463, Enseignement Math., Geneva, 2001.

\bibitem{kassel:1995}
C. Kassel, {\it Quantum groups}, Graduate Texts in Mathematics {\bf 155}, Springer-Verlag, New York, 1995.

\bibitem{kawahigashi:1995}
Y. Kawahigashi, On flatness of {O}cneanu's connections on the {D}ynkin diagrams and classification of subfactors, J. Funct. Anal. {\bf 127} (1995), 63--107.

\bibitem{kawai:1989}
T. Kawai, On the structure of fusion algebras, Phys. Lett. B {\bf 217} (1989), 47--251.

\bibitem{kostant:1984}
B. Kostant, On finite subgroups of {${\rm SU}(2)$}, simple {L}ie algebras, and the {M}c{K}ay correspondence, Proc. Nat. Acad. Sci. U.S.A. {\bf 81} (1984), 5275--5277.

\bibitem{malkin/ostrik/vybornov:2006}
A. Malkin, V. Ostrik and M. Vybornov, Quiver varieties and {L}usztig's algebra, Adv. Math. {\bf 203} (2006), 514--536.

\bibitem{mckay:1980}
J. McKay, Graphs, singularities, and finite groups, in {\it The {S}anta {C}ruz {C}onference on {F}inite {G}roups ({U}niv. {C}alifornia, {S}anta {C}ruz, {C}alif., 1979)}, Proc. Sympos. Pure Math. {\bf 37}, 183--186, Amer. Math. Soc., Providence, R.I., 1980.

\bibitem{ocneanu:1988}
A. Ocneanu, Quantized groups, string algebras and {G}alois theory for algebras, in {\it Operator algebras and applications, {V}ol.\ 2}, London Math. Soc. Lecture Note Ser. {\bf 136}, 119--172, Cambridge Univ. Press, Cambridge, 1988.

\bibitem{ocneanu:2000i}
A. Ocneanu, Paths on Coxeter diagrams: from Platonic solids and singularities to minimal models and subfactors. (Notes recorded by S. Goto), in
{\it Lectures on operator theory}, (ed. B. V. Rajarama Bhat et al.), The Fields Institute Monographs, 243--323, Amer. Math. Soc., Providence, R.I., 2000.

\bibitem{ocneanu:2000ii}
A. Ocneanu, Higher {C}oxeter Systems (2000). Talk given at MSRI. \\ http://www.msri.org/publications/ln/msri/2000/subfactors/ocneanu.

\bibitem{ocneanu:2002}
A. Ocneanu, The classification of subgroups of quantum {${\rm SU}(N)$}, in {\it Quantum symmetries in theoretical physics and mathematics (Bariloche, 2000)}, Contemp. Math. {\bf 294}, 133--159, Amer. Math. Soc., Providence, RI, 2002.

\bibitem{pugh:2008}
M. Pugh, The Ising Model and Beyond, PhD thesis, Cardiff University, 2008.

\bibitem{reid:2002}
M. Reid, La correspondance de {M}c{K}ay, Ast\'erisque 276 (2002), 53--72. S{\'e}minaire Bourbaki, Vol. 1999/2000.

\bibitem{sagan:2001}
B. E. Sagan, {\it The symmetric group}, Graduate Texts in Mathematics {\bf 203}, Springer-Verlag, New York, 2001. Representations, combinatorial algorithms, and symmetric functions.

\bibitem{dykema/voiculescu/nica:1992}
D. V. Voiculescu, K. J. Dykema and A. Nica, {\it Free random variables}, CRM Monograph Series {\bf 1}, American Mathematical Society, Providence, RI, 1992.

\bibitem{wassermann:1998}
A. Wassermann, Operator algebras and conformal field theory. {III}. {F}usion of positive energy representations of {${\rm LSU}(N)$} using bounded operators, Invent. Math. {\bf 133} (1998), 467--538.

\bibitem{xu:1998}
F. Xu, New braided endomorphisms from conformal inclusions, Comm. Math. Phys. {\bf 192} (1998), 349--403.

\end{thebibliography}
\end{document}